\documentclass[12pt]{article}
\usepackage{amssymb, amsmath, amsthm, amscd}
\usepackage[utf8]{inputenc}
\usepackage[all,cmtip]{xy}
\usepackage{enumitem}
\usepackage{setspace}
\usepackage{mathdots}

\usepackage[mathscr]{euscript}
\usepackage[colorlinks=true]{hyperref}
\hypersetup{colorlinks   = true}
\hypersetup{linkcolor=blue}
\usepackage{tikz}
\begin{document}

\mathchardef\mhyphen="2D
\newtheorem{The}{Theorem}[section]
\newtheorem{Lem}[The]{Lemma}
\newtheorem{Prop}[The]{Proposition}
\newtheorem{Cor}[The]{Corollary}
\newtheorem{Rem}[The]{Remark}
\newtheorem{Obs}[The]{Observation}
\newtheorem{SConj}[The]{Standard Conjecture}
\newtheorem{Titre}[The]{\!\!\!\! }
\newtheorem{Conj}[The]{Conjecture}
\newtheorem{Question}[The]{Question}
\newtheorem{Prob}[The]{Problem}
\newtheorem{Def}[The]{Definition}
\newtheorem{Not}[The]{Notation}
\newtheorem{Claim}[The]{Claim}
\newtheorem{Conc}[The]{Conclusion}
\newtheorem{Ex}[The]{Example}
\newtheorem{Fact}[The]{Fact}
\newtheorem{Formula}[The]{Formula}
\newtheorem{Formulae}[The]{Formulae}
\newtheorem{The-Def}[The]{Theorem and Definition}
\newtheorem{Prop-Def}[The]{Proposition and Definition}
\newtheorem{Lem-Def}[The]{Lemma and Definition}
\newtheorem{Cor-Def}[The]{Corollary and Definition}
\newtheorem{Conc-Def}[The]{Conclusion and Definition}
\newtheorem{Gen-Fact}[The]{General Fact}
\newtheorem{Terminology}[The]{Note on terminology}
\newcommand{\C}{\mathbb{C}}
\newcommand{\R}{\mathbb{R}}
\newcommand{\N}{\mathbb{N}}
\newcommand{\Z}{\mathbb{Z}}
\newcommand{\Q}{\mathbb{Q}}
\newcommand{\Proj}{\mathbb{P}}
\newcommand{\Rc}{\mathcal{R}}
\newcommand{\Oc}{\mathcal{O}}
\newcommand{\Vc}{\mathcal{V}}
\newcommand{\Id}{\operatorname{Id}}
\newcommand{\pr}{\operatorname{pr}}
\newcommand{\rk}{\operatorname{rk}}
\newcommand{\del}{\partial}
\newcommand{\delbar}{\bar{\partial}}
\newcommand{\Cdot}{{\raisebox{-0.7ex}[0pt][0pt]{\scalebox{2.0}{$\cdot$}}}}
\newcommand\nilm{\Gamma\backslash G}
\newcommand\frg{{\mathfrak g}}
\newcommand{\fg}{\mathfrak g}
\newcommand{\Oh}{\mathcal{O}}
\newcommand{\Kur}{\operatorname{Kur}}
\newcommand\gc{\frg_\mathbb{C}}
\newcommand\hisashi[1]{{\textcolor{red}{#1}}}
\newcommand\dan[1]{{\textcolor{blue}{#1}}}
\newcommand\luis[1]{{\textcolor{orange}{#1}}}

\begin{center}

{\Large\bf Adiabatic Limit and Deformations of Complex Structures}

\end{center}

\begin{center}

{\large Dan Popovici}

\end{center}

\vspace{1ex}

\noindent{\small{\bf Abstract.} Based on our recent adaptation of the adiabatic limit construction to the case of complex structures, we prove the fact that the deformation limiting manifold of any holomorphic family of Moishezon manifolds is Moishezon. Two new ingredients, hopefully of independent interest, are introduced. The first one associates with every compact complex manifold $X$, in every degree $k$, a holomorphic vector bundle over $\C$ of rank equal to the $k$-th Betti number of $X$. This vector bundle, previously given an algebraic construction in the literature, shows that the degenerating page of the Fr\"olicher spectral sequence of $X$ is the holomorphic limit, as $h\in\C^\star$ tends to $0$, of the $d_h$-cohomology of $X$, where $d_h=h\partial + \bar\partial$. A relative version of this vector bundle is then associated with every holomorphic family of compact complex manifolds. The second ingredient is a relaxation of the notion of strongly Gauduchon (sG) metric that we introduced in 2009. For a given positive integer $r$, a Gauduchon metric $\gamma$ on an $n$-dimensional compact complex manifold $X$ is said to be $E_r$-sG if $\partial\gamma^{n-1}$ represents the zero cohomology class on the $r$-th page of the Fr\"olicher spectral sequence of $X$. Strongly Gauduchon metrics coincide with $E_1$-sG metrics.}

\vspace{2ex}

\section{Introduction}\label{section:Introduction}

The main result of this paper is the following statement that first appeared in [Pop09] and [Pop10] where it was given a different, ad hoc treatment, although the general strategy and some details were the same as in the present, more conceptual, approach.

\begin{The}\label{The:moidef} Let $N$ be a positive integer and let $\pi:{\mathcal X}\rightarrow B$ be a complex analytic family of compact complex manifolds over an open ball $B\subset\C^N$ about the origin such that the fibre $X_t:=\pi^{-1}(t)$ is a Moishezon manifold for every $t\in B\setminus\{0\}$. Then $X_0:=\pi^{-1}(0)$ is again a Moishezon manifold.

\end{The}

As usual, by a complex analytic (or holomorphic) family of compact complex manifolds we mean a {\it proper holomorphic submersion} $\pi:{\mathcal X}\rightarrow B$ between two complex manifolds ${\mathcal X}$ and $B$ (cf. e.g. [Kod86]). In particular, the fibres $X_t:=\pi^{-1}(t)$ are compact complex manifolds of the same dimension. By a classical theorem of Ehresmann [Ehr47], any such family is locally (hence also globally if the base $B$ is contractible) $C^\infty$ trivial. Thus, all the fibres $X_t$ have the same underlying $C^\infty$ manifold $X$ (hence also the same De Rham cohomology groups $H^k_{DR}(X,\,\C)$ for all $k=0,\dots , 2n$), but the complex structure $J_t$ of $X_t$ depends, in general, on $t\in B$.

On the other hand, a {\it Moishezon manifold} is a compact complex manifold $Y$ for which there exists a projective manifold $\widetilde{Y}$ and a holomorphic bimeromorphic map $\mu:\widetilde{Y}\longrightarrow Y$ (cf. [Moi67]). By another classical result of [Moi67], we know that a Moishezon manifold is not K\"ahler unless it is projective.

Our Theorem \ref{The:moidef} above is a closedness result under deformations of complex structures: any deformation limit of a family of Moishezon manifolds is Moishezon. Indeed, the fibre $X_0$ can be regarded as the limit of the fibres $X_t$ when $t\in B$ tends to $0\in B$. We can suppose, without loss of generality, that $B$ is an open disc about the origin in $\C$.

\subsection{Machinery used in the paper}\label{subsection:FSS_introd} We will cast much of the proof of Theorem \ref{The:moidef} in the language of the  Fr\"olicher spectral sequence (FSS) of a given compact complex manifold $X$ with $\mbox{dim}_\C X=n$. The FSS is a finite family of complexes, called {\it pages}, that successively refine the Dolbeault cohomology of $X$ until it ``degenerates'' to the De Rham cohomology.  

Recall that the zero-th page $E_0$ of the FSS consists of the $\C$-vector spaces $E_0^{p,\,q}(X):=C^\infty_{p,\,q}(X,\,\C)$ of smooth pure-type forms on $X$ and of the type-$(0,\,1)$ differentials $d_0:=\bar\partial$ forming the Dolbeault complex: $$\dots\stackrel{d_0}{\longrightarrow}E_0^{p,\,q-1}(X)\stackrel{d_0}{\longrightarrow}E_0^{p,\,q}(X)\stackrel{d_0}{\longrightarrow}E_0^{p,\,q+1}(X)\stackrel{d_0}{\longrightarrow}\dots.$$ 

The first page $E_1$ of the Fr\"olicher spectral sequence consists of the $\C$-vector spaces $E_1^{p,\,q}(X)$ (i.e. the cohomology of the zero-th page) and of the type-$(1,\,0)$ differentials $d_1$: $$\dots\stackrel{d_1}{\longrightarrow}E_1^{p-1,\,q}(X)\stackrel{d_1}{\longrightarrow}E_1^{p,\,q}(X)\stackrel{d_1}{\longrightarrow}E_1^{p+1,\,q}(X)\stackrel{d_1}{\longrightarrow}\dots.$$

\noindent induced in cohomology by $\partial$ (i.e. $d_1([\alpha]_{\bar\partial}):=[\partial\alpha]_{\bar\partial}$).

The remaining pages are constructed inductively: the differentials $d_r=d_r^{p,\,q}:E_r^{p,\,q}(X)\longrightarrow E_r^{p+r,\,q-r+1}(X)$ are of type $(r,\,-r+1)$ for every $r$, while the $\C$-vector spaces $E_r^{p,\,q}(X):=\ker d_{r-1}^{p,\,q}/\mbox{Im}\,d_{r-1}^{p-r+1,\,q+r-2}$ on the $r^{th}$ page are defined as the cohomology of the previous page $E_{r-1}$.

A classical result of Fr\"olicher [Fro55] asserts that this spectral sequence converges to the De Rham cohomology of $X$ and degenerates after finitely many steps. This means that there are (not necessarily canonical) isomorphisms: \begin{equation}\label{eqn:Froelicher_decomposition}H^k_{DR}(X,\,\C)\simeq\bigoplus\limits_{p+q=k}E_\infty^{p,\,q}(X), \hspace{6ex} k=0,\dots , 2n,\end{equation}

\noindent where $E_\infty^{p,\,q}(X) = \dots = E_{r+2}^{p,\,q}(X) = E_{r+1}^{p,\,q}(X)= E_r^{p,\,q}(X)$ for all $p,q$ and $r\geq 1$ is the smallest positive integer such that the spectral sequence degenerates at $E_r$.

The degeneration at the $r$-th page is denoted by $E_r(X) = E_\infty(X)$.

\subsection{Brief reminder of the main construction in [Pop17]}\label{subsection:reminder} The method introduced in this paper originates in our recent adaptation to the case of complex structures (cf. [Pop17]) of the adiabatic limit construction associated with Riemannian foliations (cf., e.g., [Wit85] and [MM90]). Given a compact complex $n$-dimensional manifold $X$, for every constant $h\in \C$, we associate with the splitting $d=\partial + \bar\partial$ defining the complex structure of $X$ the following $1st$-order differential operator: $$d_h:=h\partial + \bar\partial : C^\infty_k(X,\,\C)\longrightarrow C^\infty_{k+1}(X,\,\C), \hspace{3ex} k=0,\dots , 2n,$$

\noindent acting on the space $C^\infty_k(X,\,\C)$ of smooth $k$-forms on $X$, for every degree $k$. Only positive real constants $h$ were considered in [Pop17], but we now allow $h$ to be any {\it complex} constant. In particular, $d_h$ depends on the complex structure of $X$, except when $h=1$, in which case $d_1=d$. On the other hand, $d_0=\bar\partial$.

Meanwhile, for every non-zero $h$, the linear map defined pointwise on $k$-forms by $$\theta_h:\Lambda^kT^\star X\longrightarrow\Lambda^kT^\star X,  \hspace{3ex} u=\sum\limits_{p+q=k}u^{p,\,q}\longmapsto\theta_hu:=\sum\limits_{p+q=k}h^p\,u^{p,\,q},$$

\noindent induces an automorphism of the vector bundle $\Lambda T^\star X = \oplus_{k=0}^{2n}\Lambda^kT^\star X$ and the operators $d_h$ and $d$ are related by the identity $$d_h = \theta_h d\theta_h^{-1}.$$

 This implies that $d_h^2=0$, so we can define the {\it $d_h$-cohomology} of $X$ (cf. [Pop17]) in every degree $k$ as $$H^k_{d_h}(X,\,\C):=\ker\bigg(d_h: C^\infty_k(X,\,\C)\longrightarrow C^\infty_{k+1}(X,\,\C)\bigg)\bigg/\mbox{Im}\,\bigg(d_h: C^\infty_{k-1}(X,\,\C)\longrightarrow C^\infty_k(X,\,\C)\bigg).$$

 Moreover, $\theta_h$ maps $d$-closed forms to $d_h$-closed forms and $d$-exact forms to $d_h$-exact forms, so it induces an isomorphism between the De Rham cohomology and the $d_h$-cohomology for every $h\in\C\setminus\{0\}$: $$\theta_h:H^k_{DR}(X,\,\C)\stackrel{\simeq}{\longrightarrow}H^k_{d_h}(X,\,\C), \hspace{3ex} k=0,\dots , 2n.$$

Now, if $X$ is given a Hermitian metric $\omega$, we let $d_h^\star$ be the formal adjoint of $d_h$ w.r.t. the $L^2$-inner product on differential forms induced by $\omega$. The $d_h$-Laplacian w.r.t. $\omega$ is defined in every degree $k$ in the expected way: $$\Delta_h:C^\infty_k(X,\,\C)\longrightarrow C^\infty_k(X,\,\C), \hspace{3ex} \Delta_h:=d_hd_h^\star + d_h^\star d_h.$$

 It turns out that the (non-negative, self-adjoint) $2nd$-order differential operator $\Delta_h$ is {\it elliptic} (cf. [Pop17, Lemma 2.7]). Together with the integrability of $d_h$ (i.e. $d_h^2=0$) and the compactness of $X$, this yields the Hodge isomorphism $$\ker\bigg(\Delta_h:C^\infty_k(X,\,\C)\to C^\infty_k(X,\,\C)\bigg)\simeq H^k_{d_h}(X,\,\C), \hspace{3ex} k=0,\dots , 2n,$$ for the $d_h$-cohomology. By elliptic theory, $ \Delta_h$ has a discrete spectrum $0\leq\lambda_1^{(k)}(h)\leq\lambda_2^{(k)}(h)\leq\dots\leq\lambda_j^{(k)}(h)\leq\dots$ whose only accumulation point is $\infty$. Moreover, for every $h\neq 0$, the isomorphism between the $d_h$- and the De Rham cohomologies implies that the multiplicity of $0$ as an eigenvalue of $\Delta_h$ acting on $k$-forms is the $k$-th Betti number $b_k$ of $X$.

\subsection{Constructions introduced in this paper}\label{subsection:new-ideas} There are two main ingredients that we hope are of independent interest and that we now outline. The details will occupy $\S$\ref{Froelicher_h}, respectively $\S$\ref{section:E_r-sG}.

\vspace{1ex}

\subsubsection{{\bf The Fr\"olicher approximating vector bundle (FAVB)}}\label{subsubsection:FAVB_introd} The first construction (cf. $\S$\ref{subsubsection:absolute} and $\S$\ref{subsubsection:relative}) builds on the adiabatic limit for complex structures introduced in [Pop17] and outlined above to prove that the degenerating page of the Fr\"olicher spectral sequence is the holomorphic limit, as $h\in\C$ tends to $0$, of the $d_h$-cohomology in every degree $k$.

\vspace{1ex}

{\bf (I)\, The absolute version of the FAVB}

\vspace{1ex}

Specifically, with every compact complex $n$-dimensional manifold $X$ and every degree $k\in\{0,\dots , 2n\}$, we associate a holomorphic vector bundle ${\mathcal A}^k$ of rank $b_k$ ($=$ the $k$-th Betti number of $X$) over $\C$ whose fibres are $${\mathcal A}^k_h:=H^k_{d_h}(X,\,\C) \hspace{3ex} \mbox{if} \hspace{1ex} h\in\C\setminus\{0\}, \hspace{6ex} {\mathcal A}^k_0:=\bigoplus\limits_{p+q=k}E_r^{p,\,q}(X) \hspace{3ex} \mbox{if} \hspace{1ex} h=0,$$

\noindent where $r\geq 1$ is the smallest positive integer such that the Fr\"olicher spectral sequence of $X$ degenerates at $E_r$. In other words, $E_r^{p,\,q}(X) = E_\infty^{p,\,q}(X)$ for all $p,q$.

Like any vector bundle over $\C$, ${\mathcal A}^k$ is {\it trivial}, but its main feature for us is that the canonical trivialisation over $\C\setminus\{0\}$ induced by the isomorphisms $$\theta_h:H^k_{DR}(X,\,\C)\to H^k_{d_h}(X,\,\C) \hspace{3ex} \mbox{with} \hspace{2ex} h\neq 0$$ (seen collectively as an isomorphism from the constant bundle ${\mathcal H}^k\to\C\setminus\{0\}$ of fibre $H^k_{DR}(X,\,\C)$ to ${\mathcal A}^k_{|\C\setminus\{0\}}$) {\it extends holomorphically over $0$} to a global holomorphic trivialisation over $\C$. 

In other words, $\oplus_{p+q=k}E_\infty^{p,\,q}(X)$ is the holomorphic limit of the family $(H^k_{d_h}(X,\,\C))_{h\in\C^\star}$ of vector spaces when $h$ tends to $0$.

 This fact is asserted by Corollary and Definition \ref{Cor-Def:bundle_fixed-manifold}. We call ${\mathcal A}^k$ the {\bf Fr\"olicher approximating vector bundle (FAVB)} of $X$ in degree $k$.

%\begin{The}\label{The:bundle_fixed-manifold_introd} For every $k\in\{0,\dots , 2n\}$, ${\mathcal A}^k\to\C$ is a holomorphic vector bundle of rank $b_k$.

%\end{The}

 This vector bundle has essentially been known for quite a while, although its construction has not been cast in the language of the $d_h$-cohomology employed here. The standard argument asserts the existence of an equivalence of categories between filtered complex vector spaces and $\C^{\star}$-equivariant algebraic (or holomorphic) vector bundles on the complex line $\C$. (See e.g. the discussion in [Aso06, $\S3$], which attributes the construction to [Ger66] and [Kly89].) In response to the first version of this paper, this vector bundle was revisited in [Ste20], where an algebraic construction of the FAVB was obtained as a special case of the Rees bundle construction of [Kly89]. In fact, [Ste20] shows that the FAVB coincides with the Rees bundle of the De Rham cohomology with its Hodge filtration. In particular, the FAVB is given a functorial description in [Ste20, Theorem B].

 Our treatment of the FAVB construction is analytic. It relies crucially on classical results by Kodaira and Spencer [KS60], mainly on the fundamental fact that a $C^\infty$ family of elliptic operators having kernels of the same dimension induces a $C^\infty$ vector bundle whose fibres are these kernels. Actually, the orthogonal projections onto these kernels vary in a $C^\infty$ way with the parameter. The results from [KS60] that we need are briefly reviewed in $\S$\ref{subsection:review_K-S} for the reader's convenience. 

 In our case, after fixing a Hermitian metric $\omega$ on $X$, we construct elliptic (pseudo-)differential operators depending in a $C^\infty$ way on a parameter $h\in\C$ whose kernels are equidimensional and isomorphic to the $d_h$-cohomology group $${\mathcal A}^k_h:=H^k_{d_h}(X,\,\C) \hspace{3ex} \mbox{for every} \hspace{1ex} h\in\C^\star$$ and to $${\mathcal A}^k_0:=\oplus_{p+q=k}E_r^{p,\,q}(X) \hspace{3ex} \mbox{when} \hspace{1ex} h=0.$$ 

The choice of the $C^\infty$ family of elliptic operators that defines the $C^\infty$ vector bundle structure of the FAVB ${\mathcal A}^k\longrightarrow\C$ depends on the value of $r\in \N^\star$. (See Conclusions \ref{Conc:pseudodiff_2} and \ref{Conc:pseudodiff_r}.) It is denoted by 

\vspace{1ex}

$\bullet$\, $(\Delta_h)_{h\in\C}$ (cf. $\S$\ref{subsubsection:E_1_Delta-tilde_h_def}) when $r=1$;

\vspace{1ex}

$\bullet$\, by $(\widetilde\Delta_h)_{h\in\C}$ (cf. $\S$\ref{subsubsection:E_2_Delta-tilde_h_def}) when $r=2$;

\vspace{1ex}

$\bullet$\, by $(\widetilde\Delta_h^{(r)})_{h\in\C}$ (cf. $\S$\ref{subsubsection:E_geq3_Delta-tilde_h-r_def}) when $r\geq 3$. 

\vspace{1ex}

\noindent To unify the notation, we will also let $\widetilde\Delta_h^{(1)}:=\Delta_h$ and $\widetilde\Delta_h^{(2)}:=\widetilde\Delta_h$ for every $h\in\C$.

 Since simple algebraic proofs of the FAVB construction exist in the literature for every $r$ (as mentioned above), we omit the technically involved details of the case $r\geq 3$ of this construction arising in the $d_h$-cohomology treatment given here. These details are available in full in earlier arXiv versions of this paper and partially in [PSU20, $\S2.2$, $\S7$]. Here, we only spell out the details of the simpler case $r=2$ (cf. $\S$\ref{subsubsection:E_2_Delta-tilde_h_def}). This adiabatic limit approach is needed in other parts of the paper, though, and full details will be provided in all those instances.

 \vspace{2ex}

{\bf (II)\, The relative version of the FAVB}

\vspace{1ex}

In $\S$\ref{subsubsection:relative}, with every holomorphic family $\pi:{\mathcal X}\rightarrow B$ of compact complex $n$-dimensional manifolds $X_t:=\pi^{-1}(t)$, supposed to satisfy a certain condition, over an open ball $B\subset\C^N$ about the origin we associate a holomorphic vector bundle ${\mathcal A}^k\longrightarrow\C\times B$ in every degree $k\in\{0,\dots , 2n\}$. We call it the {\it Fr\"olicher approximating vector bundle (FAVB)} of the family $\pi$ in degree $k$.

We stress right away that the vector bundles ${\mathcal A}^k\longrightarrow\C\times B$ are necessarily {\it trivial} thanks to the following

\begin{Gen-Fact}\label{Gen-Fact:bundle-trivialty_C-B} Any topological vector bundle ${\mathcal A}$ on $\C\times B$ is topologically trivial and any holomorphic vector bundle ${\mathcal A}$ on $\C\times B$ is holomorphically trivial.

\end{Gen-Fact}

 The reason for this is that $\C\times B$ is both {\it contractible} (hence ${\mathcal A}$ must be topologically trivial) and {\it Stein} (so the Grauert-Oka principle -- see [Gra58] -- applies and implies that any holomorphic vector bundle over a Stein manifold that is topologically trivial is also holomorphically trivial).

\vspace{1ex}

The construction of the relative FAVB ${\mathcal A}^k\longrightarrow\C\times B$ in degree $k$ proceeds along the following lines.

Let $X$ be the $C^\infty$ manifold that underlies the fibres $X_t$. The linear operator $$d_{h,\,t}:=h\partial_t + \bar\partial_t :C^\infty_k(X,\,\C)\to C^\infty_{k+1}(X,\,\C)$$ depends now on both $h\in\C$ and $t\in B$ (because it depends on the complex structure $J_t$ of $X_t$) and so does the pointwise linear map $$\theta_{h,\,t}:\Lambda^kT^\star X\to\Lambda^kT^\star X, \hspace{3ex} \theta_{h,\,t}\bigg(\sum_{p+q=k}u_t^{p,\,q}\bigg):= \sum_{p+q=k}h^p\,u_t^{p,\,q},$$ where the $u_t^{p,\,q}$ are the $(p,\,q)$-type components w.r.t. $J_t$ of a given $k$-form $u=\sum_{p+q=k}u_t^{p,\,q}$ on $X$.

In particular, if we fix a $C^\infty$ family $(\omega_t)_{t\in B}$ of Hermitian metrics on the fibres $(X_t)_{t\in B}$ of the family $\pi:{\mathcal X}\rightarrow B$, the corresponding $d_{h,\,t}$-Laplacian in degree $k$: \begin{eqnarray*}\Delta_{h,\,t}:C^\infty_k(X,\,\C)\to C^\infty_k(X,\,\C), \hspace{5ex} \Delta_{h,\,t}:=d_{h,\,t}d_{h,\,t}^\star + d_{h,\,t}^\star d_{h,\,t},\end{eqnarray*} depends in a $C^\infty$ way on $(h,\,t)\in\C\times B$ for every $k\in\{0,\dots , 2n\}$.

We define the fibres of the {\bf Fr\"olicher approximating vector bundle (FAVB)} ${\mathcal A}^k$ over $\C\times B$ of the family $(X_t)_{t\in B}$ in degree $k$ as \begin{eqnarray}\label{eqn:fibres_relative-FAVB_introd}\nonumber{\mathcal A}^k_{h,\,t}:=H^k_{d_{h,\,t}}(X_t,\,\C) \hspace{3ex} & \mbox{if} & \hspace{1ex} (h,\,t)\in\C^\star\times B, \\
  {\mathcal A}^k_{0,\,t}:=\bigoplus\limits_{p+q=k}E_\infty^{p,\,q}(X_t) \hspace{3ex} & \mbox{for} & \hspace{1ex} (0,\,t)\in\{0\}\times B.\end{eqnarray}

The operators in the $C^\infty$ family $(\Delta_{h,\,t})_{(h,\,t)\in\C\times B}$ (cf. $\S$\ref{subsubsection:E_1_Delta-tilde_h_def} for the absolute counterpart) that are of use when $r=1$ are generalised to elliptic pseudo-differential operators with equidimensional kernels:

\vspace{1ex}

$\bullet$\, $(\widetilde\Delta_{h,\,t})_{(h,\,t)\in\C\times B}$ (cf. $\S$\ref{subsubsection:E_2_Delta-tilde_h_def} for the absolute counterpart), that are of use when $r=2$;

\vspace{1ex}

$\bullet$\, $(\widetilde\Delta_{h,\,t}^{(r)})_{(h,\,t)\in\C\times B}$ (cf. $\S$\ref{subsubsection:E_geq3_Delta-tilde_h-r_def} for the absolute counterpart), that are of use when $r\geq 3$.

\vspace{1ex}

However, we stress that $\widetilde\Delta_{h,\,t}$ and $\widetilde\Delta_{h,\,t}^{(r)}$ (with $r\geq 3$) need not depend in a $C^\infty$ (or even a continuous) way on $t\in B$ since their definitions involve projections onto the harmonic spaces of the previous pages in the Fr\"olicher spectral sequence and these dimensions, equal to the dimensions of the corresponding Fr\"olicher cohomology spaces $\oplus_{p+q=k}E_l^{p,\,q}(X_t)$, may change when $t$ varies in $B$. 

To overcome this difficulty, we use instead (see Corollary and Definition \ref{Cor-Def:bundle_family}), for every $k\in\{0,\dots , 2n\}$, the $C^\infty$ families of elliptic (pseudo-)differential operators with equidimensional kernels: \begin{eqnarray}\label{eqn:smooth-families-Laplacians-ht_introd}\nonumber & & \bullet\, (\Delta_{h,\,t})_{(h,\,t)\in\C^\star\times B} \\
\mbox{and} & & \\
\nonumber & & \bullet\, \mbox{for every fixed}\hspace{1ex} t\in B, (\widetilde\Delta_{h,\,t}^{(r)})_{h\in\C}\end{eqnarray} to define the $C^\infty$ vector bundle structure of the relative FAVB ${\mathcal A}^k\longrightarrow\C\times B$ in degree $k$ via the Kodaira-Spencer theory [KS60].

\vspace{1ex}

 The $C^\infty$ vector bundle ${\mathcal A}^k\longrightarrow\C\times B$ has an additional holomorphic vector bundle structure. Over $\C^\star\times B$, this is defined by the cohomology isomorphisms \begin{eqnarray}\label{eqn:theta_ht_introd}\theta_{h,\,t}: H^k_{DR}(X_t,\,\C)\to H^k_{d_{h,\,t}}(X_t,\,\C) \hspace{3ex} \mbox{with} \hspace{1ex} (h,\,t)\in\C^\star\times B\end{eqnarray} viewed collectively as an isomorphism from the constant bundle ${\mathcal H}^k\longrightarrow\C^\star\times B$ of fibre $H^k_{DR}(X,\,\C)$, identified with $H^k_{DR}(X_t,\,\C)$ for every $t\in B$, to ${\mathcal A}^k_{|\C^\star\times B}$.

   The isomorphisms (\ref{eqn:theta_ht_introd}) transport the canonical trivialisation of the constant bundle ${\mathcal H}^k\longrightarrow\C^\star\times B$ to a particular global holomorphic trivialisation of ${\mathcal A}^k_{|\C^\star\times B}$. This is the concrete manifestation in our special case of the following 

   \begin{Gen-Fact}\label{Gen-Fact:bundle-trivialty_C-star-B} Any topological $\C$-vector bundle ${\mathcal A}$ on $\C^\star\times B$ is topologically trivial and any holomorphic vector bundle ${\mathcal A}$ on $\C^\star\times B$ is holomorphically trivial.

   \end{Gen-Fact}

   The reason for this is that $\C^\star\times B$ is homotopy equivalent to $\C^\star$, which in turn is homotopy equivalent to the circle $S^1$. Now, all complex vector bundles on $S^1$ are trivial (see e.g. Example 1.12 on p. 24 in the book [Hat17]). One can then use the Steinness of $\C^*\times B$ and the argument given above for the General Fact \ref{Gen-Fact:bundle-trivialty_C-B} to get the holomorphic triviality of any ${\mathcal A}\longrightarrow\C^\star\times B$ that has already been supposed holomorphic.

\vspace{1ex}

That the global holomorphic trivialisation of the vector bundle ${\mathcal A}^k\longrightarrow\C^\star\times B$ extends to a global holomorphic trivialisation of ${\mathcal A}^k$ over $\C\times B$ is asserted by Corollary and Definition \ref{Cor-Def:bundle_family}. This statement can be loosely reworded as

\begin{The}\label{The:bundle_family_introd} Fix $N\in\N^\star$. Let $\pi:{\mathcal X}\rightarrow B$ be a holomorphic family of compact complex $n$-dimensional manifolds $X_t:=\pi^{-1}(t)$ over an open ball $B\subset\C^N$ about the origin. Let $r$ be the smallest positive integer such that the Fr\"olicher spectral sequence of $X_0$ degenerates at $E_r$. Suppose that $E_r(X_t) = E_\infty(X_t)$ for all $t\in B$.

  Let $(\omega_t)_{t\in B}$ be a $C^\infty$ family of Hermitian metrics on the fibres $(X_t)_{t\in B}$.

 For every $k\in\{0,\dots , 2n\}$, the $C^\infty$ families (\ref{eqn:smooth-families-Laplacians-ht_introd}) of elliptic (pseudo-)differential operators with equidimensional kernels and the cohomology isomorphisms (\ref{eqn:theta_ht_introd}) define a (necessarily trivial) holomorphic vector bundle ${\mathcal A}^k\longrightarrow\C\times B$ of rank $b_k$ with fibres described in (\ref{eqn:fibres_relative-FAVB_introd}) and a particular global holomorphic trivialisation of it. In particular, ${\mathcal A}^k$ and its trivialisation depend on the family $(\omega_t)_{t\in B}$ of metrics.
  
\end{The}

By $b_k$ we mean the $k$-th Betti number of the fibres $X_t$, or equivalently, of the $C^\infty$ manifold $X$ underlying them. The proof of Theorem \ref{The:bundle_family_introd} uses the absolute version of the FAVB on each fibre $X_t$. The main point of this result is that it produces a particular global trivialisation (or equivalently, a particular global holomorphic frame) and certain global $C^\infty$ sections of ${\mathcal A}^k\longrightarrow\C\times B$ that are described in (ii) of Corollary and Definition \ref{Cor-Def:bundle_family}.

A key feature of the relative FAVB of Theorem \ref{The:bundle_family_introd} is that it seems {\it non-canonical} in general: its construction depends on the family of Hermitian metrics chosen on the fibres $X_t$. This ties in with the standard fact in Hodge theory stipulating that, for an arbitrary compact complex manifold $X$ (that can be chosen to be any fibre $X_t$ in our case), the isomorphisms (\ref{eqn:Froelicher_decomposition}) obtained by a classical filtration argument need not be canonical.    

\vspace{2ex}

\subsubsection{{\bf Generalised strongly Gauduchon metrics adapted to the various pages of the Fr\"olicher spectral sequence}}\label{subsubsection:E_r-sG_introd}

The second main idea of this paper is a relaxation (cf. Definition \ref{Def:E_r-sG}) of the notion of strongly Gauduchon (sG) metric introduced in [Pop09] and [Pop13]. Starting from the observation that, for every Gauduchon metric $\gamma$ on a given compact complex $n$-dimensional manifold $X$, the $(n,\,n-1)$-form $\partial\gamma^{n-1}$ is $E_r$-closed (i.e. represents an $E_r$-cohomology class on the $r$-th page of the Fr\"olicher spectral sequence of $X$) for every $r\in\N^\star$, we call $\gamma$ an {\bf $E_r$-sG metric} if $\partial\gamma^{n-1}$ is {\bf $E_r$-exact} (i.e. represents the zero $E_r$-cohomology class on the $r$-th page of the Fr\"olicher spectral sequence of $X$). Any $X$ that carries an $E_r$-sG metric is called an {\bf $E_r$-sG manifold}.

For the reader's convenience, we recall in Proposition \ref{Prop:E_r-closed-exact_conditions} how the $E_r$-closedness and $E_r$-exactness conditions translate into explicit terms. For every $r\in\N^\star$, the $E_r$-sG condition implies the $E_{r+1}$-sG condition, while the strongest of them, the $E_1$-sG condition, is equivalent to the sG condition.

The two main constructions of this paper are brought together in the following result (see Theorem \ref{The:limits_E_r_sG} for a more precise statement).

\begin{The}\label{The:limits_E_r_sG_introd} If in a holomorphic family $(X_t)_{t\in B}$ of compact complex manifolds all the fibres $X_t$ with $t\in B\setminus\{0\}$ are {\bf $\partial\bar\partial$-manifolds}, then the limiting fibre $X_0$ is an {\bf $E_r$-sG manifold}, where $r$ is the smallest positive integer such that the Fr\"olicher spectral sequence of $X_0$ degenerates at $E_r$.

\end{The}

Recall that a {\it $\partial\bar\partial$-manifold} is, by definition, a compact complex manifold $X$ that satisfies the {\it $\partial\bar\partial$-lemma} in the following sense: 

\vspace{1ex}

\noindent {\it for every $C^{\infty}$ $d$-closed pure-type form $u$ on $X$, the following exactness conditions are equivalent: \begin{equation}\label{eqn:dd-bar_def}u\in\mbox{Im}\,d \hspace{1ex} \Leftrightarrow \hspace{1ex} u\in\mbox{Im}\,\partial \hspace{1ex} \Leftrightarrow \hspace{1ex} u\in\mbox{Im}\,\bar\partial \hspace{1ex} \Leftrightarrow \hspace{1ex}  u\in\mbox{Im}\,\partial\bar\partial.\end{equation}} \noindent The $\partial\bar\partial$-property is equivalent to all the canonical linear maps $H^{p,\,q}_{BC}(X,\,\C)\longrightarrow H^{p,\,q}_A(X,\,\C)$, from the Bott-Chern to the Aeppli cohomology, being isomorphisms. Since both of these cohomologies can be computed using either smooth forms or currents, the $\partial\bar\partial$-property is also equivalent to the equivalences (\ref{eqn:dd-bar_def}) holding for every $d$-closed pure-type current on $X$.  

A standard result in Hodge theory asserts that every compact K\"ahler manifold is a $\partial\bar\partial$-manifold. Moreover, every {\it class} ${\mathcal C}$ manifold (by definition, these are the compact complex manifolds that are bimeromorphically equivalent to compact K\"ahler manifolds), hence also every Moishezon manifold, is a $\partial\bar\partial$-manifold, but the class of $\partial\bar\partial$-manifolds strictly contains the {\it class} ${\mathcal C}$. (See e.g. [Pop14] for further details.)

Another standard result in Hodge theory ensures that the Fr\"olicher spectral sequence of any $\partial\bar\partial$-manifold $X$ degenerates at $E_1$ (the best possible degeneration property): $E_1(X) = E_\infty(X)$.

\vspace{2ex}

A result stronger than Theorem \ref{The:limits_E_r_sG_introd} was proved in Proposition 4.1 of [Pop09]: any deformation limit of $\partial\bar\partial$-manifolds is a strongly Gauduchon (sG) manifold. In the present paper, we use our relative Fr\"olicher approximating vector bundle (FAVB) of Corollary and Definition \ref{Cor-Def:bundle_family} (or, alternatively, only the Kodaira-Spencer theory of [KS60] and our Laplacians -- see Remark \ref{Rem:projections_smooth-var_suffices}) to obtain the possibly weaker $E_r$-sG conclusion (for some specified $r\geq 1$) on the limiting fibre $X_0$ under the same $\partial\bar\partial$-assumption on the other fibres. (Moreover, Remark \ref{Rem:referee-suggestion} and Theorem \ref{The:limits_E_r_sG_simplified} yield the $E_3$-sG conclusion on $X_0$.) However, we show that any of these weaker conclusions on the deformation limits of $\partial\bar\partial$-manifolds yields the same optimal conclusion, captured in Theorem \ref{The:moidef}, on the deformation limits of projective and Moishezon manifolds as the one obtained in [Pop09] and [Pop10]. Moreover, the new method introduced in the present paper has the advantage of being more conceptual than the ad hoc arguments of [Pop09]. It effectively puts those arguments on a more theoretical footing via the machinery of the Fr\"olicher spectral sequence.

\vspace{2ex}

Besides Theorem \ref{The:limits_E_r_sG_introd}, the other main building block (cf. Theorem \ref{The:boundedness}) of the proof of Theorem \ref{The:moidef} is the use of a $C^\infty$ family $(\gamma_t)_{t\in B}$ of $E_r$-sG metrics on the fibres $(X_t)_{t\in B}$, whose existence is mainly guaranteed by Theorem \ref{The:limits_E_r_sG_introd}, to uniformly control the volumes of the relative (i.e. contained in the fibres) divisors that form an irreducible component of the relative Barlet space of divisors (cf. [Bar75]) associated with the family $(X_t)_{t\in B}$. Finitely many integrations by parts are used.

As the proofs show, in Theorems \ref{The:limits_E_r_sG_introd} and \ref{The:boundedness} the $\partial\bar\partial$-assumption is only needed in a select few bidegrees and with a select few of the implications (\ref{eqn:dd-bar_def}). However, we made the full $\partial\bar\partial$-assumption on $X_t$ with $t\neq 0$ to make the statements shorter and less technical.

\section{h-theory for the Fr\"olicher spectral sequence}\label{Froelicher_h}

In this section, we construct certain $C^\infty$ families of elliptic (pseudo-)differential operators (cf. $\S$\ref{subsubsection:E_1_Delta-tilde_h_def}, $\S$\ref{subsubsection:E_2_Delta-tilde_h_def} and $\S$\ref{subsubsection:E_geq3_Delta-tilde_h-r_def}) that we then use to construct in $\S$\ref{subsection:approx_vb} the Fr\"olicher approximating vector bundle (FAVB) associated with a compact complex manifold (cf. $\S$\ref{subsubsection:absolute}) and, respectively, with a holomorphic family of such manifolds (cf. $\S$\ref{subsubsection:relative}). Each of the families of operators that will be constructed is suited to one of the pages of the Fr\"olicher spectral sequence (FSS) according to which page degeneration occurs on.

\subsection{Brief review of the Kodaira-Spencer theory of $C^\infty$ families of elliptic operators}\label{subsection:review_K-S}

We will follow the presentation in [Kod86] to recall in some detail the classical results of Kodaira and Spencer in [KS60] showing how a $C^\infty$ family of elliptic (pseudo-)differential operators induces a structure of $C^\infty$ vector bundle on the family of kernels of these operators whenever the dimensions of the kernels stay constant. The notation we adopt differs in places from that of [Kod86, $\S7.1.$ and $\S7.2.$].

\vspace{2ex}

$(1)$\, Let $(X_t)_{t\in B}$ be a $C^\infty$ family of compact complex $n$-dimensional manifolds and let $(V_t)_{t\in B}$ be a $C^\infty$ family of holomorphic vector bundles over the $X_t$'s, where $B$ is, for example, an open ball about the origin in some $\C^N$. The family $(V_t)_{t\in B}$ being $C^\infty$ means that there exists a $C^\infty$ complex vector bundle $\widetilde\pi:{\mathcal V}\longrightarrow{\mathcal X}$ over the total space ${\mathcal X}$ of the family $\pi:{\mathcal X}\longrightarrow B$ of compact complex manifolds $\bigg(X_t=\pi^{-1}(t)\bigg)_{t\in B}$ such that $V_t=\widetilde\pi^{-1}(X_t)$ for every $t\in B$.

We denote by $C^\infty_{p,\,q}(X_t,\,V_t)$ the $\C$-vector space of $C^\infty$ $(p,\,q)$-forms on $X_t$ with values in $V_t$. Let $(\omega_t)_{t\in B}$ be a $C^\infty$ family of Hermitian metrics on the $X_t$'s and let $(h_t)_{t\in B}$ be a $C^\infty$ family of Hermitian metrics on the fibres of the $V_t$'s. For every $t\in B$, the metrics $\omega_t$ and $h_t$ induce an $L^2$-inner product $\langle\langle\,\,,\,\, \rangle\rangle_t$ on $C^\infty_{p,\,q}(X_t,\,V_t)$, which, in turn, defines the adjoint $\bar\partial_t^\star$ of the operator $\bar\partial_t$ associated with the holomorphic structure of $V_t$. 

 The pair $(\bar\partial_t,\,\bar\partial_t^\star)_{t\in B}$ induces the $\bar\partial$-Laplacian $$\Delta''_t:=\bar\partial_t\bar\partial_t^\star + \bar\partial_t^\star\bar\partial_t:C^\infty_{p,\,q}(X_t,\,V_t)\longrightarrow C^\infty_{p,\,q}(X_t,\,V_t),  \hspace{5ex} t\in B.$$ Denoting by ${\mathcal H}^{p,\,q}_{\Delta''_t}(X_t,\,V_t)$ the kernel of $\Delta''_t$ (the $\Delta''_t$-harmonic space), one has an $L^2_{\omega_t,\,h_t}$-orthogonal decomposition $$C^\infty_{p,\,q}(X_t,\,V_t) = {\mathcal H}^{p,\,q}_{\Delta''_t}(X_t,\,V_t)\oplus\mbox{Im}\,\Delta''_t$$ and the associated $L^2_{\omega_t,\,h_t}$-orthogonal projection $$H_t:C^\infty_{p,\,q}(X_t,\,V_t)\longrightarrow{\mathcal H}^{p,\,q}_{\Delta''_t}(X_t,\,V_t), \hspace{5ex} t\in B.$$ 

One way of expressing the fact that the $C^\infty$ family $(\Delta''_t)_{t\in B}$ of elliptic differential operators induces a natural structure of $C^\infty$ vector bundle on the family $({\mathcal H}^{p,\,q}_{\Delta''_t}(X_t,\,V_t))_{t\in B}$ of its kernels when the dimensions of these kernels are all equal is the following

\begin{The}\label{The:K-S_d-bar_t_bundle}([Kod86, Theorem 7.9.]) If $\mbox{dim}\,{\mathcal H}^{p,\,q}_{\Delta''_t}(X_t,\,V_t)$ is independent of $t\in B$, the operator $H_t$ is $C^\infty$ differentiable with respect to $t\in B$. 

\end{The} 

According to [Kod86, Definition 7.5.], the operator $H_t$ being $C^\infty$ differentiable with respect to $t\in B$ ($=$ the family $(H_t)_{t\in B}$ of operators being $C^\infty$) means that for every $C^\infty$ family $(\psi_t)_{t\in B}$ of forms $\psi_t\in C^\infty_{p,\,q}(X_t,\,V_t)$, the family $(H_t(\psi_t))_{t\in B}$ of their images under the $H_t$'s is again $C^\infty$.

\vspace{2ex}

$(2)$\, More generally, let $(V_t)_{t\in B}$ be a $C^\infty$ family of $C^\infty$ complex vector bundles over the fibres $X_t$ of a $C^\infty$ family $\pi:{\mathcal X}\longrightarrow B$ of compact complex $n$-dimensional manifolds. For every $t\in B$, we denote by $C^\infty(X_t,\,V_t)$ the space of $C^\infty$ sections of $V_t$ over $X_t$. Let $(\omega_t)_{t\in B}$ and $(h_t)_{t\in B}$ be $C^\infty$ families of Hermitian metrics on the $X_t$'s, respectively on the fibres of the $V_t$'s.

The result analogous in this more general context to Theorem \ref{The:K-S_d-bar_t_bundle} for an arbitrary $C^\infty$ family of elliptic operators with kernels of constant dimension is the following

\begin{The}\label{The:K-S_E_t_bundle}([Kod86, Theorem 7.4.]) Let $(E_t)_{t\in B}$ be a $C^\infty$ family of elliptic linear differential operators $$E_t:C^\infty(X_t,\,V_t)\longrightarrow C^\infty(X_t,\,V_t),   \hspace{5ex} t\in B.$$ If $\mbox{dim}\,\ker(E_t)$ is independent of $t\in B$, the $L^2$-orthogonal projection operator $$H_t:C^\infty(X_t,\,V_t)\longrightarrow\ker(E_t)$$ is $C^\infty$ differentiable with respect to $t\in B$. 

\end{The}

\vspace{2ex}

$(3)$\, To make even more explicit the sense in which a $C^\infty$ family $(E_t)_{t\in B}$ of elliptic differential operators induces a natural structure of $C^\infty$ vector bundle on the family $(\ker(E_t))_{t\in B}$ of its kernels whenever the dimensions of these kernels are equal to one another, we now recall another classical result of Kodaira-Spencer. 

Suppose, for example, that we are in the context of the above $(1)$. For every $t\in B$, let $H^0(X_t,\,V_t)$ be the space of global holomorphic sections of the vector bundle $V_t\longrightarrow X_t$. 

\begin{The}\label{The:K-S_basis_E_t_bundle}([Kod86, Lemma 7.8., p. 355]) If $\mbox{dim}\,H^0(X_t,\,V_t):=d$ is independent of $t\in B$, there exists $\varepsilon>0$ such that, for every $t\in B$ with $|t|<\varepsilon$, there exists a basis $$\{\varphi_t^{(1)},\dots , \varphi_t^{(d)}\}$$ of $H^0(X_t,\,V_t)$ with the property that $\varphi_t^{(q)}$ is $C^\infty$ differentiable with respect to $t$ for each $q=1,\dots , d$. 

\end{The}

The $C^\infty$ family of bases $\{\varphi_t^{(1)},\dots , \varphi_t^{(d)}\}$ given by this statement defines a local $C^\infty$ trivialisation of the induced vector bundle over $B$ whose fibres are $(H^0(X_t,\,V_t))_{t\in B}$.

\vspace{2ex}

\noindent {\it Proof of Theorem \ref{The:K-S_basis_E_t_bundle}.} We reproduce the proof given in [Kod86, p. 355] in order to stress how immediate a corollary of Theorem \ref{The:K-S_d-bar_t_bundle} this result is.

Let $\{\varphi_1,\dots , \varphi_d\}$ be a basis of $H^0(X_0,\,V_0)$. Since each holomorphic vector bundle $V_t\longrightarrow X_t$ is the restriction to $X_t$ of a same $C^\infty$ vector bundle ${\mathcal V}\longrightarrow{\mathcal X}$, there exists, for each $t\in B$, a family of sections $\psi_t^{(1)},\dots , \psi_t^{(d)}\in C^\infty(X_t,\,V_t):=C^\infty_{0,\,0}(X_t,\,V_t)$ such that:

\vspace{2ex}

(a)\, $\psi_0^{(q)} = \varphi_q$ for all $q=1,\dots , d$;

\vspace{0.5ex}

(b)\, the family $(\psi_t^{(q)})_{t\in B}$ is $C^\infty$ differentiable with respect to $t$ for every $q=1,\dots , d$.

\vspace{2ex}

Meanwhile, we have the $L^2$-orthogonal decomposition: $$C^\infty(X_t,\,V_t) = {\mathcal H}^{0,\,0}_{\Delta''_t}(X_t,\,V_t)\oplus\mbox{Im}\,\Delta''_t, \hspace{5ex} t\in B,$$ where $\displaystyle{\mathcal H}^{0,\,0}_{\Delta''_t}(X_t,\,V_t): = \ker\bigg(\Delta''_t:C^\infty(X_t,\,V_t)\longrightarrow C^\infty(X_t,\,V_t)\bigg)$

\vspace{1ex}

\hspace{16ex} $\displaystyle= \bigg\{\varphi\in C^\infty(X_t,\,V_t)\,\mid\,\bar\partial_t\varphi=0\bigg\} = H^0(X_t,\,V_t), \hspace{5ex} t\in B.$ 

\vspace{2ex}

\noindent Indeed, in bidegree $(0,\,0)$, $\bar\partial_t^\star$ vanishes identically. 

For every $t\in B$, let $H_t:C^\infty(X_t,\,V_t)\longrightarrow{\mathcal H}^{0,\,0}_{\Delta''_t}(X_t,\,V_t)$ be the $L^2$-orthogonal projection. We set: $$\varphi_t^{(q)}:=H_t(\psi_t^{(q)}), \hspace{5ex} t\in B, \hspace{1ex} q=1,\dots , d.$$ 

In particular, for $t=0$ we get: $$\varphi_0^{(q)} = \varphi_q, \hspace{5ex} q=1,\dots , d,$$ since $\psi_0^{(q)} = \varphi_q\in H^0(X_0,\,V_0) = {\mathcal H}^{0,\,0}_{\Delta''_0}(X_0,\,V_0)$ by the above (a). 

Since $\mbox{dim}\,H^0(X_t,\,V_t):=d$ is independent of $t\in B$ (by hypothesis), $H_t$ is $C^\infty$ differentiable with respect to $t\in B$ by Theorem \ref{The:K-S_d-bar_t_bundle}. Therefore, since moreover each family $(\psi_t^{(q)})_{t\in B}$ is $C^\infty$ differentiable with respect to $t$ (by the above (b)), we conclude that the family $(\varphi_t^{(q)})_{t\in B}$ is $C^\infty$ differentiable with respect to $t$ for every $q=1,\dots , d$. 

In particular, since the sections $\varphi_0^{(1)},\dots , \varphi_0^{(d)}$ are linearly independent, the sections $\varphi_t^{(1)},\dots , \varphi_t^{(d)}$ are again linearly independent for every $t\in B$ such that $|t|<\varepsilon$ if $\varepsilon>0$ is small enough. \hfill $\Box$

\vspace{2ex}

The above conclusions of the Kodaira-Spencer theory remain valid when the $C^\infty$ family of elliptic {\it differential} operators is replaced by a $C^\infty$ family of elliptic {\it pseudo-differential} operators, as shown in [Mas18]. In this paper, we will apply these conclusions (see proof of Theorem \ref{The:limits_E_r_sG}) in the following form.

\begin{Cor}\label{Cor:K-S_conclusions_to-be-applied} Let $\pi:{\mathcal X}\rightarrow B$ be a $C^\infty$ family of compact complex manifolds $X_t:=\pi^{-1}(t)$, with $t\in B$, over an open ball $B\subset\C^N$ about the origin. For every $t\in B$, let $V_t\longrightarrow X_t$ be a $C^\infty$ complex vector bundle such that $(V_t)_{t\in B}$ is a $C^\infty$ family. Let $(E_t)_{t\in B}$ be a $C^\infty$ family of elliptic linear (pseudo-)differential operators $$E_t:C^\infty(X_t,\,V_t)\longrightarrow C^\infty(X_t,\,V_t),   \hspace{5ex} t\in B.$$

If the dimension $d:=\mbox{dim}\,\ker(E_t)$ is {\bf independent of $t\in B$}, the Kodaira-Spencer theory induces a $C^\infty$ complex vector bundle ${\mathcal K}\longrightarrow B$ with the following properties:

\vspace{1ex}

(i)\, for every $t\in B$, the fibre of ${\mathcal K}$ over $t$ is ${\mathcal K}_t=\ker(E_t)$;

\vspace{1ex}

(ii)\, for every $t_0\in B$, there exists $\varepsilon_{t_0}>0$ such that a $C^\infty$ trivialisation of ${\mathcal K}$ over the open ball $B(t_0,\,\varepsilon_{t_0})\subset B$ centred at $t_0$ of radius $\varepsilon_{t_0}$ is defined by a family of bases $$\bigg(\{\varphi_t^{(1)},\dots , \varphi_t^{(d)}\}\bigg)_{t\in B(t_0,\,\varepsilon_{t_0})}$$ of the respective fibres $(\ker(E_t))_{t\in B(t_0,\,\varepsilon_{t_0})}$ with the property that $\varphi_t^{(q)}$ is $C^\infty$ differentiable with respect to $t$ for each $q=1,\dots , d$. 

\vspace{1ex}

(iii)\, for every $C^\infty$ family $s=(s_t)_{t\in B}$ of sections $s_t\in C^\infty(X_t,\,V_t)$, the following property holds. If $(\omega_t)_{t\in B}$ is a $C^\infty$ family of Hermitian metrics on the $X_t$'s, $(h_t)_{t\in B}$ is a $C^\infty$ family of Hermitian metrics on the fibres of the $V_t$'s and $$H_t:C^\infty(X_t,\,V_t)\longrightarrow\ker(E_t)$$ are the {\bf $L^2_{\omega_t,\,h_t}$-orthogonal projection operators}, then the family $$\sigma=\bigg(\sigma_t:=H_t(s_t)\bigg)_{t\in B}$$ defines a $C^\infty$ section of the vector bundle ${\mathcal K}$ over $B$.

\end{Cor}

\subsection{Constructions of approximations of the Laplacians $\Delta''$, $\widetilde\Delta$ and $\widetilde\Delta^{(r)}$ with $r\geq 3$ whose kernels are isomorphic to the $E_r$-cohomology spaces by $C^\infty$ families of operators $(\Delta_h)_{t\in B}$, $(\widetilde\Delta_h)_{t\in B}$ and $(\widetilde\Delta^{(r)}_h)_{t\in B}$ whose term for $h=0$ is the original Laplacian}\label{subsection:Laplacians_approximations} \hspace{1ex}

\vspace{1ex}

We will discuss separately the cases of the pages $E_1$, $E_2$ and $E_r$ with $r\geq 3$ of the Fr\"olicher spectral sequence.

\subsubsection{{\bf First FSS page: the Laplace-type differential operators $\Delta_h$}}\label{subsubsection:E_1_Delta-tilde_h_def} Let $(X,\,\omega)$ be an $n$-dimensional compact complex Hermitian manifold.

Recall (cf. the above $\S$\ref{subsection:reminder} or [Pop17]) that $(\Delta_h)_{h\in\C}$ is a $C^\infty$ family of elliptic differential operators such that $\Delta_0=\Delta''$. So, the $\Delta_h$'s can be regarded as an {\it approximation} (allowing for more flexibility) of the standard $\bar\partial$-Laplacian $\Delta''$. The kernel of $\Delta_0=\Delta''$ is classically isomorphic in every degree $k$ to the Dolbeault cohomology space $\oplus_{p+q=k}H^{p,\,q}_{\bar\partial}(X,\,\C)$ of degree $k$ of $X$ (thus, to the first page of the Fr\"olicher spectral sequence). Moreover, denoting by $h^{p,\,q}_{\bar\partial}(X)$, respectively $b_k(X)$, the Hodge numbers, respectively the Betti numbers of $X$, we have: $$\sum\limits_{p+q=k}h^{p,\,q}_{\bar\partial}(X)=\mbox{dim}\ker\Delta_0\geq\mbox{dim}\ker\Delta_h=b_k(X), \hspace{5ex} h\in\C\setminus\{0\},$$ in every degree $k$, where the inequality is an equality for all $k$ if and only if the Fr\"olicher spectral sequence of $X$ degenerates at the first page (a fact denoted by $E_1(X)=E_\infty(X)$).

In the case when $E_1(X)=E_\infty(X)$, $(\Delta_h)_{h\in\C}$ is the $C^\infty$ family of elliptic differential operators that we need for the construction of the FAVB of $X$. To unify the notation throughout $\S$\ref{subsubsection:E_1_Delta-tilde_h_def}, $\S$\ref{subsubsection:E_2_Delta-tilde_h_def} and $\S$\ref{subsubsection:E_geq3_Delta-tilde_h-r_def}, we will also denote $\widetilde\Delta_h^{(1)}:=\Delta_h$.

%Subsection \ref{subsubsection:E_2_Delta-tilde_h_def} lays the groundwork for the construction of the absolute version of the FAVB in $\S$\ref{subsubsection:absolute}. 

\subsubsection{{\bf Second FSS page: the pseudo-differential Laplacians $\widetilde\Delta_h$}}\label{subsubsection:E_2_Delta-tilde_h_def}

However, when the Fr\"olicher spectral sequence of $X$ degenerates at the second page (a fact denoted by $E_2(X)=E_\infty(X)$) but not at the first page, we will replace the family $(\Delta_h)_{h\in\C}$ by a new $C^\infty$ family $(\widetilde\Delta_h)_{h\in\C}$ of elliptic pseudo-differential operators that we now set out to construct.

In fact, we now introduce and analyse an approximation of the pseudo-differential Laplacian

$$\widetilde\Delta=\partial p''\partial^\star + \partial^\star p''\partial + \Delta'':\C^\infty_{p,\,q}(X,\,\C)\longrightarrow \C^\infty_{p,\,q}(X,\,\C),  \hspace{3ex} p,q=0,\dots , n,$$

\noindent introduced in [Pop16] and proved there to define a Hodge theory for the second page of the Fr\"olicher spectral sequence, namely a Hodge isomorphism

$${\mathcal H}^{p,\,q}_{\widetilde\Delta}(X,\,\C):= \ker(\widetilde\Delta:\C^\infty_{p,\,q}(X,\,\C)\longrightarrow \C^\infty_{p,\,q}(X,\,\C))\simeq E_2^{p,\,q}(X)$$

\noindent in every bidegree $(p,\,q)$. Note that $\widetilde\Delta=(\partial p'')(\partial p'')^\star + (p''\partial)^\star(p''\partial) + \Delta''$, so we will approximate $\partial p''$ and $p''\partial$ by adding to each a small $h$-multiple of its conjugate, while still approximating the second term $\Delta''$ of $\widetilde\Delta$ by $\Delta_h$.

\begin{Def}\label{Def:h-Delta-tilde} Let $(X,\,\omega)$ be a compact complex Hermitian manifold with $\mbox{dim}_\C X=n$. For every $h\in\C$ and every $k=0,\dots , 2n$, we define the operator $$\widetilde\Delta_h=(\partial p''+h\,\bar\partial p')(\partial p''+h\,\bar\partial p')^\star + (p''\partial + h\,p'\bar\partial)^\star(p''\partial + h\,p'\bar\partial) + \Delta_h:\C^\infty_k(X,\,\C)\longrightarrow \C^\infty_k(X,\,\C),$$ where \begin{eqnarray*}p' & = & p'_\omega:C^\infty_{p,\,q}(X,\,\C)\longrightarrow\ker\bigg(\Delta':\C^\infty_{p,\,q}(X,\,\C)\longrightarrow \C^\infty_{p,\,q}(X,\,\C)\bigg):={\mathcal H}^{p,\,q}_{\Delta'}(X,\,\C), \\
    p'' & = & p''_\omega:C^\infty_{p,\,q}(X,\,\C)\longrightarrow\ker\bigg(\Delta'':\C^\infty_{p,\,q}(X,\,\C)\longrightarrow \C^\infty_{p,\,q}(X,\,\C)\bigg):={\mathcal H}^{p,\,q}_{\Delta''}(X,\,\C)\end{eqnarray*} are the orthogonal projections onto the $\Delta'$-, resp. $\Delta''$-harmonic spaces of any fixed bidegree $(p,\,q)$.

  These projections are then extended by linearity to $$p'=p'_\omega:C^\infty_k(X,\,\C)\longrightarrow{\mathcal H}^k_{\Delta'}(X,\,\C), \hspace{3ex} p''=p''_\omega:C^\infty_k(X,\,\C)\longrightarrow{\mathcal H}^k_{\Delta''}(X,\,\C),$$ where ${\mathcal H}^k_{\Delta'}(X,\,\C):=\oplus_{p+q=k}{\mathcal H}^{p,\,q}_{\Delta'}(X,\,\C)$ and ${\mathcal H}^k_{\Delta''}(X,\,\C):=\oplus_{p+q=k}{\mathcal H}^{p,\,q}_{\Delta''}(X,\,\C)$.

\end{Def}

For every $h\in\C$, $\widetilde\Delta_h$ is a non-negative, self-adjoint pseudo-differential operator and $\widetilde\Delta_0 = \widetilde\Delta$. Further properties include the following.

\begin{Lem}\label{Lem:Delta_h-tilde_prop} For every $h\in\C\setminus\{0\}$, $\widetilde\Delta_h$ is an elliptic pseudo-differential operator whose kernel is \begin{eqnarray}\label{eqn:h-kernels}\nonumber\ker\widetilde\Delta_h & = & \ker(\partial p''+h\,\bar\partial p')^\star\cap\ker(p''\partial + h\,p'\bar\partial)\cap\ker d_h\cap\ker d_h^\star    \\
     & = & \ker d_h\cap\ker d_h^\star =\ker\Delta_h, \hspace{6ex} k=0,\dots , 2n.\end{eqnarray}

  Hence, the $3$-space orthogonal decompositions induced by $\widetilde\Delta_h$ and $\Delta_h$ coincide when $h\in\C\setminus\{0\}$:

\begin{equation}\label{eqn:3-space_decomp_h}C^\infty_k(X,\,\C) = \ker\widetilde\Delta_h\oplus\mbox{Im}\,d_h\oplus\mbox{Im}\,d^\star_h, \hspace{6ex} k=0,\dots , 2n,\end{equation}

  \noindent where $\ker d_h=\ker\widetilde\Delta_h\oplus\mbox{Im}\,d_h$, $\ker d_h^\star = \ker\widetilde\Delta_h\oplus\mbox{Im}\,d^\star_h$ and $\mbox{Im}\,\widetilde\Delta_h = \mbox{Im}\,d_h\oplus\mbox{Im}\,d^\star_h$.

  Consequently, we have the Hodge isomorphism: \begin{eqnarray}\label{eqn:Hodge_h-tilde}\hspace{9ex} {\mathcal H}^k_{\widetilde\Delta_h}(X,\,\C) = {\mathcal H}^k_{\Delta_h}(X,\,\C)\simeq H^k_{d_h}(X,\,\C), \hspace{3ex} k=0,\dots , 2n, \hspace{2ex} h\in\C\setminus\{0\}.\end{eqnarray}

  Moreover, the decomposition (\ref{eqn:3-space_decomp_h}) is stable under $\widetilde\Delta_h$, namely \begin{equation}\label{eqn:Delta-tilde_h_stability}\widetilde\Delta_h(\mbox{Im}\,d_h)\subset\mbox{Im}\,d_h  \hspace{3ex} \mbox{and} \hspace{3ex} \widetilde\Delta_h(\mbox{Im}\,d^\star_h)\subset\mbox{Im}\,d^\star_h.\end{equation}

\end{Lem}

\noindent {\it Proof.} The first identity in (\ref{eqn:h-kernels}) follows immediately from the fact that $\widetilde\Delta_h$ is a sum of non-negative operators of the shape $A^\star A$ and $\ker(A^\star A) = \ker A$ for every $A$, since $\langle\langle A^\star Au,\,u\rangle\rangle = ||Au||^2$.

To prove the second identity in (\ref{eqn:h-kernels}), we will prove the inclusions $\ker d_h\subset\ker(p''\partial + h\,p'\bar\partial)$ and $\ker d_h^\star\subset\ker(\partial p''+h\,\bar\partial p')^\star$.

Let $u=\sum_{r+s=k}u^{r,\,s}$ be a smooth $k$-form such that $d_h u=0$. This amounts to $h\partial u^{r,\,s} + \bar\partial u^{r+1,\,s-1} = 0$ whenever $r+s=k$. Applying $p'$ and respectively $p''$, we get

$$p'\bar\partial u^{r+1,\,s-1} = 0  \hspace{3ex} \mbox{and} \hspace{3ex} p''\partial u^{r,\,s} =0,  \hspace{3ex} \mbox{whenever} \hspace{1ex} r+s=k,$$

\noindent since $h\neq 0$, while $p'\partial =0$ and $p''\bar\partial=0$. Hence,

$$(p''\partial + h\,p'\bar\partial)\,u = \sum_{r+s=k}(p''\partial u^{r,\,s} + h\,p'\bar\partial u^{r+1,\,s-1})=0.$$

\noindent This proves the inclusion $\ker d_h\subset\ker(p''\partial + h\,p'\bar\partial)$.

The ellipticity of the (pseudo)-differential operators $\Delta_h$ and $\widetilde\Delta_h$, combined with the compactness of the manifold $X$, implies that the images of $d_h$ and $\partial p''+h\,\bar\partial p'$ are {\it closed} in $C^\infty_k(X,\,\C)$. Hence, these images coincide with the orthogonal complements of the kernels of the adjoint operators $d_h^\star$ and $(\partial p''+h\,\bar\partial p')^\star$. Therefore, proving the inclusion $\ker d_h^\star\subset\ker(\partial p''+h\,\bar\partial p')^\star$ is equivalent to proving the inclusion $\mbox{Im}\,(\partial p''+h\,\bar\partial p')\subset\mbox{Im}\,d_h$. (Actually, the closedness of these images is not needed here, we would have taken closures otherwise.)

Let $u=\partial p''v + h\,\bar\partial p'v$ be a smooth $k$-form lying in the image of $\partial p''+h\,\bar\partial p'$. Since $\partial p'=0$ and $\bar\partial p''=0$, while $h\neq 0$, we get

$$u = (h\partial)\,(\frac{1}{h}\,p''v + h\,p'v) + \bar\partial\,(\frac{1}{h}\,p''v + h\,p'v) = d_h\,(\frac{1}{h}\,p''v + h\,p'v)\in\mbox{Im}\,d_h.$$

This completes the proof of (\ref{eqn:h-kernels}).

Since $\Delta_h$ commutes with both $d_h$ and $d_h^\star$, to prove (\ref{eqn:Delta-tilde_h_stability}) it suffices to prove the stability of $\mbox{Im}\,d_h$ and $\mbox{Im}\,d_h^\star$ under $\widetilde\Delta_h - \Delta_h$. Now, since $(p''\partial + hp'\bar\partial)\,d_h=0$ (immediate verification), we get

$$(\widetilde\Delta_h - \Delta_h)\,d_h = (\partial p'' + h\bar\partial p')(\partial p'' + h\bar\partial p')^\star(h\partial+\bar\partial).$$

\noindent Since $\mbox{Im}\,(\partial p'' + h\bar\partial p')\subset\mbox{Im}\,d_h$ (as seen above), we get $(\widetilde\Delta_h - \Delta_h)(\mbox{Im}\,d_h)\subset\mbox{Im}\,d_h$. Similarly, an immediate verification shows that $(\partial p'' + h\bar\partial p')^\star d_h^\star=0$. Consequently, 

$$(\widetilde\Delta_h - \Delta_h)\,d_h^\star = (p''\partial + hp'\bar\partial)^\star(p''\partial + hp'\bar\partial)d_h^\star.$$

\noindent Meanwhile, $\mbox{Im}\,(p''\partial + hp'\bar\partial)^\star\subset\mbox{Im}\,d_h^\star$ (since this is equivalent to the inclusion $\ker d_h\subset\ker(p''\partial + h\,p'\bar\partial)$ that was proved above). Therefore, $(\widetilde\Delta_h - \Delta_h)(\mbox{Im}\,d_h^\star)\subset\mbox{Im}\,d_h^\star$. The proof of (\ref{eqn:Delta-tilde_h_stability}) is complete.

The remaining statements follow from the standard elliptic theory as in [Pop17].  \hfill $\Box$

\vspace{2ex}

We sum up these conclusions in the following statement, where the properties of the operators $\Delta_h$ of $\S$\ref{subsubsection:E_1_Delta-tilde_h_def} are repeated for the sake of comparison with those of the operators $\widetilde\Delta_h$ introduced in this $\S$\ref{subsubsection:E_2_Delta-tilde_h_def}.

 \begin{Conc}\label{Conc:pseudodiff_2} Let $(X,\,\omega)$ be a compact complex Hermitian manifold with $\mbox{dim}_\C X=n$. For every degree $k\in\{0,\dots , 2n\}$, there exist $C^\infty$ families of elliptic differential operators $(\Delta_h)_{h\in\C}$ and, respectively, elliptic pseudo-differential operators $(\widetilde\Delta_h)_{h\in\C}$ from $C^\infty_k(X,\,\C)$ to $C^\infty_k(X,\,\C)$ such that

    \vspace{1ex}

   (i)\, $\Delta_0 = \Delta''$ and $\widetilde\Delta_0 = \widetilde\Delta$;

   \vspace{1ex}

   (ii)\, ${\mathcal H}^k_{\Delta_h}(X,\,\C) = {\mathcal H}^k_{\widetilde\Delta_h}(X,\,\C)\simeq H^k_{d_h}(X,\,\C)$ \hspace{3ex} for all $h\in\C\setminus\{0\}$;

   \vspace{1ex}

   (iii)\, ${\mathcal H}^k_{\Delta_0}(X,\,\C)\simeq\bigoplus_{p+q=k}H^{p,\,q}_{\bar\partial}(X,\,\C)$ \hspace{3ex} and \hspace{3ex} ${\mathcal H}^k_{\widetilde\Delta_0}(X,\,\C)\simeq\bigoplus_{p+q=k}E_2^{p,\,q}(X)$.

\end{Conc}

 \noindent {\it Proof.} Only the latter part of (iii) still needs a proof. Since $\widetilde\Delta$ preserves the pure type of forms and since the kernel of $\widetilde\Delta:C^\infty_{p,\,q}(X,\,\C)\longrightarrow C^\infty_{p,\,q}(X,\,\C)$ is isomorphic to $E_2^{p,\,q}(X,\,\C)$ for every bidegree $(p,\,q)$ (cf. [Pop16, Theorem 1.1]), the isomorphism follows.  \hfill $\Box$

 \vspace{2ex}

 Thus, Conclusion \ref{Conc:pseudodiff_2} expresses the fact that $\Delta_h$ smoothly approximates $\Delta_0 = \Delta''$, while $\widetilde\Delta_h$ smoothly approximates $\widetilde\Delta_0 = \widetilde\Delta$.

To unify the notation throughout $\S$\ref{subsubsection:E_1_Delta-tilde_h_def}, $\S$\ref{subsubsection:E_2_Delta-tilde_h_def} and $\S$\ref{subsubsection:E_geq3_Delta-tilde_h-r_def}, we will also denote $\widetilde\Delta_h^{(2)}:=\widetilde\Delta_h$.

\subsubsection{{\bf Higher FSS pages: the pseudo-differential Laplacians $\widetilde\Delta^{(r)}_h$ with $r\geq 3$}}\label{subsubsection:E_geq3_Delta-tilde_h-r_def}

 In a similar fashion, one gets the following analogue for $r\geq 3$ of Conclusion \ref{Conc:pseudodiff_2}. Again, the properties of the operators $\Delta_h$ of $\S$\ref{subsubsection:E_1_Delta-tilde_h_def} are repeated for the sake of comparison with those of the operators $\widetilde\Delta^{(r)}_h$ introduced in this $\S$\ref{subsubsection:E_geq3_Delta-tilde_h-r_def}.

\begin{Conc}\label{Conc:pseudodiff_r} Let $(X,\,\omega)$ be a compact complex Hermitian manifold with $\mbox{dim}_\C X=n$. For every integer $r\geq 3$ and every degree $k\in\{0,\dots , 2n\}$, there exist $C^\infty$ families of elliptic differential operators $(\Delta_h)_{h\in\C}$ (independent of $r$) and, respectively, elliptic pseudo-differential operators $(\widetilde\Delta^{(r)}_h)_{h\in\C}$ from $C^\infty_k(X,\,\C)$ to $C^\infty_k(X,\,\C)$ such that

    \vspace{1ex}

   (i)\, $\Delta_0 = \Delta''$;

   \vspace{1ex}

   (ii)\, ${\mathcal H}^k_{\Delta_h}(X,\,\C) = {\mathcal H}^k_{\widetilde\Delta^{(r)}_h}(X,\,\C)\simeq H^k_{d_h}(X,\,\C)$ \hspace{3ex} for all $h\in\C\setminus\{0\}$;

   \vspace{1ex}

   (iii)\, ${\mathcal H}^k_{\Delta_0}(X,\,\C)\simeq\bigoplus_{p+q=k}H^{p,\,q}_{\bar\partial}(X,\,\C)$ \hspace{3ex} and \hspace{3ex} ${\mathcal H}^k_{\widetilde\Delta^{(r)}_0}(X,\,\C)\simeq\bigoplus_{p+q=k}E_r^{p,\,q}(X)$.

\end{Conc}

 Thus, $\widetilde\Delta^{(r)}_h$ smoothly approximates an elliptic operator $\widetilde\Delta^{(r)}_0$ whose kernel is isomorphic to $\oplus_{p+q=k}E_r^{p,\,q}(X)$.

 \vspace{2ex}
 
Conclusion \ref{Conc:pseudodiff_r} is needed only for our proof of the existence of the FAVB in $\S$\ref{subsubsection:absolute}. The proof of Conclusion \ref{Conc:pseudodiff_r} can be found in earlier arXiv versions of this paper. Parts of it are also to be found in [PSU20, $\S2.2$, $\S7$]. We skip these technical details due to the existence of the simpler algebraic approaches to the FAVB of [Aso06], [Ger66], [Kly89], [Ste20] mentioned in the introduction.

 \subsection{The Fr\"olicher approximating vector bundle (FAVB)}\label{subsection:approx_vb}

To fix the notation, recall the following statement that was proved in [CFGU97]. We will use the following terminology that was also used in [PSU20, Proposition 2.3].

\begin{Prop}\label{Prop:E_r-closed-exact_conditions} (i)\, Fix $r\geq 1$. A form $\alpha\in C^\infty_{p,\,q}(X,\,\C)$ is {\bf $E_r$-closed} (i.e. $\alpha$ represents an $E_r$-cohomology class) if and only if there exist forms $u_l\in C^\infty_{p+l,\,q-l}(X,\,\C)$ with $l\in\{1,\dots , r-1\}$ satisfying the following $r$ equations: \begin{eqnarray}\label{eqn:tower_E_r-closedness}\nonumber \bar\partial\alpha & = & 0 \\
    \nonumber \partial\alpha & = & \bar\partial u_1 \\
    \nonumber \partial u_1 & = & \bar\partial u_2 \\
    \nonumber \vdots & & \\
    \nonumber \partial u_{r-2} & = & \bar\partial u_{r-1}.\end{eqnarray}

  \vspace{1ex}

  (When $r=1$, the above equations reduce to $\bar\partial\alpha=0$.)

\vspace{1ex}

  (ii)\, Fix $r\geq 1$. A form $\alpha\in C^\infty_{p,\,q}(X,\,\C)$ is {\bf $E_r$-exact} (i.e. $\alpha$ represents the {\it zero} $E_r$-cohomology class) if and only if there exist forms $\zeta_{r-2}\in C^\infty_{p-1,\,q}(X,\,\C)$ and $\xi_0\in C^\infty_{p,\,q-1}(X,\,\C)$ such that $$\alpha=\partial\zeta_{r-2} + \bar\partial\xi_0,$$

\noindent with $\xi_0$ arbitrary and $\zeta_{r-2}$ satisfying the following additional condition (which is empty when $r=1$, denoting $\zeta_{-1}=0$, and reduces to requiring that $\zeta_{r-2} = \zeta_0$ be $\bar\partial$-closed when $r=2$.)

  \vspace{1ex}
  
  There exist $C^\infty$ forms $v^{(r-2)}_0, v^{(r-2)}_1,\dots , v^{(r-2)}_{r-3}$ satisfying the following $(r-1)$ equations: \begin{eqnarray}\label{eqn:tower_E_r-exactness_l}\nonumber \bar\partial\zeta_{r-2} & = & \partial v^{(r-2)}_{r-3}  \\
    \nonumber \bar\partial v^{(r-2)}_{r-3} & = & \partial v^{(r-2)}_{r-4} \\
    \nonumber \vdots & & \\
    \nonumber \bar\partial v^{(r-2)}_1 & = & \partial v^{(r-2)}_0 \\
               \bar\partial v^{(r-2)}_0 & = & 0,\end{eqnarray}
  
\noindent with the convention that any form $v^{(r-2)}_l$ with $l<0$ vanishes.

(Note that, thanks to (i), equations (\ref{eqn:tower_E_r-exactness_l}), when read from bottom to top, express precisely the condition that the form $v^{(r-2)}_0\in C^\infty_{p-r+1,\,q+r-2}(X,\,\C)$ be $E_{r-1}$-closed. Moreover, the form $\partial\zeta_{r-2}$ featuring on the r.h.s. of the above expression for $\alpha$ represents the $E_{r-1}$-class $(-1)^rd_{r-1}(\{v^{(r-2)}_0\}_{E_{r-1}})$.)

\vspace{1ex}

\end{Prop}

\noindent {\it Proof.} See [CFGU97].  \hfill $\Box$

\vspace{2ex}

We now return to the map $\theta_h:\Lambda^kT^\star X\longrightarrow\Lambda^kT^\star X,$ $\theta_h(\sum_{p+q=k}u^{p,\,q})=\sum_{p+q=k}h^p\,u^{p,\,q}$, with $h\in\C$ fixed.  When $h=0$, $\theta_0(\sum_{p+q=k}u^{p,\,q}) = u^{0,\,k}$. In the following statement, we notice that this projection onto the $(0,\,k)$-component at the level of forms induces the analogous projection in cohomology, that will still be denoted by $\theta_0$, relative to the (non-canonical) splitting $H^k_{DR}(X,\,\C)\simeq\oplus_{p+q=k}E^{p,\,q}_{\infty}(X)$ provided by the Fr\"olicher spectral sequence of $X$. Thus, unlike the splitting, the projection in cohomology is canonical.

\begin{Lem}\label{Lem:theta_0_cohom} For every $k\in\{0,\dots , n\}$, the canonical linear map: \begin{equation}\label{eqn:theta_0_cohom}\theta_0:H^k_{DR}(X,\,\C)\longrightarrow E^{0,\,k}_{\infty}(X), \hspace{3ex} \{\alpha\}_{DR}\longmapsto\{\alpha^{0,\,k}\}_{E_\infty} = \{\theta_0\alpha\}_{E_\infty},\end{equation} is {\bf well defined} and {\bf surjective}.

\end{Lem}

\noindent {\it Proof.} Let $r\in\N^\star$ be the smallest positive integer $l$ such that the Fr\"olicher spectral sequence of $X$ degenerates at $E_l$. In particular, $E^{0,\,k}_{\infty}(X) = E^{0,\,k}_r(X)$. 

To show well-definedness, we have to show two things, namely that

\vspace{1ex}

 (a)\, $\alpha^{0,\,k}=\theta_0\alpha$ is $E_r$-closed for every $d$-closed $k$-form $\alpha$. (This will show that $\alpha^{0,\,k}=\theta_0\alpha$ represents an $E_r$-cohomology class, or equivalently an $E_\infty$-cohomology class.)

\vspace{1ex}

 (b)\, for any De Rham cohomologous $d$-closed $k$-forms $\alpha$ and $\beta$, their $(0,\,k)$-components $\alpha^{0,\,k}$ and $\beta^{0,\,k}$ are $E_r$-cohomologous. (This will show that the $E_\infty$-cohomology class of $\alpha^{0,\,k}=\theta_0\alpha$ does not depend on the choice of representative of the De Rham class $ \{\alpha\}_{DR}$.)

\vspace{1ex}

To prove (a), let $\alpha\in C^\infty_k(X,\,\C)$ be $d$-closed. Identifying the pure-type components, we see that the condition $d\alpha = 0$ is equivalent to the following tower of $(k+2)$ equations: \begin{eqnarray}\label{eqn:tower_theta-map_proof}\nonumber\partial\alpha^{k,\,0} & = & 0 \\ \nonumber\partial\alpha^{k-1,\,1} & = & -\bar\partial\alpha^{k,\,0} \\ \nonumber & \vdots & \\ \nonumber\partial\alpha^{0,\,k} & = & -\bar\partial\alpha^{1,\,k-1} \\ \bar\partial\alpha^{0,\,k} & = & 0.\end{eqnarray} When read from bottom to top, this tower of equations implies that $\alpha^{0,\,k}$ is $E_l$-closed for every $l\geq k+2$. Indeed, the equation $\partial\alpha^{k,\,0} = 0 = \bar\partial(0)$ can be continued indefinitely with $\partial(0) = \bar\partial(0)$ repeated as many times as needed. (Note that $\partial\alpha^{k,\,0}$ is of type $(k+1,\,0)$, so it vanishes if and only if it is $\bar\partial$-exact.)

Now, if $k+2\geq r$, any $E_{k+2}$-closed form is also $E_r$-closed. So, $\alpha^{0,\,k}$ is $E_r$-closed in this case. If $k+2 < r$, we have already noticed above that $\alpha^{0,\,k}$ is $E_r$-closed. Thus, $\alpha^{0,\,k}$ is always $E_r$-closed.

\vspace{1ex}

To prove (b), let $\alpha, \beta\in C^\infty_k(X,\,\C)$ such that $d\alpha = d\beta = 0$ and $\alpha = \beta + d\gamma$ for some $\gamma\in C^\infty_{k-1}(X,\,\C)$. The last identity implies that $\alpha^{0,\,k} - \beta^{0,\,k} = \bar\partial\gamma^{0,\,k-1}$. Thus, being $\bar\partial$-exact (equivalently, $E_1$-exact), $\alpha^{0,\,k} - \beta^{0,\,k}$ is also $E_l$-exact for every $l\geq 1$, hence $E_r$-exact, i.e. $E_\infty$-exact. Therefore, $\{\alpha^{0,\,k}\}_{E_\infty} = \{\beta^{0,\,k}\}_{E_\infty}$.

\vspace{1ex}

To show surjectivity, let $\{\alpha^{0,\,k}\}_{E_r}\in E^{0,\,k}_r(X)$. Pick an arbitrary representative $\alpha^{0,\,k}\in C^\infty_{0,\,k}(X,\,\C)$ of this class. It is necessarily $E_r$-closed. Thus, if $r\geq k+2$, $\alpha^{0,\,k}$ is also $E_{k+2}$-closed. This means that there exist smooth pure-type forms $\alpha^{1,\,k-1}, \alpha^{2,\,k-2}, \dots , \alpha^{k-1,\,1}, \alpha^{k,\,0}$ of the shown types that, together with $\alpha^{0,\,k}$, satisfy the tower (\ref{eqn:tower_theta-map_proof}) of $(k+2)$ equations. This expresses the fact that the smooth $k$-form $\alpha: = \alpha^{k,\,0} + \dots + \alpha^{0,\,k}$ is $d$-closed. It is obvious, by construction, that $\theta_0(\{\alpha\}_{DR}) = \{\alpha^{0,\,k}\}_{E_\infty}$.

If $r\leq k+1$, then $E^{0,\,k}_r(X) = E^{0,\,k}_\infty(X) = E^{0,\,k}_{k+2}(X)$ and the $E_r$-closed forms coincide with the $E_{k+2}$-closed forms. Hence, we still get forms $\alpha^{l,\,k-l}$ as above satisfying the tower of equations (\ref{eqn:tower_theta-map_proof}) and the conclusion is the same. \hfill $\Box$

\subsubsection{The FAVB in the absolute case}\label{subsubsection:absolute}

As a first application of the pseudo-differential operators $\widetilde\Delta^{(r)}_h$ constructed for every $r\in\N^\star$ in $\S$\ref{subsubsection:E_1_Delta-tilde_h_def}, $\S$\ref{subsubsection:E_2_Delta-tilde_h_def} and $\S$\ref{subsubsection:E_geq3_Delta-tilde_h-r_def}, we associate with $X$ a holomorphic vector bundle over $\C$ (constructed algebraically in [Aso06], [Ger66], [Kly89], [Ste20]) whose fibre above $0$ is defined by the page in the Fr\"olicher spectral sequence of $X$ on which degeneration occurs. Thus, in the next statement, $E^{p,\,q}_r(X) = E^{p,\,q}_\infty(X)$ for all $p,q$.

 \begin{Cor-Def}\label{Cor-Def:bundle_fixed-manifold} Let $X$ be a compact complex manifold with $\mbox{dim}_\C X=n$. Let $r$ be the smallest positive integer such that the Fr\"olicher spectral sequence of $X$ degenerates at $E_r$. Fix a Hermitian metric $\omega$ on $X$.

\vspace{1ex}

  For every $k\in\{0,\dots , 2n\}$, the $C^\infty$ family $(\widetilde\Delta^{(r)}_h)_{h\in\C}$ of elliptic (pseudo-)differential operators constructed in $\S$\ref{subsubsection:E_1_Delta-tilde_h_def}, $\S$\ref{subsubsection:E_2_Delta-tilde_h_def} and $\S$\ref{subsubsection:E_geq3_Delta-tilde_h-r_def} induces a $C^\infty$ complex vector bundle ${\mathcal A}^k\longrightarrow\C$, of rank equal to the $k$-th Betti number $b_k$ of $X$, such that:

\vspace{1ex}

$\bullet$ its structure is described in Corollary \ref{Cor:K-S_conclusions_to-be-applied} of the Kodaira-Spencer theory after replacing $(E_t)_{t\in B}$ with $(\widetilde\Delta^{(r)}_h)_{h\in\C}$;

\vspace{1ex}

$\bullet$ its fibres are $${\mathcal A}^k_h=H^k_{d_h}(X,\,\C) \hspace{3ex} \mbox{if} \hspace{1ex} h\in\C\setminus\{0\} \hspace{3ex} \mbox{and} \hspace{3ex} {\mathcal A}^k_0=\bigoplus\limits_{p+q=k}E^{p,\,q}_r(X) \hspace{3ex} \mbox{if} \hspace{1ex} h=0;$$ 

\vspace{1ex}

$\bullet$ its restriction to $\C\setminus\{0\}$ is isomorphic to the constant vector bundle ${\mathcal H}^k_{|\C^\star}\longrightarrow\C\setminus\{0\}$ of fibre $H^k_{DR}(X,\,\C)$ under the holomorphic vector bundle isomorphism $\theta=(\theta_h)_{h\in\C^\star}:{\mathcal H}^k_{|\C^\star}\longrightarrow{\mathcal A}^k_{|\C^\star}$.

\vspace{1ex}

 The vector bundle ${\mathcal A}^k\longrightarrow\C$ has an extra structure as a {\bf holomorphic} vector bundle and is called the {\bf Fr\"olicher approximating vector bundle (FAVB)} of $X$ in degree $k$.

\end{Cor-Def}

 \noindent {\it Proof.} Recall that $\mbox{dim}_\C H^k_{d_h}(X,\,\C)=b_k$ for every $h\neq 0$. Fix any Hermitian metric $\omega$ on $X$.

 If $r=1$, the dimension of ${\mathcal A}^k_0=\oplus_{p+q=k}E^{p,\,q}_1(X,\,\C)$ equals $b_k$ and the $\C$-vector space ${\mathcal A}^k_0$ is isomorphic to the kernel of $\Delta'' = \Delta_0:C^\infty_k(X,\,\C)\longrightarrow C^\infty_k(X,\,\C)$. Thus, the $C^\infty$ family $(\Delta_h)_{h\in\C}$ of elliptic differential operators has the property that the dimension of the kernel of $\Delta_h:C^\infty_k(X,\,\C)\longrightarrow C^\infty_k(X,\,\C)$ is independent of $h\in\C$. Corollary \ref{Cor:K-S_conclusions_to-be-applied} of the classical Kodaira-Spencer theory [KS60] ensures that the harmonic spaces ${\mathcal H}^k_{\Delta_h}(X,\,\C)$ depend in a $C^\infty$ way on $h\in\C$ and that there exists a $C^\infty$ vector bundle on $\C$ whose fibres are these spaces. Equivalently, there exists a $C^\infty$ vector bundle ${\mathcal A}^k$ on $\C$ whose structure is described in Corollary \ref{Cor:K-S_conclusions_to-be-applied} and whose fibres are the vector spaces ${\mathcal A}^k_h$ to which the vector spaces ${\mathcal H}^k_{\Delta_h}(X,\,\C)$ are isomorphic.

 If $r=2$, the dimension of ${\mathcal A}^k_0=\oplus_{p+q=k}E^{p,\,q}_2(X,\,\C)$ equals $b_k$ and the fibre ${\mathcal A}^k_0$ is isomorphic to the kernel of $\widetilde\Delta = \widetilde\Delta_0:C^\infty_k(X,\,\C)\longrightarrow C^\infty_k(X,\,\C)$ by Theorem 1.1 in [Pop16]. Corollary \ref{Cor:K-S_conclusions_to-be-applied} of the classical Kodaira-Spencer theory [KS60] still applies to the $C^\infty$ family $(\widetilde\Delta_h)_{h\in\C}$ of elliptic pseudo-differential operators (cf. argument in [Mas18] for the case $h=0$), whose kernels have dimension independent of $h\in\C$ (and equal to $b_k$, see Conclusion \ref{Conc:pseudodiff_2}), to ensure that the harmonic spaces ${\mathcal H}^k_{\widetilde\Delta_h}(X,\,\C)$ depend in a $C^\infty$ way on $h\in\C$ and that there exists a $C^\infty$ vector bundle ${\mathcal A}^k$ on $\C$ whose structure is described in Corollary \ref{Cor:K-S_conclusions_to-be-applied} and whose fibres are the vector spaces ${\mathcal A}^k_h$ to which the harmonic spaces ${\mathcal H}^k_{\widetilde\Delta_h}(X,\,\C)$ are isomorphic for all $h\in\C$ (cf. Conclusion \ref{Conc:pseudodiff_2}).

 If $r\geq 3$, the dimension of ${\mathcal A}^k_0=\oplus_{p+q=k}E^{p,\,q}_r(X,\,\C)$ equals $b_k$ and the fibre ${\mathcal A}^k_0$ is isomorphic to the kernel of $\widetilde\Delta^{(r)}_0:C^\infty_k(X,\,\C)\to C^\infty_k(X,\,\C)$ (cf. Conclusion \ref{Conc:pseudodiff_r}). Corollary \ref{Cor:K-S_conclusions_to-be-applied} of the classical Kodaira-Spencer theory [KS60] still applies to the $C^\infty$ family $(\widetilde\Delta^{(r)}_h)_{h\in\C}$ of elliptic pseudo-differential operators (cf. argument in [Mas18] for the case of $\widetilde\Delta$), whose kernels have dimension independent of $h\in\C$ (and equal to $b_k$), to ensure that the harmonic spaces ${\mathcal H}^k_{\widetilde\Delta^{(r)}_h}(X,\,\C)$ depend in a $C^\infty$ way on $h\in\C$ and that there exists a $C^\infty$ vector bundle ${\mathcal A}^k$ on $\C$ whose structure is described in Corollary \ref{Cor:K-S_conclusions_to-be-applied} and whose fibres are the vector spaces ${\mathcal A}^k_h$ to which the harmonic spaces ${\mathcal H}^k_{\widetilde\Delta^{(r)}_h}(X,\,\C)$ are isomorphic for all $h\in\C$ (cf. Conclusion \ref{Conc:pseudodiff_r}). 

\vspace{1ex}

 Meanwhile, we know from [Pop17, Lemma 2.5] (see also Introduction) that for every $h\neq 0$, the linear map $\theta_h: H^k_{DR}(X,\,\C)\longrightarrow H^k_{d_h}(X,\,\C)$ defined by $\theta_h(\{u\}_{DR}) = \{\theta_h u\}_{d_h}$ is an isomorphism of $\C$-vector spaces. Since $\theta_h$ depends {\it holomorphically} on $h$ and the space $H^k_{DR}(X,\,\C)$ is independent of $h$, we infer that the restriction to $\C\setminus\{0\}$ of the $C^\infty$ vector bundle ${\mathcal A}^k\longrightarrow\C$ constructed above has an extra holomorphic structure obtained as the image of the holomorphic vector bundle structure of the constant bundle ${\mathcal H}^k_{|\C^\star}\longrightarrow\C\setminus\{0\}$ under the vector bundle isomorphism $\theta=(\theta_h)_{h\in\C^\star}$. 

This holomorphic structure of ${\mathcal A}^k$ over $\C\setminus\{0\}$ extends holomorphically across $0$ since the underlying $C^\infty$ vector bundle structure does, as we have seen above. Indeed, (ii) of Corollary \ref{Cor:K-S_conclusions_to-be-applied} yields a global $C^\infty$ trivialisation of the $C^\infty$ vector bundle ${\mathcal A}^k\longrightarrow\C$ in the form of a family of bases \begin{eqnarray}\label{eqn:absolute_family_bases}\bigg(\bigg\{\varphi_h^{(1)},\dots , \varphi_h^{(b_k)}\bigg\}\bigg)_{h\in\C}\end{eqnarray} of the respective $\C$-vector spaces $({\mathcal A}^k_h)_{h\in\C}$. Specifically, the $\varphi_h^{(q)}$'s are the $L^2_\omega$-orthogonal projections onto $\ker\widetilde\Delta^{(r)}_h\simeq{\mathcal A}^k_h$ (i.e. the images under the $L^2_\omega$-orthogonal projection operator $H_h^{(k)}:C^\infty_k(X,\,\C)\longrightarrow\ker\widetilde\Delta^{(r)}_h\subset C^\infty_k(X,\,\C)$) of a fixed $\C$-basis $(s_j^{(k)})_{j\in\N}$ of $C^\infty_k(X,\,\C)$. The family of operators $(H_h^{(k)})_{h\in\C}$ being $C^\infty$ (up to $h=0$) -- see the Kodaira-Spencer Theorem \ref{The:K-S_E_t_bundle} -- we infer that, for every $q$, the dependence of $\varphi_h^{(q)}$ on $h\in\C$ is $C^\infty$ on the whole of $\C$. Meanwhile, this dependence on $h\in\C^\star$ is also holomorphic on $\C^\star$ thanks to the holomorphic dependence on $h\in\C^\star$ of $\theta_h$. It must then be holomorphic on the whole of $\C$. \hfill $\Box$

\subsubsection{The FAVB in the relative case}\label{subsubsection:relative}

We will now define the Fr\"olicher approximating vector bundles of a holomorphic family $(X_t)_{t\in B}$ of compact complex $n$-dimensional manifolds induced by a proper holomorphic submersion $\pi:{\mathcal X}\longrightarrow B$ whose base $B\subset\C^N$ is an open ball about the origin in some complex Euclidean vector space. 

By the classical Ehresmann Theorem, the differential structure of the fibres $X_t$ is independent of $t\in B$, hence so is the Poincar\'e differential $d$. However, its splitting as $d=\partial_t + \bar\partial_t$ depends on the complex structure of $X_t$ that varies with $t$. Thus, the differential operators $d_h$ depend on $t$ (except when $h=1$), so we put $$d_{h,t}:=h\partial_t + \bar\partial_t:C^\infty_k(X,\,\C)\longrightarrow C^\infty_{k+1}(X,\,\C),  \hspace{3ex} h\in\C, \hspace{1ex} t\in B, \hspace{1ex} k\in\{0,\dots , 2n\},$$

\noindent where $X$ is the $C^\infty$ manifold underlying the fibres $X_t$. 

Consequently, if $(\omega_t)_{t\in B}$ is a $C^\infty$ family of Hermitian metrics on the fibres $(X_t)_{t\in B}$, the $1$-parameter families of operators constructed in $\S$\ref{subsubsection:E_1_Delta-tilde_h_def}, $\S$\ref{subsubsection:E_2_Delta-tilde_h_def} and $\S$\ref{subsubsection:E_geq3_Delta-tilde_h-r_def} depend now on the extra variable $t$, so they become, for every $r\in\N^\star$, a (possibly non-continuous) $2$-parameter family $(\widetilde\Delta_{h,\,t}^{(r)})_{(h,\,t)\in\C\times B}$ of elliptic (pseudo-)differential operators. As in $\S$\ref{subsection:Laplacians_approximations}, we will sometimes use the unifying notation $$\Delta_{h,\,t}:=\widetilde\Delta_{h,\,t}^{(1)} \hspace{3ex}\mbox{and}\hspace{3ex} \widetilde\Delta_{h,\,t}:=\widetilde\Delta_{h,\,t}^{(2)}.$$ 

For example, when $r=1$, we have a $C^\infty$ family $(\Delta_{h,\,t})_{(h,\,t)\in\C\times B}$ of elliptic differential operators defined in every degree $k$ by analogy with the absolute case: $$\Delta_{h,\,t}=d_{h,\,t}d_{h,\,t}^\star + d_{h,\,t}^\star d_{h,\,t}:C^\infty_k(X,\,\C)\longrightarrow C^\infty_k(X,\,\C),$$ 

\noindent where the formal adjoint $d_{h,\,t}^\star$ is computed w.r.t. the metric $\omega_t$. 

\vspace{1ex}

A key fact that we will exploit is the following

\begin{Cor}\label{Cor:Laplacians_same-kernel} Let $(\omega_t)_{t\in B}$ be a $C^\infty$ family of Hermitian metrics on the fibres $(X_t)_{t\in B}$ of a holomorphic family of compact complex $n$-dimensional manifolds. 

Then, for every $k\in\{0,\dots , 2n\}$ and every positive integer $r$, the Laplacians $\Delta_{h,\,t}:C^\infty_k(X,\,\C)\longrightarrow C^\infty_k(X,\,\C)$ and $\widetilde\Delta_{h,\,t}^{(r)}:C^\infty_k(X,\,\C)\longrightarrow C^\infty_k(X,\,\C)$ have the same kernel: \begin{eqnarray}\label{eqn:Laplacians_same-kernel}\ker\Delta_{h,\,t} = \ker\widetilde\Delta_{h,\,t}^{(r)}\end{eqnarray} for every $(h,\,t)\in\C^\star\times B$.

\end{Cor}

\noindent {\it Proof.} This is part of (ii) of Conclusions \ref{Conc:pseudodiff_2} and \ref{Conc:pseudodiff_r} when $(X,\,\omega)$ is replaced by $(X_t,\,\omega_t)$. \hfill $\Box$

\vspace{2ex}

Likewise, the pointwise linear maps $\theta_h$ (which are isomorphisms when $h\neq 0$) depend on $t$ (because the splitting of $k$-forms into pure-type-forms depends on the complex structure of $X_t$), so we put $$\theta_{h,t}:\Lambda^k T^\star X\longrightarrow\Lambda^k T^\star X, \hspace{3ex} u=\sum\limits_{p+q=k}u^{p,\,q}_t\mapsto\theta_{h,t}u:=\sum\limits_{p+q=k}h^p\,u^{p,\,q}_t.$$

 When $h\neq 0$, this induces an {\it isomorphism} in cohomology: \begin{equation}\label{eqn:theta_ht_rel-case}\theta_{h,t}: H^k_{DR}(X,\,\C)\longrightarrow H^k_{d_{h,t}}(X_t,\,\C), \hspace{6ex} \theta_{h,t}(\{u\}_{DR}) = \{\theta_{h,t}u\}_{d_{h,t}},\end{equation} for every $t\in B$, since $\theta_{h,t}d = d_{h,t}\theta_{h,t}$. When $h=0$, we saw in Lemma \ref{Lem:theta_0_cohom} that $\theta_{0,t}$ induces a {\it surjective} linear map: \begin{equation}\label{eqn:theta_h0_rel-case}\theta_{0,t}: H^k_{DR}(X,\,\C)\longrightarrow E^{0,\,k}_\infty(X_t), \hspace{6ex} \theta_{0,t}(\{u\}_{DR}) = \{u^{0,\,k}_t\}_{E_\infty},\end{equation} for every $t\in B$, where $u^{0,\,k}_t$ is the component of type $(0,\,k)$ of $u$ w.r.t. the complex structure of $X_t$.

Recall that the degeneration at $E_1$ of the Fr\"olicher spectral sequence is a {\it deformation open} property of compact complex manifolds. Thus, if $E_1(X_0)=E_\infty(X_0)$, then $E_1(X_t)=E_\infty(X_t)$ for every $t\in B$, after possibly shrinking $B$ about $0$. (This follows at once from the upper semicontinuity of the Hodge numbers $h^{p,\,q}(t)$ and from the invariance of the Betti numbers $b_k$ of the fibres $X_t$.) 

However, when $r\geq 2$, the degeneration at $E_r$ of the Fr\"olicher spectral sequence is not deformation open, so we will have to make an appropriate assumption with regard to this in Corollary and Definition \ref{Cor-Def:bundle_family}.

Note also that, for every $t\in B$, the fibre at $(0,\,t)$ of the vector bundle ${\mathcal A}^k\longrightarrow\C\times B$ constructed therein consists of the direct sum of all the $E^{p,\,q}_\infty(X_t)$ with $p+q=k$, which is different from the space $F^pH^k(X_t)$ for a fixed $p$ of the standard fibrewise Hodge filtration. 

\vspace{2ex}

The next statement is a further illustration of the Kodaira-Spencer theory of [KS60] and [Kod86] whereby $C^\infty$ families of elliptic (pseudo-)differential operators equip the families of their kernels with $C^\infty$ vector bundle structures.

\begin{Cor-Def}\label{Cor-Def:bundle_family} Fix $N\in\N^\star$. Let $\pi:{\mathcal X}\longrightarrow B$ be a holomorphic family of compact complex $n$-folds over an open ball $B\subset\C^N$ about the origin. Let $r$ be the smallest positive integer such that the Fr\"olicher spectral sequence of $X_0$ degenerates at $E_r$. Suppose that $E_r(X_t) = E_\infty(X_t)$ for all $t\in B$.

Let $(\omega_t)_{t\in B}$ be a $C^\infty$ family of Hermitian metrics on the fibres $(X_t)_{t\in B}$.

\vspace{1ex}

 (i)\, For every $k\in\{0,\dots , 2n\}$, the $C^\infty$ families:

\vspace{1ex}

$\bullet$ $(\Delta_{h,\,t})_{(h,\,t)\in\C^\star\times B}$ 

\vspace{1ex}

\noindent and 

\vspace{1ex}

$\bullet$ for every fixed $t\in B$, $(\widetilde\Delta_{h,\,t}^{(r)})_{h\in\C}$ 

\vspace{1ex}

\noindent of elliptic (pseudo-)differential operators constructed in $\S$\ref{subsubsection:E_1_Delta-tilde_h_def}, $\S$\ref{subsubsection:E_2_Delta-tilde_h_def} and $\S$\ref{subsubsection:E_geq3_Delta-tilde_h-r_def} induce a $C^\infty$ complex vector bundle ${\mathcal A}^k\longrightarrow\C\times B$ of rank equal to the $k$-th Betti number $b_k$ of $X$ ($=$ of any fibre $X_t$) such that:

\vspace{1ex}

(a)\, its structure is described in Corollary \ref{Cor:K-S_conclusions_to-be-applied} of the Kodaira-Spencer theory after replacing $(E_t)_{t\in B}$ with the above families of operators;

\vspace{1ex}

(b)\, its fibres are $${\mathcal A}^k_{h,\,t}=H^k_{d_{h, t}}(X_t,\,\C) \hspace{2ex} \mbox{for} \hspace{1ex} (h,\,t)\in\C^\star\times B, \hspace{4ex}  {\mathcal A}^k_{0,\,t}=\bigoplus\limits_{p+q=k}E^{p,\,q}_r(X_t) \hspace{2ex} \mbox{for} \hspace{1ex} (0,\,t)\in\{0\}\times B;$$ 

(c)\, its restriction to $\C^\star\times B$ is isomorphic to the restriction of the constant vector bundle ${\mathcal H}^k\longrightarrow\C\times B$ of fibre $H^k_{DR}(X,\,\C)$ under the holomorphic vector bundle isomorphism $\theta=(\theta_{h,t})_{(h,\,t)\in\C^\star\times B}:{\mathcal H}^k_{|\C^\star\times B}\longrightarrow{\mathcal A}^k_{|\C^\star\times B}$.

\vspace{1ex}

 (ii)\, For every $k\in\{0,\dots , 2n\}$ and every $C^\infty$ family $s=(s_{h,t})_{(h,\,t)\in\C\times B}$ of $C^\infty$ $k$-forms on the smooth manifold $X$ underlying the fibres $X_t$ such that \begin{eqnarray*}d_{h,t}(s_{h,t}) = 0, & &  \mbox{for all}\hspace{1ex} (h,\,t)\in\C^\star\times B; \\
  s_{0,\,t} \hspace{1ex} \mbox{is} \hspace{1ex} E_r(X_t)\mbox{-closed,} & & \mbox{for all}\hspace{1ex} t\in B,\end{eqnarray*} the family $\sigma=(\sigma_{h,t})_{(h,\,t)\in\C\times B}$ defined by \begin{eqnarray*}\sigma_{h,t} & = & \{s_{h,t}\}_{d_{h,t}}\in H^k_{d_{h,\,t}}(X_t,\,\C)={\mathcal A}^k_{h,\,t}, \hspace{3ex} \mbox{if} \hspace{3ex} (h,\,t)\in\C^\star\times B,  \\  
 \sigma_{0,t} & = & \{s_{0,t}\}_{E_r(X_t)}\in\bigoplus\limits_{p+q=k}E_r^{p,\,q}(X_t) = {\mathcal A}^k_{0,\,t}, \hspace{3ex} \mbox{if} \hspace{3ex} (h,\,t) = (0,\,t)\in\{0\}\times B,\end{eqnarray*} defines a $C^\infty$ section of the vector bundle ${\mathcal A}^k$ over $\C\times B$.

\vspace{1ex}

 (iii)\, The $C^\infty$ vector bundle ${\mathcal A}^k\longrightarrow\C\times B$ has an extra structure as a {\bf holomorphic} vector bundle and is called the {\bf Fr\"olicher approximating vector bundle (FAVB)} of the family $(X_t)_{t\in B}$ in degree $k$.

\end{Cor-Def}

\noindent {\it Proof.} {\bf (i) and (iii).} We know that $\mbox{dim}_\C H^k_{d_{h, t}}(X_t,\,\C)=b_k$ for all $h\neq 0$ and $t\in B$. Moreover, thanks to the $E_r$-degeneration assumption on every fibre $X_t$, $\mbox{dim}_\C\oplus_{p+q=k}E^{p,\,q}_r(X_t) = b_k$ for all $t\in B$. Thus, $\mbox{dim}_\C{\mathcal A}^k_{h,\,t}=b_k$ for all $(h,\,t)\in\C\times B$. 

%Now, the given $C^\infty$ family $(\omega_t)_{t\in B}$ of Hermitian metrics on the fibres $(X_t)_{t\in B}$ induces a $C^\infty$ family $(\Delta_{h,\,t})_{(h,\,t)\in\C^\star\times B}$ of elliptic differential operators defined in every degree $k$ by analogy with the absolute case as $$\Delta_{h,\,t}=d_{h,\,t}d_{h,\,t}^\star + d_{h,\,t}^\star d_{h,\,t}:C^\infty_k(X,\,\C)\longrightarrow C^\infty_k(X,\,\C),$$ 

%\noindent where the formal adjoint $d_{h,\,t}^\star$ is computed w.r.t. the metric $\omega_t$. 

For every $k$, we fix a $\C$-basis $(s_j^{(k)})_{j\in\N}$ of the infinite-dimensional $\C$-vector space $C^\infty_k(X,\,\C)$ and, for every $(h,\,t)\in\C\times B$, we let \begin{eqnarray*}H^{(k)}_{h,\,t}:C^\infty_k(X,\,\C)\longrightarrow{\mathcal H}^k_{\widetilde\Delta^{(r)}_{h,\,t}}(X_t,\,\C)\subset C^\infty_k(X,\,\C)\end{eqnarray*} be the $L^2_{\omega_t}$-orthogonal projection onto $\ker\widetilde\Delta_{h,\,t}^{(r)}$ in degree $k$. Recall that, by Corollary \ref{Cor:Laplacians_same-kernel}, we have $\ker\widetilde\Delta_{h,\,t}^{(r)} = \ker\Delta_{h,\,t}$ for every $(h,\,t)\in\C^\star\times B$.

As explained in (II) of $\S$\ref{subsubsection:FAVB_introd}, the operator $\widetilde\Delta_{h,\,t}^{(r)}$ need not depend continuously on $t\in B$, but the families (\ref{eqn:smooth-families-Laplacians-ht_introd}) of operators are $C^\infty$. Hence, due to the equidimensionality of the kernels of these operators, the Kodaira-Spencer theory (see Theorem \ref{The:K-S_E_t_bundle}) ensures that the families \begin{eqnarray}\label{eqn:smooth-families-projections-ht_proof}\nonumber & & \bullet\, (H_{h,\,t}^{(k)})_{(h,\,t)\in\C^\star\times B} \\
\mbox{and} & & \\
\nonumber & & \bullet\, \mbox{for every fixed}\hspace{1ex} t\in B, (H_{h,\,t}^{(k)})_{h\in\C}\end{eqnarray} of orthogonal projectors are $C^\infty$.  

\vspace{1ex}

$\bullet$ By (ii) of Conclusion \ref{Conc:pseudodiff_2} with $\Delta_h$ replaced by $\Delta_{h,\,t}$, the kernels $\ker\Delta_{h,\,t}$ are isomorphic to the vector spaces ${\mathcal A}^k_{h,\,t}$, hence they have a dimension independent of $(h,\,t)\in\C^\star\times B$ (and equal to $b_k$). This implies, thanks to the classical Kodaira-Spencer theory of [KS60], that there exists a $C^\infty$ complex vector bundle \begin{eqnarray}\label{eqn:relative-FAVB_proof_1}{\mathcal A}^k\longrightarrow\C^\star\times B\end{eqnarray} of rank $b_k$ whose structure is defined by the $C^\infty$ family of elliptic differential operators $(\Delta_{h,\,t})_{(h,\,t)\in\C^\star\times B}$ as described in Corollary \ref{Cor:K-S_conclusions_to-be-applied} after replacing $(E_t)_{t\in B}$ with $(\Delta_{h,\,t})_{(h,\,t)\in\C^\star\times B}$. The fibres of this vector bundle are the vector spaces ${\mathcal A}^k_{h,\,t}$ to which the harmonic spaces $\ker\Delta_{h,\,t}$ are respectively isomorphic for all $(h,\,t)\in\C^\star\times B$. 

The $C^\infty$ vector bundle (\ref{eqn:relative-FAVB_proof_1}) has an extra {\it holomorphic} vector bundle structure obtained by transporting the holomorphic vector bundle structure of ${\mathcal H}^k_{|\C^\star\times B}$ under the vector bundle isomorphism $$\theta=(\theta_{h,t})_{(h,\,t)\in\C^\star\times B}:{\mathcal H}^k_{|\C^\star\times B}\longrightarrow{\mathcal A}^k_{|\C^\star\times B}$$ viewed as a section of $\mbox{End}\,({\mathcal H}^k_{|\C^\star\times B},\,{\mathcal A}^k_{|\C^\star\times B})$. Indeed, $\theta_{h,t}$ depends in a holomorphic way on $(h,\,t)\in\C^\star\times B$.

 From this and from (ii) of Corollary \ref{Cor:K-S_conclusions_to-be-applied} we infer that a global trivialisation of the holomorphic vector bundle ${\mathcal A}^k\longrightarrow\C^\star\times B$ is defined by the family of bases \begin{eqnarray}\label{eqn:families_bases_1}\bigg(\bigg\{\varphi_{h,\,t}^{(1)},\dots , \varphi_{h,\,t}^{(b_k)}\bigg\}\bigg)_{(h,\,t)\in\C^\star\times B}\end{eqnarray} of the respective $\C$-vector spaces $({\mathcal A}^k_{h,\,t})_{(h,\,t)\in\C^\star\times B}$, where the $\varphi_{h,\,t}^{(q)}$'s are the $L^2_{\omega_t}$-orthogonal projections onto $\ker\Delta_{h,\,t}\simeq{\mathcal A}^k_{h,\,t}$ (i.e. the images under the operator $H_{h,\,t}^{(k)}$) of the fixed $\C$-basis $(s_j^{(k)})_{j\in\N}$ of $C^\infty_k(X,\,\C)$. Thus, for every $q=1,\dots , b_k$, the $\C$-valued $k$-form $\varphi_{h,\,t}^{(q)}\in\ker\Delta_{h,\,t}$ on $X$ depends holomorphically on $(h,\,t)\in\C^\star\times B$. 

\vspace{1ex}

Note that no assumption on the spectral sequence is necessary to construct the holomorphic vector bundle ${\mathcal A}^k\longrightarrow\C^\star\times B$ of (\ref{eqn:relative-FAVB_proof_1}). As pointed out in General Fact \ref{Gen-Fact:bundle-trivialty_C-star-B}, this vector bundle is necessarily trivial, but our construction has produced the particular trivialisation defined by the global frame (\ref{eqn:families_bases_1}) associated with the given family $(\omega_t)_{t\in B}$ of metrics.  

\vspace{1ex}

$\bullet$ On the other hand, for every fixed $t\in B$, we know from the absolute case of Corollary and Definition \ref{Cor-Def:bundle_fixed-manifold} that $\C\ni h\longmapsto{\mathcal A}^k_{h,\,t}$ is a {\it holomorphic} vector bundle (of rank $b_k$) over $\C$ whose $C^\infty$ vector bundle structure is defined, via the classical Kodaira-Spencer theory of [KS60], by the $C^\infty$ family $(\widetilde\Delta_{h,\,t}^{(r)})_{h\in\C}$ of elliptic (pseudo-)differential operators with equidimensional kernels as described in Corollary \ref{Cor:K-S_conclusions_to-be-applied} after replacing $(E_t)_{t\in B}$ with $(\widetilde\Delta_{h,\,t}^{(r)})_{h\in\C}$. 

From this and from (ii) of Corollary \ref{Cor:K-S_conclusions_to-be-applied} we infer that, for every fixed $t\in B$, a global trivialisation of the holomorphic vector bundle \begin{eqnarray}\label{eqn:relative-FAVB_proof_2}{\mathcal A}^k_{\bullet,\,t}\longrightarrow\C\end{eqnarray} is defined by the family of bases \begin{eqnarray}\label{eqn:families_bases_2}\bigg(\bigg\{\varphi_{h,\,t}^{(1)},\dots , \varphi_{h,\,t}^{(b_k)}\bigg\}\bigg)_{h\in\C}\end{eqnarray} of the respective $\C$-vector spaces $({\mathcal A}^k_{h,\,t})_{h\in\C}$, where the $\varphi_{h,\,t}^{(q)}$'s are the $L^2_{\omega_t}$-orthogonal projections onto $\ker\widetilde\Delta_{h,\,t}^{(r)}\simeq{\mathcal A}^k_{h,\,t}$ (i.e. the images under the operator $H_{h,\,t}^{(k)}$) of the fixed $\C$-basis $(s_j^{(k)})_{j\in\N}$ of $C^\infty_k(X,\,\C)$. Thus, for every $q=1,\dots , b_k$, the $\C$-valued $k$-form $\varphi_{h,\,t}^{(q)}\in\ker\widetilde\Delta_{h,\,t}^{(r)}$ on $X$ depends holomorphically on $h\in\C$. 

An argument similar to the one that gave General Fact \ref{Gen-Fact:bundle-trivialty_C-B} shows that its conclusion remains valid when $\C\times B$ is replaced by $\C$. Thus, the holomorphic vector bundle ${\mathcal A}^k_{\bullet,\,t}\longrightarrow\C$ of (\ref{eqn:relative-FAVB_proof_2}) is necessarily trivial for every fixed $t\in B$. However, our construction has produced the particular trivialisation defined by the global frame (\ref{eqn:families_bases_2}) associated with the given family $(\omega_t)_{t\in B}$ of metrics.

\vspace{1ex}

$\bullet$ All that now remains to do is to put together the conclusions we drew above about the holomorphic vector bundle (\ref{eqn:relative-FAVB_proof_1}) and the family, indexed over $t\in B$, of holomorphic vector bundles (\ref{eqn:relative-FAVB_proof_2}).

For every $(h,\,t)\in\C^\star\times B$, since $\ker\widetilde\Delta_{h,\,t}^{(r)} = \ker\Delta_{h,\,t}$ (see Corollary \ref{Cor:Laplacians_same-kernel}), the $\C$-valued $k$-forms $\varphi_{h,\,t}^{(q)}$ of (\ref{eqn:families_bases_1}) coincide with the $\C$-valued $k$-forms $\varphi_{h,\,t}^{(q)}$ of (\ref{eqn:families_bases_2}). Indeed, the two families of forms are the images of the same $\C$-basis $(s_j^{(k)})_{j\in\N}$ of $C^\infty_k(X,\,\C)$ under the same family $H_{h,\,t}^{(k)}$ of $L^2_{\omega_t}$-orthogonal projectors.

Now, the forms $\varphi_{h,\,t}^{(q)}$ in (\ref{eqn:families_bases_1}) are defined, a priori, for $(h,\,t)\in\C^\star\times B$. However, the forms $\varphi_{h,\,t}^{(q)}$ in (\ref{eqn:families_bases_2}) are defined, for every fixed $t\in B$, for all $h\in\C$ (including $h=0$). Since these two families of forms coincide for every $(h,\,t)\in\C^\star\times B$, we conclude that we actually have a family of bases \begin{eqnarray}\label{eqn:families_bases_3}\bigg(\bigg\{\varphi_{h,\,t}^{(1)},\dots , \varphi_{h,\,t}^{(b_k)}\bigg\}\bigg)_{(h,\,t)\in\C\times B}\end{eqnarray} of the respective $\C$-vector spaces $({\mathcal A}^k_{h,\,t})_{(h,\,t)\in\C\times B}$. 

Moreover, for every $q=1,\dots , b_k$, the form $\varphi_{h,\,t}^{(q)}$ depends in a holomorphic way on $(h,\,t)\in\C\times B$ because the maps $\C^\star\times B\ni(h,\,t)\mapsto\varphi_{h,\,t}^{(q)}$ and, for every $t\in B$, $\C\ni h\mapsto\varphi_{h,\,t}^{(q)}$ are holomorphic.

We conclude that (\ref{eqn:families_bases_3}) is a global holomorphic frame of a (necessarily trivial) holomorphic vector bundle \begin{eqnarray*}\C\times B\ni(h,\,t)\longmapsto{\mathcal A}^k_{h,\,t}\end{eqnarray*} with the properties described in the statement under (i) and (iii).  

\vspace{2ex}

{\bf (ii).} This statement is an immediate application of (iii) of Corollary \ref{Cor:K-S_conclusions_to-be-applied} combined with the Hodge isomorphisms of Conclusions \ref{Conc:pseudodiff_2} and \ref{Conc:pseudodiff_r}.  \hfill $\Box$

\section{$E_r$-sG manifolds and deformations of complex structures}\label{section:E_r-sG}

In this section, we apply the Fr\"olicher approximating vector bundle (FAVB) constructed in $\S$\ref{subsection:approx_vb} to the study of limits of $\partial\bar\partial$-manifolds and then of real $(1,\,1)$-cohomology classes under holomorphic deformations. Together with a uniform volume control result, this will prove the main theorem \ref{The:moidef}.

\subsection{$E_r$-sG manifolds}\label{subsection:E_r-sG}

We begin by generalising the notion of {\it strongly Gauduchon (sG)} metric introduced in [Pop09] and [Pop13]. Recall that a Gauduchon metric on a compact complex $n$-dimensional manifold $X$ is a positive definite, $C^\infty$ $(1,\,1)$-form $\gamma$ on $X$ such that $\partial\bar\partial\gamma^{n-1}=0$ (or, equivalently, $\partial\gamma^{n-1}$ is $\bar\partial$-closed). Thanks to [Gau77], such metrics always exist. If the stronger requirement that $\partial\gamma^{n-1}$ be $\bar\partial$-exact ($= E_1$-exact w.r.t. the Fr\"olicher spectral sequence) is imposed, $\gamma$ is said to be {\it strongly Gauduchon (sG)} (cf. [Pop09] and [Pop13]). We will relax this definition by requiring $E_r$-exactness instead, for a possibly larger $r\geq 1$.

Finally, note that for any Gauduchon metric $\gamma$ on $X$, the $(n,\,n-1)$-form $\partial\gamma^{n-1}$ is {\it $E_r$-closed} for every $r\in\N^\star$. Indeed, in (i) of Proposition \ref{Prop:E_r-closed-exact_conditions} we can choose $u_1= \dots = u_{r-1} = 0$.

\begin{Def}\label{Def:E_r-sG} Let $\gamma$ be a Gauduchon metric on a compact complex manifold $X$ with $\mbox{dim}_\C X=n$. Fix an arbitrary integer $r\geq 1$.

  \vspace{1ex}
  
 (i)\, We say that $\gamma$ is an {\bf $E_r$-sG metric} if  $\partial\gamma^{n-1}$ is $E_r$-exact.

  (ii)\, A compact complex manifold $X$ is said to be an {\bf $E_r$-sG manifold} if an $E_r$-sG metric exists on $X$.

  (iii)\, A compact complex manifold $X$ is said to be an {\bf $E_r$-sGG manifold} if every Gauduchon metric on $X$ is an $E_r$-sG metric.

\end{Def}  

\vspace{2ex}

The term chosen in the last definition is a nod to the notion of sGG manifold that we introduced jointly with L. Ugarte in [PU14] as any compact complex manifold on which every Gauduchon metric is strongly Gauduchon. It follows from the above definitions that the $E_1$-sG property is equivalent to the sG property and that the following implications hold for any Hermitian metric $\gamma$ and every $r\in\N^\star$: $$\gamma \hspace{1ex} \mbox{is} \hspace{1ex} E_1\mbox{-sG} \implies \gamma \hspace{1ex} \mbox{is} \hspace{1ex} E_2\mbox{-sG} \implies\dots\implies \gamma \hspace{1ex} \mbox{is} \hspace{1ex} E_r\mbox{-sG} \implies \gamma \hspace{1ex} \mbox{is} \hspace{1ex} E_{r+1}\mbox{-sG} \implies\dots.$$

 Actually, for bidegree reasons, if a Hermitian metric $\gamma$ is $E_r$-sG for some integer $r\geq 1$, only the case $r\leq 3$ is relevant. Indeed, if $(p,\,q) = (n,\,n-1)$, the tower of relations (\ref{eqn:tower_E_r-exactness_l}) reduces to its first two lines since $\zeta_{r-2}$ is of bidegree $(n-1,\,n-1)$, hence $v^{(r-2)}_{r-3}$ is of bidegree $(n-2,\,n)$, hence $\bar\partial v^{(r-2)}_{r-3}=0$ for bidegree reasons, so $v^{(r-2)}_{r-4}, \dots, v^{(r-2)}_0$ can all be chosen to be zero. 

\vspace{2ex}

We now notice that the $E_r$-sG property is open under deformations of the complex structure. 

\begin{Lem}\label{Lem:E_r-sG_openness} Let $\pi:{\mathcal X}\longrightarrow B$ be a $C^\infty$ family of compact complex $n$-dimensional manifolds over an open ball $B\subset\C^N$ about the origin. Fix an integer $r\geq 1$.

 If $\gamma_0$ is an $E_r$-sG metric on $X_0:=\pi^{-1}(0)$, after possibly shrinking $B$ about $0$ there exists a $C^\infty$ family $(\gamma_t)_{t\in B}$ of $E_r$-sG metrics on the respective fibres $X_t:=\pi^{-1}(t)$ whose element for $t=0$ is the original $\gamma_0$. 

Moreover, this family can be chosen such that $\partial_t\gamma_t^{n-1} = \bar\partial_t\Gamma_t^{n,\,n-2} + \partial_t\zeta_{r-2,\,t}$ for all $t$, with $J_t$-type $(n,\,n-2)$-forms $\Gamma_t^{n,\,n-2}$ and $J_t$-type $(n-1,\,n-1)$-forms $\zeta_{r-2,\,t}$ depending in a $C^\infty$ way on $t$.

\end{Lem}

The forms $\Gamma_t^{n,\,n-2}$, $\zeta_{r-2,\,t}$ and the induced $v^{(r-2)}_{k,\,t}$ (with $0\leq k\leq r-3$) satisfying the tower of relations (\ref{eqn:tower_E_r-exactness_l}) that are (non-uniquely) associated with an $E_r$-sG metric $\gamma_t$ will be called {\bf potentials} of $\gamma_t$. So, the above lemma says that not only can any $E_r$-sG metric $\gamma_0$ on $X_0$ be deformed in a smooth way to $E_r$-sG metrics $\gamma_t$ on the nearby fibres $X_t$, but so can its potentials.

\vspace{2ex}

\noindent {\it Proof of Lemma \ref{Lem:E_r-sG_openness}.} By (ii) of Proposition \ref{Prop:E_r-closed-exact_conditions}, the $E_r$-sG assumption on $\gamma_0$ implies the existence of a $J_0$-type $(n,\,n-2)$-form $\Gamma_0^{n,\,n-2}$ and of a $J_0$-type $(n-1,\,n-1)$-form $\zeta_{r-2,\,0}$ such that $\partial_0\gamma_0^{n-1} = \bar\partial_0\Gamma_0^{n,\,n-2} + \partial_0\zeta_{r-2,\,0}$ and such that \begin{eqnarray}\label{eqn:tower_E_r-exactness_bis}\bar\partial_0\zeta_{r-2,\,0} = \partial_0 v_{r-3,\,0}^{(r-2)}, \hspace{3ex} \mbox{and} \hspace{3ex} \bar\partial_0 v_{r-3,\,0}^{(r-2)} = 0,\end{eqnarray}

\noindent for some $J_0$-type $(n-2,\,n)$-form $v_{r-3,\,0}^{(r-2)}$. (As already pointed out, for bidegree reasons, the general tower (\ref{eqn:tower_E_r-exactness_l}) reduces to (\ref{eqn:tower_E_r-exactness_bis}) in this case.) 

We get $\partial_0(\gamma_0^{n-1} - \zeta_{r-2,\,0} - \overline{\zeta_{r-2,\,0}}) = \bar\partial_0(\Gamma_0^{n,\,n-2} - \overline{v_{r-3,\,0}^{(r-2)}})$, so the $(2n-2)$-form $$\Omega:= -(\Gamma_0^{n,\,n-2} - \overline{v_{r-3,\,0}^{(r-2)}}) + (\gamma_0^{n-1} - \zeta_{r-2,\,0} - \overline{\zeta_{r-2,\,0}}) - \overline{(\Gamma_0^{n,\,n-2} - \overline{v_{r-3,\,0}^{(r-2)}})}$$

\noindent is real and $d$-closed and its $J_0$-pure-type components $\Omega^{n,\,n-2}_0, \Omega^{n-1,\,n-1}_0, \Omega^{n-2,\,n}_0$ are given by the respective paratheses, with their respective signs, on the right of the above identity defining $\Omega$. 

If $\Omega^{n,\,n-2}_t, \Omega^{n-1,\,n-1}_t, \Omega^{n-2,\,n}_t$ stand for the $J_t$-pure-type components of $\Omega$ for any $t\in B$, they all depend in a $C^\infty$ way on $t$. On the other hand, deforming identities (\ref{eqn:tower_E_r-exactness_bis}) in a $C^\infty$ way when the complex structure $J_0$ deforms to $J_t$, we find (non-unique) $C^\infty$ families of $J_t$-type $(n-1,\,n-1)$-forms $(\zeta_{r-2,\,t})_{t\in B}$ and $J_t$-type $(n-2,\,n)$-forms $(v_{r-3,\,t}^{(r-2)})_{t\in B}$, whose elements for $t=0$ are $\zeta_{r-2,\,0}$, respectively $v_{r-3,\,0}^{(r-2)}$, such that $\bar\partial_t\zeta_{r-2,\,t} = \partial_t v_{r-3,\,t}^{(r-2)}$ and $\bar\partial_t v_{r-3,\,t}^{(r-2)} = 0$ for $t\in B$. Then, the $J_t$-type $(n-1,\,n-1)$-form $\Omega^{n-1,\,n-1}_t + \zeta_{r-2,\,t} + \overline{\zeta_{r-2,\,t}}$ depends in a $C^\infty$ way on $t\in B$. When $t=0$, it equals $\gamma_0^{n-1}$, so it is positive definite. By continuity, it remains positive definite for all $t\in B$ sufficiently close to $0\in B$, so it has a unique $(n-1)$-st root and the root is positive definite. In other words, there exists a unique $C^\infty$ positive definite $J_t$-type $(1,\,1)$-form $\gamma_t$ such that $$\gamma_t^{n-1} = \Omega^{n-1,\,n-1}_t + \zeta_{r-2,\,t} + \overline{\zeta_{r-2,\,t}} >0, \hspace{3ex} t\in B,$$

\noindent after possibly shrinking $B$ about $0$. By construction, $\gamma_t$ depends in a $C^\infty$ way on $t$.

If we set $\Gamma_t^{n,\,n-2}:= -\Omega^{n,\,n-2}_t + \overline{v_{r-3,\,t}^{(r-2)}}$ for all $t\in B$ close to $0$, we get $\partial_t\gamma_t^{n-1} = \bar\partial_t\Gamma_t^{n,\,n-2} + \partial_t\zeta_{r-2,\,t}$. Since $\bar\partial_t\zeta_{r-2,\,t} = \partial_t v_{r-3,\,t}^{(r-2)}$ and $\bar\partial_t v_{r-3,\,t}^{(r-2)} = 0$, we conclude that $\gamma_t$ is an $E_r$-sG metric for the complex structure $J_t$ for all $t\in B$ close to $0$.   \hfill $\Box$

\vspace{3ex}

 We are now in a position to prove the first main result of this paper on the deformation limits of a specific class of compact complex manifolds (cf. Theorem \ref{The:limits_E_r_sG_introd} and the comments thereafter.) While it is one of the two building blocks that will yield a proof of Theorem \ref{The:moidef}, we hope that it also holds an independent interest. Note that although $\partial_t\gamma_t^{n-1}$ is of type $(n,\,n-1)$ on $X_t$, its $E_r$-class $\{\partial_t\gamma_t^{n-1}\}_{E_r(X_t)}$ is viewed as an element of the fibre ${\mathcal A}^{2n-1}_{0,\,t} = E_\infty^{n,\,n-1}(X_t)\oplus E_\infty^{n-1,\,n}(X_t)$ of the FAVB ${\mathcal A}^{2n-1}\longrightarrow\C\times B$ at $(0,\,t)$, rather than just as an element of $E_\infty^{n,\,n-1}(X_t)$. Thus, the space $F^nH^{2n-1}(X_t)$ of the standard fibrewise Hodge filtration has no role to play here.

\begin{The}\label{The:limits_E_r_sG} Fix $N\in\N^\star$. Let $\pi:{\mathcal X}\longrightarrow B$ be a holomorphic family of compact complex manifolds over an open ball $B\subset\C^N$ about the origin. Suppose that the fibre $X_t:=\pi^{-1}(t)$ is a {\bf $\partial\bar\partial$-manifold} for all $t\in B\setminus\{0\}$.

  Then, the fibre $X_0:=\pi^{-1}(0)$ is an {\bf $E_r$-sG manifold}, where $r$ is the smallest positive integer such that the Fr\"olicher spectral sequence of $X_0$ degenerates at $E_r$.

  Furthermore, $X_0$ is even an {\bf $E_r$-sGG manifold}.

\end{The}  

\noindent {\it Proof.} Let $\gamma_0$ be an arbitrary Gauduchon metric on $X_0$. It is known that, after possibly shrinking $B$ about $0$, $\gamma_0$ can be extended to a $C^\infty$ family $(\gamma_t)_{t\in\ B}$ of $C^\infty$ $2$-forms on $X$ ($=$ the $C^\infty$ manifold underlying the complex manifolds $X_t$) such that $\gamma_t$ is a Gauduchon metric on $X_t$ for every $t\in B$ (see, e.g., [Pop13, section 3]). Let $n$ be the complex dimension of the fibres $X_t$.

The Gauduchon property of the $\gamma_t$'s implies that $d_{h,\,t}(\partial_t\gamma_t^{n-1})=0$ for all $(h,\,t)\in\C^\star\times B$ and that $\partial_t\gamma_t^{n-1}$ is $E_r(X_t)$-closed for all $t\in B$. Thus, the following object is well defined:

\[
 \sigma(h,\,t): = \left\{
\begin{array}{ll}
  \{\partial_t\gamma_t^{n-1}\}_{d_{h,\,t}}\in H^{2n-1}_{d_{h,\,t}}(X_t,\,\C)={\mathcal A}^{2n-1}_{h,\,t}, & \mbox{if} \hspace{3ex} (h,\,t)\in\C^\star\times B,  \\
  \{\partial_t\gamma_t^{n-1}\}_{E_r(X_t)}\in\bigoplus\limits_{p+q=2n-1}E_r^{p,\,q}(X_t) = {\mathcal A}^{2n-1}_{0,\,t}, & \mbox{if} \hspace{3ex} (h,\,t) = (0,\,t)\in\{0\}\times B,  \\
\end{array}
\right.
\]

\noindent where ${\mathcal A}^{2n-1}\longrightarrow\C\times B$ is the Fr\"olicher approximating vector bundle of the family $(X_t)_{t\in B}$ in degree $2n-1$ defined in Corollary and Definition \ref{Cor-Def:bundle_family}. 

%That result applies in this setting since, for every $t\neq 0$, the $\partial\bar\partial$-assumption on $X_t$ implies $E_1(X_t) = E_\infty(X_t)$, hence also $E_r(X_t) = E_\infty(X_t)$.

Note that the $\partial\bar\partial$-assumption on the fibres $X_t$ with $t\neq 0$ implies that the Fr\"olicher spectral sequence of each of these fibres degenerates at $E_1$, hence also at $E_l$ for every $l\geq 1$ (including $l=r$). Thus, the assumption of Corollary and Definition \ref{Cor-Def:bundle_family} is satisfied and that result yields the holomorphic vector bundle ${\mathcal A}^{2n-1}\longrightarrow\C\times B$ of rank $b_{2n-1}=b_1$ ($=$ the $(2n-1)$-st, respectively the first Betti numbers of $X$, that are equal by Poincar\'e duality).

This last fact, in turn, implies that $\sigma$ is a global $C^\infty$ section of ${\mathcal A}^{2n-1}$ on $\C\times B$. Indeed, $\partial_t$ varies holomorphically with $t\in B$, $\gamma_t^{n-1}$ varies in a $C^\infty$ way with $t\in B$, while the vector space ${\mathcal A}^{2n-1}_{h,\,t}$ varies holomorphically with $(h,\,t)\in\C\times B$. So, (ii) of Corollary and Definition \ref{Cor-Def:bundle_family} applies.

Meanwhile, the $\partial\bar\partial$-assumption on every $X_t$ with $t\in B^\star$ implies that the $d$-closed $\partial_t$-exact $(n,\,n-1)$-form $\partial_t\gamma_t^{n-1}$ is $(\partial_t\bar\partial_t)$-exact, hence also $d_{h,\,t}$-exact for every $h\in\C$. (Indeed, if $\partial_t\gamma_t^{n-1} = \partial_t\bar\partial_t u_t$, then $\partial_t\gamma_t^{n-1} = d_{h,\,t}(-\partial_t u_t)$.) This translates to $\sigma(h,\,t)= \{\partial_t\gamma_t^{n-1}\}_{d_{h,\,t}} = 0\in{\mathcal A}^{2n-1}_{h,\,t}$ for all $(h,\,t)\in\C^\star\times B^\star$. (We even have $\sigma(h,\,t)=0$ for all $(h,\,t)\in\C\times B^\star$.)

Thus, the restriction of $\sigma$ to $\C^\star\times B^\star$ (and even the restriction to $\C\times B^\star$) is identically zero. Then, by continuity, $\sigma$ must be identically zero on $\C\times B$. In particular, $$\sigma(0,\,t) = \{\partial_t\gamma_t^{n-1}\}_{E_r(X_t)} = 0\in{\mathcal A}^{2n-1}_{0,\,t}  \hspace{3ex} \mbox{for all} \hspace{1ex} t\in B,$$

\noindent which means precisely that $\partial_t\gamma_t^{n-1}$ is $E_r(X_t)$-exact for every $t\in B$. In other words, $\gamma_t$ is an $E_r$-sG metric on $X_t$ for every $t\in B$, including $t=0$. In particular, $X_0$ is an $E_r$-sG manifold and even an $E_r$-sGG manifold since the Gauduchon metric $\gamma_0$ was chosen arbitrarily on $X_0$ in the first place. \hfill $\Box$

\begin{Rem}\label{Rem:projections_smooth-var_suffices} The (holomorphic vector bundle structure of the) relative FAVB ${\mathcal A}^{2n-1}\longrightarrow\C\times B$ is not needed in the proof of the above Theorem \ref{The:limits_E_r_sG}. Indeed, the information that the families (\ref {eqn:smooth-families-projections-ht_proof}) of orthogonal projectors are $C^\infty$ suffices.

\end{Rem}

\noindent {\it Proof.} We saw in the proof of Theorem \ref{The:limits_E_r_sG} that, for every $(h,\,t)\in\C^\star\times B^\star$, we have $\sigma(h,\,t)=\{\partial_t\gamma_t^{n-1}\}_{d_{h,\,t}} = 0\in{\mathcal A}^{2n-1}_{h,\,t}$. This amounts to $H^{2n-1}_{h,\,t}(\partial_t\gamma_t^{n-1}) = 0$ for all $(h,\,t)\in\C^\star\times B^\star$.

    In particular, for every fixed $h\in\C^\star$, we have: \begin{eqnarray}\label{eqn:projections_smooth-var_suffices_proof_0}H^{2n-1}_{h,\,t}(\partial_t\gamma_t^{n-1}) = 0,   \hspace{6ex} \mbox{for all}\hspace{1ex} t\in B\setminus\{0\}.\end{eqnarray} Since the family of operators $(H^{2n-1}_{h,\,t})_{t\in B}$ is $C^\infty$ for every fixed $h\in\C^\star$ (because even the family of operators $(H^{2n-1}_{h,\,t})_{(h,\,t)\in\C^\star\times B}$, the first one in (\ref {eqn:smooth-families-projections-ht_proof}) for $k=2n-1$, is $C^\infty$) and since the family of forms $(\partial_t\gamma_t^{n-1} )_{t\in B}$ is $C^\infty$, we infer that the family of forms $$\bigg(H^{2n-1}_{h,\,t}(\partial_t\gamma_t^{n-1})\bigg)_{t\in B}$$ is $C^\infty$ for every fixed $h\in\C^\star$. Hence, by continuity at $t=0$, from (\ref{eqn:projections_smooth-var_suffices_proof_0}) we get: \begin{eqnarray}\label{eqn:projections_smooth-var_suffices_proof_1}H^{2n-1}_{h,\,0}(\partial_0\gamma_0^{n-1}) = 0,   \hspace{6ex} \mbox{for all}\hspace{1ex} h\in\C^\star.\end{eqnarray}

    Now, the family of operators $(H^{2n-1}_{h,\,0})_{h\in\C}$ is $C^\infty$ because the family of operators $(H^{2n-1}_{h,\,t})_{h\in\C}$, the second one in (\ref {eqn:smooth-families-projections-ht_proof}) for $k=2n-1$, is $C^\infty$ for every fixed $t\in B$, including $t=0$. Consequently, the family of forms $$\bigg(H^{2n-1}_{h,\,0}(\partial_0\gamma_0^{n-1})\bigg)_{h\in\C}$$ is $C^\infty$. Together with (\ref{eqn:projections_smooth-var_suffices_proof_1}), the continuity at $h=0$ implies that \begin{eqnarray*}H^{2n-1}_{0,\,0}(\partial_0\gamma_0^{n-1}) = 0.\end{eqnarray*}

    This is equivalent to \begin{eqnarray*}\sigma(0,\,0) = \{\partial_0\gamma_0^{n-1}\}_{E_r(X_0)} = 0\in{\mathcal A}^{2n-1}_{0,\,0},\end{eqnarray*} which proves Theorem \ref{The:limits_E_r_sG}, hence also Remark \ref{Rem:projections_smooth-var_suffices}.  \hfill $\Box$

\begin{Rem}\label{Rem:referee-suggestion}\footnote{The author is grateful for a helpful suggestion leading to this simplification.} As was observed, the following slightly weaker version of Theorem \ref{The:limits_E_r_sG} can be obtained by a drastically simplified variant of the above arguments. 

\end{Rem}

\begin{The}\label{The:limits_E_r_sG_simplified} Fix $N\in\N^\star$. Let $\pi:{\mathcal X}\longrightarrow B$ be a holomorphic family of compact complex manifolds over an open ball $B\subset\C^N$ about the origin. Suppose that the fibre $X_t:=\pi^{-1}(t)$ is a {\bf $\partial\bar\partial$-manifold} for all $t\in B\setminus\{0\}$.

  Then, the fibre $X_0:=\pi^{-1}(0)$ is an {\bf $E_3$-sGG manifold}.

\end{The}  

\noindent {\it Proof.} It runs in three stages.

\vspace{1ex}

$\bullet$ {\it Stage $1$} consists in the following general

\begin{Claim}\label{Claim:stage-1_referee} Let $X$ be a compact complex manifold and $(p,\,q)$ any bidegree. For every $\alpha\in C^\infty_{p,\,q}(X,\,\C)$, the following implication holds: \begin{eqnarray*}\alpha\in\mbox{Im}\,d \implies \alpha \hspace{1ex} \mbox{is}\hspace{1ex} E_r\mbox{-exact} \hspace{3ex}\mbox{for all}\hspace{1ex} r\geq p+1.\end{eqnarray*}

\end{Claim}

\noindent {\it Proof of Claim.} If $\alpha$ is $d$-exact, there exists a $C^\infty$ form $\beta$ of degree $p+q-1$ such that $\alpha = d\beta$. Equating the pure-type parts on either side of this equality and using the bidegree $(p,\,q)$ of $\alpha$, we get, among others, the following $(p+1)$ equalities: \begin{eqnarray*} & & \alpha = \partial\beta^{p-1,\,q} + \bar\partial\beta^{p,\,q-1} \\
  & & \bar\partial\beta^{p-1,\,q} = -\partial\beta^{p-2,\,q+1}, \dots , \bar\partial\beta^{p-l,\,q+l-1} = -\partial\beta^{p-l-1,\,q+l}, \dots , \bar\partial\beta^{0,\,q+p-1} = 0.\end{eqnarray*} This amounts to $\alpha$ being $E_{p+1}$-exact (cf. (ii) of Proposition \ref{Prop:E_r-closed-exact_conditions}), hence also $E_r$-exact for every $r\geq p+1$. \hfill $\Box$

\vspace{1ex}

$\bullet$ {\it Stage $2$} consists in the following general

\begin{Claim}\label{Claim:stage-2_referee} Let $X$ be a compact complex $n$-dimensional manifold and let $\gamma$ be a Gauduchon metric on $X$. The following implication holds: \begin{eqnarray*}\partial\gamma^{n-1}\in\mbox{Im}\,d \implies \gamma \hspace{1ex} \mbox{is an}\hspace{1ex} E_3\mbox{-sG metric}.\end{eqnarray*}

\end{Claim}

\noindent {\it Proof of Claim.} We apply Claim \ref{Claim:stage-1_referee} and its proof to the $(n,\,n-1)$-form $\alpha:=\partial\gamma^{n-1}$. We infer that $\alpha$ is $E_3$-exact (a fact amounting to the metric $\gamma$ being $E_3$-sG) since, for bidegree reasons, there can be at most three non-trivial equalities among those in the proof of Claim \ref{Claim:stage-1_referee}.  \hfill $\Box$

\vspace{1ex}

$\bullet$ {\it Stage $3$: end of proof of Theorem \ref{The:limits_E_r_sG_simplified}}

Let $\gamma_0$ be any Gauduchon metric on $X_0$ and let $(\gamma_t)_{t\in B}$ be any $C^\infty$ family of Gauduchon metrics on the respective fibres $(X_t)_{t\in B}$ whose member for $t=0$ is $\gamma_0$. For every $t\in B$, let $\alpha_t:=\partial_t\gamma_t^{n-1}$, where $n$ is the complex dimension of the fibres $X_t$. Then, $\alpha_t$ depends in a $C^\infty$ way on $t\in B$.

Now, thanks to the Gauduchon property of $\gamma_t$, the $\partial_t$-exact $(n,\,n-1)$-form $\alpha_t$ is $d$-closed for every $t\in B$. Therefore, for every $t\in B\setminus\{0\}$, the $\partial\bar\partial$-hypothesis on $X_t$ implies that $\alpha_t$ is $d$-exact. In other words, the corresponding De Rham classes vanish: \begin{eqnarray*}\{\alpha_t\}_{DR} = 0\in H^{2n-1}_{DR}(X,\,\C), \hspace{5ex} t\in B\setminus\{0\},\end{eqnarray*} where $X$ is the $C^\infty$ manifold underlying all the fibres $X_t$. The dependence on $t\in B$ of $\{\alpha_t\}_{DR}$ being continuous (even $C^\infty$, since the dependence of $\alpha_t$ is), we infer that \begin{eqnarray*}\{\alpha_0\}_{DR} = 0\in H^{2n-1}_{DR}(X,\,\C).\end{eqnarray*} This means that $\alpha_0=\partial_0\gamma_0^{n-1}\in\mbox{Im}\,d$. Hence, by Claim \ref{Claim:stage-2_referee} at {\it Stage $2$}, the Gauduchon metric $\gamma_0$ must be $E_3$-sG.

Since $\gamma_0$ was an arbitrary Gauduchon metric on $X_0$, we are done.  \hfill $\Box$

\subsection{Deformation limits of real $(1,\,1)$-classes}\label{subsection:def-limits_11-classes} By $H^{p,\,q}_{DR}(X,\,\C)$ we mean the space of De Rham cohomology classes of degree $p+q$ that can be represented by a ($d$-closed) pure-type $(p,\,q)$-form. These classes are said to be of {\it type $(p,\,q)$}.

The next statement will play a key role despite its simplicity. It gives a criterion for a {\it real} De Rham $2$-class to be of {\it type $(1,\,1)$} on a possibly {\it non-$\partial\bar\partial$-manifold} that is analogous to the familiar criterion on $\partial\bar\partial$-manifolds requiring the vanishing of the projection onto $H^{0,\,2}(X,\,\C)$ in the canonical Hodge decomposition $H^2_{DR}(X,\,\C) \simeq H^{2,\,0}(X,\,\C)\oplus H^{1,\,1}(X,\,\C)\oplus H^{0,\,2}(X,\,\C)$. On an arbitrary $X$, there is no Hodge decomposition, but its role is played in a certain sense by the non-canonical isomorphism $H^2_{DR}(X,\,\C)\simeq E_\infty^{2,\,0}(X)\oplus E_\infty^{1,\,1}(X)\oplus E_\infty^{0,\,2}(X)$, as the following result shows.

\begin{Lem}\label{Lem:11_02_DR-Einf} Let $\{\alpha\}_{DR}\in H^2_{DR}(X,\,\R)$ be a {\bf real} class. The following equivalence holds: \begin{equation}\label{eqn:11_02_DR-Einf}\{\alpha\}_{DR}\in H^{1,\,1}_{DR}(X,\,\C) \iff \theta_0(\{\alpha\}_{DR}) = 0\in E_\infty^{0,\,2}(X),\end{equation} where $\theta_0:H^2_{DR}(X,\,\C)\longrightarrow E_\infty^{0,\,2}(X)$ is the map defined in Lemma \ref{Lem:theta_0_cohom} by $\theta_0(\{\alpha\}_{DR}) = \{\alpha^{0,\,2}\}_{E_\infty}$.

\end{Lem}

\noindent {\it Proof.} ``$\implies$'' If $\{\alpha\}_{DR}\in H^{1,\,1}_{DR}(X,\,\C)$, there exists a $d$-closed form $u^{1,\,1}\in C^\infty_{1,\,1}(X,\,\C)$ such that $\{\alpha\}_{DR} = \{u^{1,\,1}\}_{DR}$. Then, $\theta_0u^{1,\,1} = 0$, hence $\theta_0(\{\alpha\}_{DR}) = \{\theta_0u^{1,\,1}\}_{DR} = 0\in E_\infty^{0,\,2}(X)$.   

\vspace{1ex}

``$\Longleftarrow$" Since the class $\{\alpha\}_{DR}$ is {\it real}, it can be represented by a {\it real} form $\alpha = \alpha^{2,\,0} + \alpha^{1,\,1} + \alpha^{0,\,2}$. The condition $\alpha = \overline\alpha$ translates to $\alpha^{1,\,1} = \overline{\alpha^{1,\,1}}$ and $\alpha^{2,\,0} = \overline{\alpha^{0,\,2}}$, while the condition $d\alpha=0$ (that $\alpha$ satisfies since it represents a De Rham cohomology class) for the real form $\alpha$ translates to either of the following two equivalent conditions being satisfied: \begin{equation}\label{eqn:d-closedness_proof_11}(\partial\alpha^{2,\,0} = 0 \hspace{1ex}\mbox{and}\hspace{1ex} \partial\alpha^{1,\,1} + \bar\partial\alpha^{2,\,0} = 0) \iff (\bar\partial\alpha^{0,\,2} = 0 \hspace{1ex}\mbox{and}\hspace{1ex} \partial\alpha^{0,\,2} + \bar\partial\alpha^{1,\,1} = 0).\end{equation} 

On the other hand, $\theta_0\alpha = \alpha^{0,\,2}$, so the hypothesis $\theta_0(\{\alpha\}_{DR}) = 0$ amounts to $\{\alpha^{0,\,2}\}_{E_\infty} = 0$. This is equivalent to $\alpha^{0,\,2}$ being $E_r$-exact, where $r$ is the smallest positive integer $l$ such that the Fr\"olicher spectral sequence of $X$ degenerates at $E_l$. However, for bidegree reasons, a $(0,\,q)$-form is $E_r$-exact if and only if it is $\bar\partial$-exact. (See characterisation of $E_r$-exactness in (ii) of Proposition \ref{Prop:E_r-closed-exact_conditions}. In an arbitrary bidegree, $\bar\partial$-exactness, which coincides with $E_1$-exactness, is a stronger property than $E_r$-exactness when $r\geq 2$.) Thus, our assumption $\theta_0(\{\alpha\}_{DR}) = 0$ translates to the existence of a form $u^{0,\,1}\in C^\infty_{0,\,1}(X,\,\C)$ such that $$\alpha^{0,\,2} = \bar\partial u^{0,\,1}.$$

Conjugating the above identity, we get $\alpha^{2,\,0} = \partial u^{1,\,0}$, where we put $u^{1,\,0}:=\overline{u^{0,\,1}}$. This yields: $$\alpha^{2,\,0} + \alpha^{0,\,2} = du - (\bar\partial u^{1,\,0} + \partial u^{0,\,1}), \hspace{5ex} \mbox{where} \hspace{1ex} u:=u^{1,\,0} + u^{0,\,1},$$ hence finally $$\alpha - du = \alpha^{1,\,1} - (\bar\partial u^{1,\,0} + \partial u^{0,\,1}).$$ This shows that $\alpha - du$ is a representative of bidegree $(1,\,1)$ of the De Rham cohomology class $\{\alpha\}_{DR}$, proving that $\{\alpha\}_{DR}\in H^{1,\,1}_{DR}(X,\,\C)$.  \hfill $\Box$

\vspace{3ex}

We can now prove the following 

\begin{The}\label{The:limits_11-classes} Let $\pi:{\mathcal X}\longrightarrow B$ be a holomorphic family of compact complex manifolds over an open ball $B\subset\C^N$ about the origin. Suppose that the fibre $X_t:=\pi^{-1}(t)$ is a {\bf $\partial\bar\partial$-manifold} for all $t\in B\setminus\{0\}$. Let $\{\alpha\}_{DR}\in H^2_{DR}(X,\,\R)$ be a {\bf real} class.

 If $\{\alpha\}_{DR}\in H^{1,\,1}_{DR}(X_t,\,\C)$ for every $t\in B\setminus\{0\}$, then $\{\alpha\}_{DR}\in H^{1,\,1}_{DR}(X_0,\,\C)$.

\end{The}  

\noindent {\it Proof.} Corollary and Definition \ref{Cor-Def:bundle_family} applies in this setting since, for every $t\neq 0$, the $\partial\bar\partial$-assumption on $X_t$ implies $E_1(X_t) = E_\infty(X_t)$, hence also $E_r(X_t) = E_\infty(X_t)$, where $r$ is the smallest positive integer such that $E_r(X_0) = E_\infty(X_0)$.

Let $\theta:{\mathcal H}^2\longrightarrow{\mathcal A}^2$ be the holomorphic vector bundle morphism from the constant bundle of fibre $H^2_{DR}(X,\,\C)$ on $\C\times B$ to the {\it Fr\"olicher approximating vector bundle} (constructed in Corollary and Definition \ref{Cor-Def:bundle_family}) ${\mathcal A}^2\longrightarrow\C\times B$ of the family $(X_t)_{t\in B}$ in degree $2$ defined by the family of linear maps: $$\theta_{h,\,t}:H^2_{DR}(X,\,\C)\longrightarrow{\mathcal A}^2_{h,\,t}, \hspace{5ex} (h,\,t)\in\C\times B.$$ (See (\ref{eqn:theta_ht_rel-case}) and (\ref{eqn:theta_h0_rel-case}).) In particular, $\theta_{h,\,t}$ depends holomorphically on $(h,\,t)\in\C\times B$.    

By Lemma \ref{Lem:11_02_DR-Einf}, the hypothesis $\{\alpha\}_{DR}\in H^{1,\,1}_{DR}(X_t,\,\C)$ for every $t\in B\setminus\{0\}$ translates to $$\theta_{0,\,t}(\{\alpha\}_{DR}) = 0\in{\mathcal A}^2_{0,\,t}, \hspace{5ex} t\in B\setminus\{0\}.$$ Since $\theta_{0,\,0}(\{\alpha\}_{DR}) = \lim\limits_{t\to 0}\theta_{0,\,t}(\{\alpha\}_{DR})$, we get $$\theta_{0,\,0}(\{\alpha\}_{DR}) = 0\in{\mathcal A}^2_{0,\,0} = E_\infty^{2,\,0}(X_0)\oplus E_\infty^{1,\,1}(X_0)\oplus E_\infty^{0,\,2}(X_0).$$ We know from Lemma \ref{Lem:theta_0_cohom} that $\theta_{0,\,0}(\{\alpha\}_{DR})\in E_\infty^{0,\,2}(X_0)$, so $\theta_{0,\,0}(\{\alpha\}_{DR}) = 0\in E_\infty^{0,\,2}(X_0)$.

By Lemma \ref{Lem:11_02_DR-Einf}, this is equivalent to $\{\alpha\}_{DR}\in H^{1,\,1}_{DR}(X_0,\,\C)$ and we are done.  \hfill $\Box$

\vspace{2ex}

We now prove the following complement to Theorem \ref{The:limits_11-classes} that will be needed later on.

\begin{Prop}\label{Prop:smooth-rep_d-closed-parts_11-classes} Let $\pi:{\mathcal X}\longrightarrow B$ be a holomorphic family of compact complex manifolds over an open ball $B\subset\C^N$ about the origin. Suppose that the fibre $X_t:=\pi^{-1}(t)$ is a {\bf $\partial\bar\partial$-manifold} for all $t\in B\setminus\{0\}$. Let $\{\alpha\}_{DR}\in H^2_{DR}(X,\,\R)$ be a {\bf real} class that is of $J_t$-type $(1,\,1)$ (i.e. $\{\alpha\}_{DR}\in H^{1,\,1}_{DR}(X_t,\,\C)$) for every $t\in B$.

  Then, after possibly shrinking $B$ about $0$, for any complex vector line $L\subset\C^N$ and any {\bf real curve} $\Gamma\subset D:=B\cap L$ passing through $0$, the class $\{\alpha\}_{DR}$ can be represented by a $C^\infty$ family $(\alpha_t)_{t\in\Gamma}$ of smooth $d$-closed $2$-forms on $X$ such that, for every $t\in\Gamma$, all the components of pure $J_t$-type of $\alpha_t$ are $d$-closed. 

\end{Prop}

Here and elsewhere, by a {\it real curve} we mean a {\it real-analytic real curve}.

\vspace{1ex}  

We need a few (well-known) preliminaries before proving this statement.

\vspace{1ex}

$\bullet$ The first (standard and obvious) observation is that, for any compact complex manifold $X$ and any $k$, a De Rham $k$-class $\{\alpha\}_{DR}\in H^k_{DR}(X,\,\C)$ can be represented by a form all of whose pure-type components are $d$-closed (we will say in this case that $\{\alpha\}_{DR}$ has property $(\star)$) if and only if $\{\alpha\}_{DR}$ lies in the image of the canonical linear map: \begin{eqnarray*}T:\bigoplus_{p+q=k}H^{p,\,q}_{BC}(X,\,\C)\longrightarrow H^k_{DR}(X,\,\C), \hspace{3ex} \bigg([\alpha^{p,\,q}]_{BC}\bigg)_{p+q=k}\longmapsto\bigg\{\sum\limits_{p+q=k}\alpha^{p,\,q}\bigg\}_{DR}.\end{eqnarray*} This map is neither injective, nor surjective on an arbitrary $X$, but it is bijective if $X$ is a $\partial\bar\partial$-manifold.

Meanwhile, if $k=2$ and $\{\alpha\}_{DR}$ is of type $(1,\,1)$, then $\{\alpha\}_{DR}$ has property $(\star)$. Thus, the class $\{\alpha\}_{DR}\in H^2_{DR}(X,\,\R)$ given in Proposition \ref{Prop:smooth-rep_d-closed-parts_11-classes} lies in the image of the canonical linear map \begin{eqnarray*}T_t:\bigoplus_{p+q=2}H^{p,\,q}_{BC}(X_t,\,\C)\longrightarrow H^2_{DR}(X,\,\C)\end{eqnarray*} for every $t\in B$. Morever, by the $\partial\bar\partial$-assumption, the map $T_t$ is bijective for every $t\in B\setminus\{0\}$.

Let $(\gamma_t)_{t\in B}$ be an arbitrary $C^\infty$ family of Hermitian metrics on the respective fibres $(X_t)_{t\in B}$. If ${\mathcal H}^{p,\,q}_{\Delta_{BC}}(X_t,\,\C)$ is the kernel of the Bott-Chern Laplacian $\Delta_{BC,\,t}:C^\infty_{p,\,q}(X_t,\,\C)\longrightarrow C^\infty_{p,\,q}(X_t,\,\C)$ induced by the metric $\gamma_t$, we have the Hodge isomorphisms  \begin{eqnarray*}{\mathcal H}^{p,\,q}_{\Delta_{BC}}(X_t,\,\C)\simeq H^{p,\,q}_{BC}(X_t,\,\C),  \hspace{5ex} t\in B.\end{eqnarray*} Hence, we get linear maps \begin{eqnarray}\label{eqn:T_t_tilde_def}\widetilde{T}_t:\bigoplus_{p+q=2}{\mathcal H}^{p,\,q}_{\Delta_{BC}}(X_t,\,\C)\longrightarrow H^2_{DR}(X,\,\C),  \hspace{5ex} t\in B,\end{eqnarray} such that the class $\{\alpha\}_{DR}\in H^2_{DR}(X,\,\R)$ given in Proposition \ref{Prop:smooth-rep_d-closed-parts_11-classes} lies in the image of
  $\widetilde{T}_t$ for every $t\in B$ and $\widetilde{T}_t$ is bijective for every $t\in B\setminus\{0\}$.

\vspace{1ex}

$\bullet$ We now recall the following classical result of Grauert's.

\begin{The}([Gra58])\footnote{Years ago, J.-P. Demailly kindly pointed out to the author this result of Grauert's and its Corollary \ref{Cor:real-analytic-variation_families}.}\label{The:Grauert_real-analytic} Let $E\longrightarrow X$ be a {\bf real-analytic} vector bundle over a {\bf real-analytic} manifold $X$. Then, the space of {\bf real-analytic} sections of $E$ is {\bf dense} in the space of $C^\infty$ sections of $E$.

\end{The}

Grauert proved this using the technique of Stein tubular neighbourhoods in the complexified manifold $\widetilde{X}$. As a consequence, we get the following

\begin{Cor}\label{Cor:real-analytic-variation_families} Let $\pi:{\mathcal X}\longrightarrow B$ be a holomorphic family of compact complex manifolds $X_t:=\pi^{-1}(t)$, with $t\in B$, over an open ball $B\subset\C^N$ about the origin for some $N\in\N^\star$.

  \vspace{1ex}

  (1)\, There exists a {\bf real-analytic} family $(\gamma_t)_{t\in B}$ of Hermitian metrics on the respective fibres $X_t$.

  \vspace{1ex}

  (2)\, Taking adjoints and Laplacians w.r.t. to the $\gamma_t$'s, the familiar differential operators: $$\partial_t^\star, \hspace{1ex} \bar\partial_t^\star, \hspace{1ex}, \Delta_t, \hspace{1ex}, \Delta'_t, \hspace{1ex}, \Delta''_t, \hspace{1ex}, \Delta_{BC,\,t}, \hspace{1ex}, \Delta_{A,\,t},$$ vary in a {\bf real-analytic} way with $t\in B$.

  \vspace{1ex}

  (3)\, For any bidegree $(p,\,q)$, any {\bf real-analytic} family $(P_t)_{t\in B}$ of {\bf elliptic} differential operators $P_t:C^\infty_{p,\,q}(X_t,\,\C)\longrightarrow C^\infty_{p,\,q}(X_t,\,\C)$ and for any {\bf Jordan curve} $C\subset\C$ that contains $0\in\C$ in its interior and does not meet the spectrum of $P_0$, there exists a small neighbourhood $B_0\subset B$ of $0$ in $\C^N$ such that $C$ does not meet the spectrum of $P_t$ for any $t\in B_0$ and the vector bundle given by the Kodaira-Spencer theory (see e.g. [Kod86, $\S7.1$]): $$B_0\ni t\mapsto \bigoplus\limits_{\lambda(t)\in\mbox{int}\,(C)}E_{\lambda(t)}(P_t)$$ is {\bf real-analytic}, where $E_{\lambda(t)}(P_t)$ is the $\lambda (t)$-eigenspace of $P_t$.

\end{Cor}

\noindent {\it Proof.} Part (1) follows at once from Grauert's Theorem \ref{The:Grauert_real-analytic} and immediately implies part (2). Part (3) follows from parts (1) and (2) and from the Cauchy integral formula given in the Kodaira-Spencer lemma $7.42$ on p. 340 of [Kod86], by integrating w.r.t. $\lambda$ the Green operators $(P_t - \lambda\,\mbox{Id})^{-1}$ on the Jordan curve.  \hfill $\Box$

\vspace{2ex}

Using Corollary \ref{Cor:real-analytic-variation_families}, we can prove the following

\begin{Prop}\label{Prop:limiting-position_real-anal}\footnote{This statement and its ensuing proof were pointed out years ago to the author by P. Deligne, to whom the former is very grateful.} Let $V\longrightarrow D\subset\C$ be a {\bf real-analytic} $\C$-vector bundle over an open disc about $0$ in the complex plane. Suppose $V$ is equipped with a {\bf real-analytic} fibre metric and that $H:V\longrightarrow V$ is a {\bf real-analytic} Hermitian endomorphism of $V$ (i.e. a {\bf real-analytic} family of self-adjoint operators $H_t:V_t\longrightarrow V_t$) such that $H_t\geq 0$ for all $t\in D$.

Suppose that $H_0=0$ (so, $\ker H_0 = V_0$) and that $\mbox{dim}\,\ker H_t < \mbox{dim}\, V_t$ for all $t\in D\setminus\{0\}$.

Then, for any {\bf real curve} $\Gamma\subset D$ passing through $0$, $\ker H_t$ has a {\bf limiting position} when $\Gamma\setminus\{0\}\ni t\to 0$.

\end{Prop}

By $\ker H_t$ having a {\it limiting position} as $t\in \Gamma\setminus\{0\}$ converges to $0$ we mean that there exists a holomorphic vector bundle over $U\cap D$, where $U$ is an open neighbourhood of $\Gamma$ in $\C$, whose fibre at each $t\in\Gamma\setminus\{0\}$ is $\ker H_t$. The {\it limiting position} of $\ker H_t$ is the fibre at $t=0$ of this vector bundle.

\vspace{2ex}

\noindent {\it Proof of Proposition \ref{Prop:limiting-position_real-anal}.} Restrict $V$ to $\Gamma$ and complexify to get a {\it holomorphic} vector bundle $\widetilde{V}\longrightarrow U\cap D$, where $U$ is an open neighbourhood of $\Gamma$ in $\C$. Similarly, restrict $H$ to $\Gamma$ and complexify to get a {\it holomorphic} endomorphism $\widetilde{H}:\widetilde{V}\longrightarrow\widetilde{V}$. In particular, $$\widetilde{V}_{|\Gamma} = V_{|\Gamma}  \hspace{3ex} \mbox{and} \hspace{3ex} \widetilde{H}_{|\Gamma} = H_{|\Gamma}.$$

Then, $\ker\widetilde{H}$ is a {\it coherent} subsheaf of the locally free sheaf ${\mathcal O}(\widetilde{V})$ (because the kernel of a morphism of coherent sheaves is coherent). Hence, $\ker\widetilde{H}$ is also {\it torsion-free} (because any coherent subsheaf of a torsion-free sheaf is torsion-free). Similarly, $\widetilde{V}/\ker\widetilde{H}$ is a {\it torsion-free coherent sheaf} on $U\cap D$.

Now, every torsion-free coherent sheaf is locally free outside an analytic subset of codimension $\geq 2$. (See e.g. [Kob87, Corollary 5.5.15].) Since the complex dimension of $U\cap D$ is $1$, we get that $\ker\widetilde{H}$ and $\widetilde{V}/\ker\widetilde{H}$ are {\it locally free} on $U\cap D$ and $\widetilde{V}/\ker\widetilde{H}$ is {\it locally a direct factor}.

Moreover, the fibre at $t$ of the holomorphic vector bundle $\ker\widetilde{H}$ is $\ker H_t$ for every $t\in\Gamma$. It follows that the fibre at $t=0$ of $\ker\widetilde{H}$ is a limiting position for $\ker H_t$ when $t\in\Gamma$ approaches $0$. \hfill $\Box$

\vspace{2ex}

$\bullet$ In our case (the setting of Proposition \ref{Prop:smooth-rep_d-closed-parts_11-classes}), we first fix a real-analytic family $(\gamma_t)_{t\in B}$ of Hermitian metrics on the respective fibres $(X_t)_{t\in B}$. Then, by letting $L\subset\C^N$ be a complex vector line and after possibly shrinking $B$ about $0$, we get a sufficiently small open disc $D = B\cap L$ about the origin in $L\simeq\C$.

\vspace{1ex}

The main difficulty in proving Proposition \ref{Prop:smooth-rep_d-closed-parts_11-classes} stems from the fact that the Bott-Chern harmonic spaces need not form a vector bundle near $t=0$ because their dimensions may jump (upwards) at $t=0$. For this reason, we first embed them into a vector bundle $t\mapsto V_t^{p,\,q}$ obtained from the small eigenvalues of the Bott-Chern Laplacian. The harmonic spaces at $t\neq 0$ are then realised as kernels of a real-analytic Hermitian endomorphism of this bundle. Then, Proposition \ref{Prop:limiting-position_real-anal} enables us to extract a limiting holomorphic family of kernels along curves, which replaces the discontinuous family of harmonic spaces. We now carry out this construction in detail.  

\vspace{1ex}

For every bidegree $(p,\,q)\in\{(2,\,0),\,(1,\,1),\,(0,\,2)\}$ and every $t\in D$, let $E^{p,\,q}_{\Delta_{BC,\,t}}(\lambda)$ be the $\lambda$-eigenspace of the Bott-Chern Laplacian $\Delta_{BC,\,t}:C^\infty_{p,\,q}(X_t,\,\C)\longrightarrow C^\infty_{p,\,q}(X_t,\,\C)$ induced by the metric $\gamma_t$ and let $H^{p,\,q}_t$ be the restriction of $\Delta_{BC,\,t}$ to $\oplus_{0\leq\lambda<\varepsilon}E^{p,\,q}_{\Delta_{BC,\,t}}(\lambda)$ for some small $\varepsilon>0$. (We may assume that $\mbox{dim}\,\ker\Delta_{BC,\,0}> \mbox{dim}\,\ker\Delta_{BC,\,t}$ for $t\neq 0$; otherwise, there is nothing to prove.)

%a complex vector subspace ${\mathcal H}_0^{p,\,q}\subset {\mathcal H}_{\Delta_{BC}}^{p,\,q}(X_0,\,\C)$ with the above-mentioned properties by taking

With a view to applying Proposition \ref{Prop:limiting-position_real-anal}, for every $(p,\,q)$, we consider:

\vspace{1ex}

\noindent (i)\, for some $\varepsilon>0$, the real-analytic $\C$-vector bundle $V^{p,\,q}$: $$D\ni t\mapsto \bigoplus\limits_{0\leq\lambda<\varepsilon}E^{p,\,q}_{\Delta_{BC,\,t}}(\lambda):=V^{p,\,q}_t$$ equipped with the real-analytic fibre metric induced by the pointwise inner product associated with $(\gamma_t)_{t\in D}$. (Thus, in order to apply $(3)$ of Corollary \ref{Cor:real-analytic-variation_families}, we take as our Jordan curve $C_\varepsilon$ the circle of radius $\varepsilon$ centred at $0$ in $\C$.);

\vspace{1ex}

\noindent (ii)\, the real-analytic Hermitian endomorphism $H^{p,\,q} = (H^{p,\,q}_t)_{t\in D}:V^{p,\,q}\longrightarrow V^{p,\,q}$.

\vspace{1ex}

Then, we take $\varepsilon>0$ so small that:

\vspace{1ex}

$\cdot$ $0$ is the only eigenvalue of $\Delta_{BC,\,0}$ in the interval $[0,\,\varepsilon)$ (i.e. $V^{p,\,q}_0 = \ker H^{p,\,q}_0 = {\mathcal H}^{p,\,q}_{\Delta_{BC}}(X_0,\,\C)$);

\vspace{1ex}

$\cdot$ the circle $C_\varepsilon\subset\C$ does not meet the spectrum of any $\Delta_{BC,\,t}$ with $t\in D$ (shrink $B$, hence also $D$, about $0$ if necessary);

\vspace{1ex}

$\cdot$ the number of eigenvalues (counted with multiplicities) of $\Delta_{BC,\,t}$ lying inside the circle $C_\varepsilon$ is independent of $t\in D$.

\vspace{1ex}

Our choice of $\varepsilon$ (which can be made the same for all $(p,\,q)$ with $p+q=2$) implies that $H^{p,\,q}_0=0$.

In particular, $V:=\oplus_{p+q=2}V^{p,\,q}$ is a real-analytic $\C$-vector bundle over $D\subset B\subset\C^N$. It is equipped with a real-analytic fibre metric induced by the pointwise inner product associated with $(\gamma_t)_{t\in D}$.

\vspace{1ex}

We then pick an arbitrary real curve $\Gamma\subset D$ through $0$ and use Proposition \ref{Prop:limiting-position_real-anal} to define, for every $(p,\,q)$ with $p+q=2$, the complex vector subspace ${\mathcal H}_0^{p,\,q}\subset {\mathcal H}_{\Delta_{BC}}^{p,\,q}(X_0,\,\C)$ as the limiting position of $\ker H_t^{p,\,q} = {\mathcal H}^{p,\,q}_{\Delta_{BC}}(X_t,\,\C):={\mathcal H}_t^{p,\,q}$ as $t\in\Gamma$ approaches $0$. Actually, Proposition \ref{Prop:limiting-position_real-anal} gives a holomorphic vector bundle ${\mathcal H} = {\mathcal H}^{2,\,0}\oplus{\mathcal H}^{1,\,1}\oplus{\mathcal H}^{0,\,2}$ over $U\cap D$, where $U$ is an open neighbourhood of $\Gamma$ in $L\simeq\C$, whose fibre at each $t\in\Gamma\setminus\{0\}$ is ${\mathcal H}_t = {\mathcal H}^{2,\,0}_{\Delta_{BC}}(X_t,\,\C)\oplus{\mathcal H}^{1,\,1}_{\Delta_{BC}}(X_t,\,\C)\oplus{\mathcal H}^{0,\,2}_{\Delta_{BC}}(X_t,\,\C)$ (the direct sum of the $\Delta_{BC,\,t}$-harmonic spaces) and whose fibre at $t=0$ is ${\mathcal H}_0 = {\mathcal H}^{2,\,0}_0\oplus{\mathcal H}^{1,\,1}_0\oplus{\mathcal H}^{0,\,2}_0$.

In particular, for every bidegree $(p,\,q)$ with $p+q=2$, we have: \begin{eqnarray*}{\mathcal H}_0^{p,\,q}\subset {\mathcal H}_{\Delta_{BC}}^{p,\,q}(X_0,\,\C) = V_0^{p,\,q} \hspace{3ex}\mbox{and}\hspace{3ex} {\mathcal H}_t^{p,\,q} = {\mathcal H}^{p,\,q}_{\Delta_{BC}}(X_t,\,\C)\subset V_t^{p,\,q} \hspace{1ex}\mbox{for}\hspace{1ex} t\in\Gamma\setminus\{0\}.\end{eqnarray*} In the non-trivial case studied here (i.e. when $\mbox{dim}\,\ker\Delta_{BC,\,0}> \mbox{dim}\,\ker\Delta_{BC,\,t}$ for $t\neq 0$), both of the above inclusions are strict.

%\vspace{1ex}
 
%Finally, set $V:=\oplus_{p+q=2}V^{p,\,q}$, a real-analytic $\C$-vector bundle over $D\subset B\subset\C^N$ equipped with the real-analytic fibre metric induced by the pointwise inner product associated with $(\gamma_t)_{t\in D}$. Set $H:=\oplus_{p+q=2}H^{p,\,q}:V\longrightarrow V$, a real-analytic Hermitian endomorphism. Our choice of $\varepsilon$ (which can be made to work for all $(p,\,q)$) implies that $H_0=0$. 

\vspace{2ex}

$\bullet$ {\it End of proof of Proposition \ref{Prop:smooth-rep_d-closed-parts_11-classes}.} Let $\{\alpha\}_{DR}\in H^2_{DR}(X,\,\R)$ be a real class that is of $J_t$-type $(1,\,1)$ for every $t\in B$. In particular, there exists a real representative $\widetilde\alpha_0 = \widetilde\alpha_0^{1,\,1}\in\{\alpha\}_{DR}$ of $J_0$-type $(1,\,1)$. Now, $\widetilde\alpha_0$ differs from the $\Delta_{BC,\,0}$-harmonic representative $\widetilde\alpha_{0,\,BC}$ of its Bott-Chern class $[\widetilde\alpha_0]_{BC}\in H^{1,\,1}_{BC}(X_0,\,\C)$ by a $\partial_0\bar\partial_0$-exact form, hence by a $d$-exact form. Thus, $\widetilde\alpha_{0,\,BC}$ still represents $\{\alpha\}_{DR}$ and is still real of type $(1,\,1)$. Therefore, after possibly replacing $\widetilde\alpha_0$ with $\widetilde\alpha_{0,\,BC}$, we may assume that \begin{eqnarray*}\widetilde\alpha_0 = \widetilde\alpha_0^{1,\,1}\in{\mathcal H}^{1,\,1}_{\Delta_{BC}}(X_0,\,\C)\cap\{\alpha\}_{DR}\subset V^{1,\,1}_0\subset V_0.\end{eqnarray*}

We will now perform two successive projections: first onto $V_t$, then onto ${\cal H}_t$.

For every $t\in D$, let \begin{eqnarray*}\widetilde\alpha_t = \widetilde\alpha_t^{2,\,0} + \widetilde\alpha_t^{1,\,1} + \widetilde\alpha_t^{0,\,2}\in V_t=V^{2,\,0}_t\oplus V^{1,\,1}_t\oplus V^{0,\,2}_t, \hspace{5ex} t\in D,\end{eqnarray*} be the $L^2_{\gamma_t}$-orthogonal projection of $\widetilde\alpha_0\in C^\infty_2(X,\,\C)$ onto $V_t$. Thanks to the Kodaira-Spencer theory, the forms $\widetilde\alpha_t^{p,\,q}$ vary in a $C^\infty$ way with $t\in B$ for every $(p,\,q)$ with $p+q=2$. Moreover, at $t=0$, this projection coincides with the form $\widetilde\alpha_0 = \widetilde\alpha_0^{1,\,1}$ we started off with (in particular, $\widetilde\alpha_0^{2,\,0} = 0$ and $\widetilde\alpha_0^{0,\,2} = 0$) since $V^{p,\,q}_0 = {\mathcal H}^{p,\,q}_{\Delta_{BC}}(X_0,\,\C)$ for all $(p,\,q)$.

Recall that ${\mathcal H}_0$ is the limiting position of \begin{eqnarray*}\ker H_t = \ker \Delta_{BC,\,t} = {\mathcal H}_t = {\mathcal H}^{2,\,0}_{\Delta_{BC}}(X_t,\,\C)\oplus{\mathcal H}^{1,\,1}_{\Delta_{BC}}(X_t,\,\C)\oplus{\mathcal H}^{0,\,2}_{\Delta_{BC}}(X_t,\,\C)\end{eqnarray*} as $t\in\Gamma$ approaches $0$. Since ${\cal H}\to U\cap D$ is a holomorphic vector bundle equipped with the induced Hermitian metric, the orthogonal projectors onto the fibres $({\mathcal H}_t)_{t\in U\cap D}$ vary in a $C^\infty$ way with $t$, by the Kodaira-Spencer theory. Therefore, if \begin{eqnarray*}\alpha_t = \alpha_t^{2,\,0} + \alpha_t^{1,\,1} + \alpha_t^{0,\,2}\in{\mathcal H}_t={\mathcal H}_t^{2,\,0}\oplus {\mathcal H}_t^{1,\,1}\oplus {\mathcal H}_t^{0,\,2}, \hspace{5ex} t\in U\cap D,\end{eqnarray*} are the $L^2_{\gamma_t}$-orthogonal projections of the $C^\infty$ family $(\widetilde\alpha_t)_{t\in U\cap D}$ onto these fibres, we get a $C^\infty$ family $(\alpha_t)_{t\in U\cap D}$ of $2$-forms on the $C^\infty$ manifold $X$ underlying the fibres $X_t$ and a $C^\infty$ family $(\alpha_t)_{t\in\Gamma}$ of representatives of the original class $\{\alpha\}_{DR}$ all of which have the property $(\star)$ (i.e. the property that all their pure-type components are $d$-closed).

To see this last claim, note that \begin{eqnarray*}\{\alpha_t^{2,\,0} + \alpha_t^{1,\,1} + \alpha_t^{0,\,2}\}_{DR} =  \{\alpha\}_{DR}\end{eqnarray*} for all $t\in\Gamma\setminus\{0\}$, thanks to $\{\alpha\}_{DR}\in H^2_{DR}(X,\,\R)$ lying in the image of the {\it bijective} linear map \begin{eqnarray*}\widetilde{T}_t:{\cal H}_t\longrightarrow H^2_{DR}(X,\,\C)\end{eqnarray*} introduced in (\ref{eqn:T_t_tilde_def}) for every $t\in\Gamma\setminus\{0\}$. In other words, the $\Delta_{BC,\,t}$-harmonic form $\alpha_t$ is exactly the unique representative of $\{\alpha\}_{DR}$ inside ${\cal H}_t$ for each $t\in\Gamma\setminus\{0\}$.

Since $\alpha_t^{2,\,0} + \alpha_t^{1,\,1} + \alpha_t^{0,\,2}$ converges to $\alpha_0^{2,\,0} + \alpha_0^{1,\,1} + \alpha_0^{0,\,2}$ as $t\rightarrow 0$ (by the continuity at $0$ of the family $(\alpha_t^{p,\,q})_{t\in U\cap D}$), we get that $\{\alpha\}_{DR} = \{\alpha_t^{2,\,0} + \alpha_t^{1,\,1} + \alpha_t^{0,\,2}\}_{DR}$ converges to $\{\alpha_0^{2,\,0} + \alpha_0^{1,\,1} + \alpha_0^{0,\,2}\}_{DR}$. Thus, $\{\alpha_0^{2,\,0} + \alpha_0^{1,\,1} + \alpha_0^{0,\,2}\}_{DR} = \{\alpha\}_{DR}$. Moreover, for each $(p,\,q)$ with $p+q=2$, we have by construction: $\alpha_t^{p,\,q}\in{\mathcal H}_t^{p,\,q}=\ker \Delta_{BC,\,t}$, hence $d\alpha_t^{p,\,q} = 0$, for all $t\in\Gamma\setminus\{0\}$. The continuity at $0$ of the family $(\alpha_t^{p,\,q})_{t\in\Gamma}$ then yields $d\alpha_0^{p,\,q} = 0$ for each $(p,\,q)$ with $p+q=2$.    \hfill $\Box$

\subsection{Deformation limits of Moishezon manifolds}\label{subsection:def-limits_Moishezon}

We shall now show that the $E_r$-sG property of the limiting fibre $X_0$ proved in Theorem \ref{The:limits_E_r_sG} suffices to prove that any deformation limit of Moishezon manifolds is again Moishezon (cf. Theorem \ref{The:moidef} and the main result in [Pop10]). The result that, together with Theorem \ref{The:limits_E_r_sG}, will prove this fact is the following

\begin{The}\label{The:boundedness} Let $\pi:{\mathcal X}\longrightarrow B$ be a holomorphic family of compact complex $n$-dimensional manifolds over an open ball $B\subset\C^N$ about the origin such that the fibre $X_t:=\pi^{-1}(t)$ is a {\bf $\partial\bar\partial$-manifold} for all $t\in B\setminus\{0\}$. Let $X$ be the $C^\infty$ manifold that underlies the fibres $(X_t)_{t\in B}$ and let $J_t$ be the complex structure of $X_t$.

  Suppose there exists a $C^\infty$ family $(\widetilde\omega_t)_{t\in B}$ of $d$-closed, smooth, real $2$-forms on $X$ such that, for every $t\in B$, the $J_t$-pure-type components of $\widetilde\omega_t$ are $d$-closed. Fix an integer $r\geq 1$ and suppose there exists a $C^\infty$ family $(\gamma_t)_{t\in B}$ of $E_r$-sG metrics on the fibres $(X_t)_{t\in B}$ with potentials depending in a $C^\infty$ way on $t$.

  \vspace{1ex}

  (i)\, If, for every $t\in B^\star$, there exists a K\"ahler metric $\omega_t$ on $X_t$ that is De Rham-cohomologous to $\widetilde\omega_t$, then there exists a constant $C>0$ independent of $t\in B^\star$ such that the $\gamma_t$-masses of the metrics $\omega_t$ are uniformly bounded above by $C$:

  $$0\leq M_{\gamma_t}(\omega_t):=\int\limits_X\omega_t\wedge\gamma_t^{n-1}< C<+\infty,  \hspace{3ex} t\in B^\star.$$

 \noindent In particular, there exists a sequence of points $t_j\in B^\star$ converging to $0\in B$ and a $d$-closed positive $J_0$-$(1,\,1)$-current $T$ on $X_0$ such that $\omega_{t_j}$ converges in the weak topology of currents to $T$ as $j\rightarrow +\infty$.

 \vspace{1ex}

 (ii)\, If, for every $t\in B^\star$, there exists an effective analytic $(n-1)$-cycle $Z_t = \sum_l n_l(t)\,Z_l(t)$ on $X_t$ (i.e. a finite linear combination with integer coefficients $n_l(t)\in\N^\star$ of irreducible analytic subsets $Z_l(t)\subset X_t$ of codimension $1$) that is De Rham-cohomologous to $\widetilde\omega_t$, then there exists a constant $C>0$ independent of $t\in B^\star$ such that the $\gamma_t$-volumes of the cycles $Z_t$ are uniformly bounded above by $C$:

 $$0\leq v_{\gamma_t}(Z_t):=\int\limits_X[Z_t]\wedge\gamma_t^{n-1}< C<+\infty,  \hspace{3ex} t\in B^\star.$$

\end{The}

\noindent {\it Proof.} We will prove (ii). The proof of (i) is very similar and we will indicate the minor differences after the proof of (ii). The method is almost the same as the one in [Pop10].

Since the positive $(1,\,1)$-current $[Z_t]= \sum_l n_l(t)\,[Z_l(t)]$ (a linear combination of the currents $[Z_l(t)]$ of integration on the hypersurfaces $Z_t$) on $X_t$ is De Rham cohomologous to $\widetilde\omega_t$ for every $t\in B^\star$, there exists a {\it real} current $\beta'_t$ of degree $1$ on $X$ such that

\begin{equation}\label{eqn:beta'_t_def}\widetilde\omega_t = [Z_t] + d\beta'_t, \hspace{3ex} t\in B^\star.\end{equation}

\noindent This implies that \begin{equation}\label{eqn:omega-tilde-02_exactness}\bar\partial_t\beta^{'0,\,1}_t = \widetilde\omega_t^{0,\,2}, \hspace{3ex} t\in B^\star.\end{equation}

\noindent In particular, $\widetilde\omega_t^{0,\,2}$ is $\bar\partial_t$-exact for every $t\in B^\star$, so it can be regarded as the right-hand side term of equation (\ref{eqn:omega-tilde-02_exactness}) whose unknown is $\beta^{'0,\,1}_t$.

For every $t\in B^\star$, let $\beta_t^{0,\,1}$ be the minimal $L^2_{\gamma_t}$-norm solution of equation (\ref{eqn:omega-tilde-02_exactness}). Thus, $\beta_t^{0,\,1}$ is the $C^\infty$ $J_t$-type $(0,\,1)$-form given by the Neumann formula

\begin{equation}\label{eqn:Neumann_beta_01_proof}\beta_t^{0,\,1} = \Delta^{''-1}_t\bar\partial_t^\star\widetilde\omega_t^{0,\,2},  \hspace{3ex} t\in B^\star,\end{equation}

\noindent where $\Delta^{''-1}_t$ is the Green operator of the $\bar\partial$-Laplacian $\Delta_t'':=\bar\partial_t\bar\partial_t^\star + \bar\partial_t^\star\bar\partial_t$ induced by the metric $\gamma_t$ on the forms of $X_t$. The difficulty we are faced with is that the family of operators $(\Delta^{''-1}_t)_{t\in B^\star}$, hence also the family of forms $(\beta^{0,\,1}_t)_{t\in B^\star} $, need not extend in a continuous way to $t=0$ if the Hodge number $h^{0,\,1}(t)$ of $X_t$ jumps at $t=0$ (i.e. if $h^{0,\,1}(0)>h^{0,\,1}(t)$ for $t\in B^\star$ close to $0$).

As in [Pop10], the way around this goes through the use of special metrics on the fibres $X_t$. Set

$$\beta_t^{1,\,0}:=\overline{\beta_t^{0,\,1}} \hspace{3ex} \mbox{and} \hspace{3ex} \beta_t:=\beta_t^{1,\,0} + \beta_t^{0,\,1},  \hspace{3ex} t\in B^\star.$$

\noindent Since $\widetilde\omega_t$ is real, this and equation (\ref{eqn:omega-tilde-02_exactness}) satisfied by $\beta_t^{0,\,1}$ imply that $\widetilde\omega_t - [Z_t] - d\beta_t$ is a $J_t$-type $(1,\,1)$-current. Since this current is $d$-exact (it equals $d(\beta_t'-\beta_t)$) and since every fibre $X_t$ with $t\in B^\star$ is supposed to be a $\partial\bar\partial$-manifold, we infer that the current $\widetilde\omega_t - [Z_t] - d\beta_t$ is $\partial_t\bar\partial_t$-exact. (See analogue of (\ref{eqn:dd-bar_def}) for currents and the comment in the Introduction on its equivalence to the smooth-form version of the $\partial\bar\partial$-hypothesis.) Hence, there exists a family of distributions $(R_t)_{t\in B^\star}$ on $(X_t)_{t\in B^\star}$ such that

\begin{equation}\label{eqn:distributions_def}\widetilde\omega_t = [Z_t] + d\beta_t + \partial_t\bar\partial_tR_t \hspace{3ex} \mbox{on} \hspace{1ex} X_t \hspace{3ex} \mbox{for all} \hspace{1ex} t\in B^\star.\end{equation}

Consequently, for the $\gamma_t$-volume of the divisor $Z_t$ we get:

\begin{eqnarray}\label{eqn:volume-estimate_1}v_{\gamma_t}(Z_t):=\int\limits_X[Z_t]\wedge\gamma_t^{n-1} = \int\limits_X\widetilde\omega_t\wedge\gamma_t^{n-1} - \int\limits_X d\beta_t\wedge\gamma_t^{n-1}, \hspace{3ex} t\in B^\star,\end{eqnarray}

\noindent since $\int_X\partial_t\bar\partial_t R_t\wedge\gamma_t^{n-1} = 0$ thanks to the Gauduchon property of $\gamma_t$ and to integration by parts. Now, the families of forms $(\widetilde\omega_t)_{t\in B}$ and $(\gamma_t^{n-1})_{t\in B}$ depend in a $C^\infty$ way on $t$ up to $t=0$, so the quantity $\int_X\widetilde\omega_t\wedge\gamma_t^{n-1}$ is bounded as $t\in B^\star$ converges to $0\in B$. Thus, we are left with proving the boundedness of the quantity $\int_X d\beta_t\wedge\gamma_t^{n-1} = \int_X\partial_t\beta_t^{0,\,1}\wedge\gamma_t^{n-1} + \int_X\bar\partial_t\beta_t^{1,\,0}\wedge\gamma_t^{n-1}$ whose two terms are conjugated to each other. Consequently, it suffices to prove the boundedness of the quantity $$I_t:=\int\limits_X\partial_t\beta_t^{0,\,1}\wedge\gamma_t^{n-1} = \int\limits_X\beta_t^{0,\,1}\wedge\partial_t\gamma_t^{n-1}, \hspace{3ex} t\in B^\star,$$

\noindent as $t$ approaches $0\in B$.

So far, the proof has been identical to the one in [Pop10]. The assumption made on the $C^\infty$ family $(\gamma_t)_{t\in B}$ of $E_r$-sG metrics implies the existence of $C^\infty$ families of $J_t$-type $(n,\,n-2)$-forms $(\Gamma_t^{n,\,n-2})_{t\in B}$ and of $J_t$-type $(n-1,\,n-1)$-forms $(\zeta_{r-2,\,t})_{t\in B}$ such that \begin{eqnarray}\label{eqn:E_r_exactness_gamma_t}\partial_t\gamma_t^{n-1} = \bar\partial_t\Gamma_t^{n,\,n-2} + \partial_t\zeta_{r-2,\,t}, \hspace{3ex} t\in B,\end{eqnarray} \noindent and \begin{eqnarray}\label{eqn:tower_E_r-exactness_l_proof}\bar\partial_t\zeta_{r-2,\,t} & = & \partial_t v^{(r-2)}_{r-3,\,t}  \\
    \nonumber \bar\partial_t v^{(r-2)}_{r-3,\,t} & = &  0.\end{eqnarray}

\noindent (We have already noticed that, for bidegree reasons, tower (\ref{eqn:tower_E_r-exactness_l}) reduces to its first two rows when we start off in bidegree $(n,\,n-1)$.)

 On the other hand, $\bar\partial_t(\partial_t\beta_t^{0,\,1}) = - \partial_t(\bar\partial_t\beta_t^{0,\,1}) = - \partial_t\widetilde{\omega}_t^{0,\,2} = 0$, the last identity being a consequence of the $d$-closedness assumption made on the $J_t$-pure-type components of $\widetilde{\omega}_t$. The $\partial\bar\partial$-assumption on $X_t$ for every $t\in B^\star$ implies that the $J_t$-type $(1,\,1)$-form $\partial_t\beta_t^{0,\,1}$ is $\bar\partial_t$-exact (since it is already $d$-closed and $\partial_t$-exact), so there exist $J_t$-type $(1,\,0)$-forms $(u_t)_{t\in B^\star}$ such that \begin{equation}\label{eqn:u_t_def}\partial_t\beta_t^{0,\,1} = \bar\partial_t u_t,  \hspace{3ex} t\in B^\star.\end{equation}

 \noindent This, in turn, implies that the $J_t$-type $(2,\,0)$-form $\partial_t u_t$ is $\bar\partial_t$-closed, hence $d$-closed. The $\partial\bar\partial$-assumption on $X_t$ for every $t\in B^\star$ implies that $\partial_t u_t$ is $\bar\partial_t$-exact, hence zero, for bidegree reasons. Thus \begin{equation}\label{eqn:del-u_t_zero}\partial_t u_t = 0,  \hspace{3ex} t\in B^\star.\end{equation}

 Putting (\ref{eqn:E_r_exactness_gamma_t}), (\ref{eqn:tower_E_r-exactness_l_proof}), (\ref{eqn:u_t_def}) and (\ref{eqn:del-u_t_zero}) together and integrating by parts several times, we get: \begin{eqnarray}\label{eqn:I_t_computation_1}\nonumber I_t & = & \int\limits_X\bar\partial_t\beta_t^{0,\,1}\wedge\Gamma_t^{n,\,n-2} + \int\limits_X\partial_t\beta_t^{0,\,1}\wedge\zeta_{r-2,\,t} = \int\limits_X\widetilde{\omega}_t^{0,\,2}\wedge\Gamma_t^{n,\,n-2} + \int\limits_X \bar\partial_t u_t\wedge\zeta_{r-2,\,t} \\
   \nonumber & = & \int\limits_X\widetilde{\omega}_t^{0,\,2}\wedge\Gamma_t^{n,\,n-2} + \int\limits_X u_t\wedge\bar\partial_t\zeta_{r-2,\,t} = \int\limits_X\widetilde{\omega}_t^{0,\,2}\wedge\Gamma_t^{n,\,n-2} + \int\limits_X u_t\wedge\partial_t v^{(r-2)}_{r-3,\,t} \\
   \nonumber & = & \int\limits_X\widetilde{\omega}_t^{0,\,2}\wedge\Gamma_t^{n,\,n-2} + \int\limits_X \partial_t u_t\wedge v^{(r-2)}_{r-3,\,t} = \int\limits_X\widetilde{\omega}_t^{0,\,2}\wedge\Gamma_t^{n,\,n-2}, \hspace{3ex} t\in B^\star.\end{eqnarray}

 \noindent Since the families of forms $(\Gamma_t^{n,\,n-2})_{t\in B}$ and $(\widetilde{\omega}_t^{0,\,2})_{t\in B}$ vary in a $C^\infty$ way with $t$ up to $t=0\in B$, we infer that the quantities $(I_t)_{t\in B^\star}$ are bounded as $t\in B^\star$ converges to $0\in B$. This completes the proof of (ii). 

 The proof of (i) is identical to that of (ii), except for the fact that $[Z_t]$ has to be replaced by $\omega_t$ in (\ref{eqn:beta'_t_def}), (\ref{eqn:distributions_def}) and (\ref{eqn:volume-estimate_1}), while $\beta_t'$ and $R_t$ are smooth.  \hfill $\Box$

 \vspace{3ex}

 Before continuing, we make a very simple observation.

 \begin{Lem}\label{Lem:volumes-equal_ddbar-G} Let $X$ be a compact complex $n$-dimensional {\bf $\partial\bar\partial$-manifold} and let $\omega$ be a Gauduchon metric on $X$. (More generally, suppose that there exists a {\bf strongly Gauduchon} metric $\omega$ on $X$.)

   Then, any two divisors $Z,Z'$ on $X$ whose currents of integration lie in the same De Rham cohomology class $\{[Z]\}_{DR} = \{[Z']\}_{DR}$ have equal $\omega$-volumes: $$v_\omega(Z) = \int\limits_X[Z]\wedge\omega^{n-1} = \int\limits_X[Z']\wedge\omega^{n-1} = v_\omega(Z').$$

\end{Lem}   

 \noindent {\it Proof.} As noticed in [Pop13], any Gauduchon metric on a $\partial\bar\partial$-manifold is strongly Gauduchon ($E_1$-sG in the language of this paper, see Definition \ref{Def:E_r-sG}). Meanwhile, by [Pop13, Proposition 4.2], the metric $\omega$ being strongly Gauduchon is equivalent to the form $\omega^{n-1}$ being the component $\Omega^{n-1,\,n-1}$ of type $(n-1,\,n-1)$ of some $C^\infty$ real $(2n-2)$-form $\Omega$ on $X$ such that $d\Omega=0$.

 Since the currents $[Z]$ and $[Z']$ are of bidegree $(1,\,1)$, we have $[Z]\wedge\Omega^{n-1,\,n-1} = [Z]\wedge\Omega$ and the analogous equality for $[Z']$. Meanwhile, the hypothesis $\{[Z]\}_{DR} = \{[Z']\}_{DR}$ translates to the existence of a real $1$-current $S$ on $X$ such that $[Z] = [Z'] +dS$.

 We get: \begin{eqnarray*}v_\omega(Z) = \int\limits_X[Z]\wedge\Omega = \int\limits_X[Z']\wedge\Omega + \int\limits_XS\wedge d\Omega = \int\limits_X[Z']\wedge\Omega = v_\omega(Z')\end{eqnarray*} after using Stokes's theorem and the equality $d\Omega=0$.   \hfill $\Box$

\vspace{3ex}

 We are now in a position to show that our results obtained above combine with classical results on the relative Barlet space of divisors to yield the key Theorem \ref{The:dd-bar_def} from which the main result of this paper, Theorem \ref{The:moidef}, follows by standard arguments. We first recall a few well-known facts about (relative) cycles that will be used.

 The relative Barlet space ${\mathcal C}^{n-1}({\mathcal X}/B)$ (cf. [Bar75]) of effective analytic divisors ($=$ effective analytic $(n-1)$-cycles) $Z_t$ contained in the fibres $X_t$ (whose complex dimension is denoted by $n$) of a holomorphic family $\pi:{\mathcal X}\rightarrow B$ of compact complex manifolds is a {\it closed analytic subset} of the (absolute) Barlet space ${\mathcal C}^{n-1}({\mathcal X})$ of compact $(n-1)$-cycles on ${\mathcal X}$ (see e.g. [BM14, th\'eor\`eme 8.2.2., p. 481]). Moreover, the canonical projection \begin{eqnarray*}\mu_{n-1}:{\mathcal C}^{n-1}({\mathcal X}/B)\rightarrow B, \hspace{3ex} \mu_{n-1}(Z_t)=t,\end{eqnarray*} mapping every relative divisor $Z_t$ of every fibre to the unique point $t\in B$ such that the fibre $X_t$ above $t$ contains the support $|Z_t|$ of $Z_t$, is {\it holomorphic} (see e.g. [BM14, remarque (i), p. 484]).

 Recall that ${\mathcal C}({\mathcal X}):=\cup_p{\mathcal C}^p({\mathcal X})$ is the Chow scheme of ${\mathcal X}$ (which, by definition, parametrises the compactly supported analytic $p$-cycles of ${\mathcal X}$ for all $p$, namely the finite formal linear combinations $\sum_ln_l\,Z_l$ of irreducible {\it compact} $p$-dimensional subvarieties $Z_l$ of ${\mathcal X}$ with positive integers $n_l$ as coefficients) that Barlet endowed with a natural structure of a Banach analytic set whose irreducible components are finite-dimensional analytic sets (cf. [Bar75]). Moreover, any irreducible component $S$ of ${\mathcal C}({\mathcal X})$ arises as an analytic family of compact cycles $(Z_s)_{s\in S}$ parametrised by $S$, while giving an analytic family $(Z_s)_{s\in S}$ of compact cycles of dimension $p$ on ${\mathcal X}$ is equivalent to giving an analytic subset

$${\mathcal Z}=\{(s,\, z)\in S\times {\mathcal X}\, /\, z\in |Z_s|\}\subset S\times{\mathcal X},$$

\noindent where $|Z_s|$ denotes the support of the cycle $Z_s$, such that the restriction to ${\mathcal Z}$ of the natural projection on $S$ is proper, surjective and has fibres of pure dimension $p$ (cf. [Bar75, th\'eor\`eme 1, p. 38]). 

Recall Lieberman's strengthened form ([Lie78, Theorem 1.1]) of Bishop's Theorem [Bis64]: a subset $S\subset{\mathcal C}({\mathcal X})$ is {\it compact} if and only if the supports $|Z_s|$, $s\in S$, all lie in a same {\it compact} subset of ${\mathcal X}$ and the $\widetilde{\omega}$-volume of the $p$-cycle $Z_s$, namely the quantity $$v_{\widetilde{\omega}}(Z_s):=\int\limits_{{\mathcal X}}[Z_s]\wedge\widetilde{\omega}^p=\sum\limits_l n_{s,\,l}\,\int\limits_{Y_{s,\,l}}\widetilde{\omega}^p,$$ where $[Z_s]=\sum_l n_{s,\,l}[Y_{s,\,l}]$ is the current of integration over the cycle $Z_s:=\sum_l n_{s,\,l}Y_{s,\,l}$, is {\it uniformly bounded} when $s$ ranges over $S$ for some (hence any) Hermitian metric $\widetilde{\omega}$ on ${\mathcal X}$.

The proof of this Bishop-Lieberman result uses the continuity of the volume map \begin{eqnarray*}S\ni s\mapsto v_{\widetilde{\omega}}(Z_s)\end{eqnarray*} for every irreducible component $S=(Z_s)_{s\in S}$ of ${\mathcal C}({\mathcal X})$.

While the irreducible components of the Barlet space ${\mathcal C}^p(X)$ of $p$-cycles need not be compact for an arbitrary $p$ on a general compact complex manifold $X$ (cf. [Lie78]), compactness of the irreducible components of the Barlet space ${\mathcal C}^{n-1}(X)$ of {\it divisors} of $X$ always holds if $X$ is compact (see e.g. [CP94, Remark 2.18.]). Thus, the absolute case of the following Theorem \ref{The:dd-bar_def} (i.e. when $B$ is reduced to a point) is well known and no special assumption is necessary.

%by piecing together the above results, the following statement that easily implies the main result of this paper, Theorem \ref{The:moidef}, via classical results on the relative Barlet space ${\mathcal C}^{n-1}({\mathcal X}/B)$ (cf. [Bar75]) of effective analytic divisors $Z_t$ contained in the fibres $X_t$ of a holomorphic family $(X_t)_{t\in B}$ of compact complex manifolds. 

\begin{The}\label{The:dd-bar_def} Let $\pi:{\mathcal X}\rightarrow B$ be a complex analytic family of compact complex $n$-dimensional manifolds over an open ball $B\subset\C^N$ about the origin such that the fibre $X_t:=\pi^{-1}(t)$ is a {\bf $\partial\bar\partial$-manifold} for every $t\in B\setminus\{0\}$. Then, the canonical holomorphic projection \begin{eqnarray*}\mu_{n-1}:{\mathcal C}^{n-1}({\mathcal X}/B)\rightarrow B, \hspace{3ex} \mu_{n-1}(Z_t)=t,\end{eqnarray*} mapping every divisor $Z_t\subset X_t$ contained in some fibre $X_t$ to the base point $t\in B$, has the property that its restrictions to the irreducible components of ${\mathcal C}^{n-1}({\mathcal X}/B)$ are {\bf proper}.

\end{The}

\noindent {\it Proof.} By Theorem \ref{The:limits_E_r_sG}, $X_0$ is an $E_r$-sG manifold, where $r\in\N^\star$ is the smallest positive integer such that $E_r(X_0) = E_{\infty}(X_0)$. Therefore, thanks to Lemma \ref{Lem:E_r-sG_openness}, after possibly shrinking $B$ about $0$, there exists a $C^\infty$ family $(\gamma_t)_{t\in B}$ of $E_r$-sG metrics on the fibres $(X_t)_{t\in B}$ whose potentials depend in a $C^\infty$ way on $t\in B$.

To show properness over $B$ of an arbitrary irreducible component $S\subset{\mathcal C}^{n-1}({\mathcal X}/B)$, one has to show that for every compact subset $K\subset B$, $\mu_{n-1}^{-1}(K)\cap S$ is a compact subset of ${\mathcal C}^{n-1}({\mathcal X}/B)$. If $(Z_s)_{s\in S}$ is the analytic family of divisors parametrised by $S$ (such that $Z_s\subset X_{\mu_{n-1}(s)}$, $s\in S$), this amounts to proving, by the Bishop-Lieberman results of [Bis64] and [Lie78] recalled above, that the volumes \begin{eqnarray*}v_{\gamma_s}(Z_s)=\int\limits_X[Z_s]\wedge\gamma_s^{n-1}\end{eqnarray*} are uniformly bounded when $s$ ranges over $\mu_{n-1}^{-1}(K)\cap S$. We have set for convenience $\gamma_s=\gamma_{\mu_{n-1}(s)}$.

A standard observation is that the De Rham cohomology class $\{[Z_s]\}_{DR}\in H^2(X,\,\R)$ of the current $[Z_s]=\sum_l n_{s,\,l}[Y_{s,\,l}]$ of integration over any divisor $Z_s=\sum_l n_{s,\,l}Y_{s,\,l}$ (where $n_{s,\,l}\in\N^\star$ and the $Y_{s,\,l}$'s are compact irreducible analytic hypersurfaces of $X_{\mu_{n-1}(s)}$) is integral. Therefore, the continuous, integral-class-valued map $$S\ni s\mapsto\{[Z_s]\}_{DR}\in H^2(X,\,\Z)$$ is constant.

As mentioned above, it is known that the compactness of the fibres $X_t$ implies that the absolute Barlet space ${\mathcal C}^{n-1}(X_t)$ of divisors of every fibre $X_t$ has {\it compact} irreducible components. Thus, since the volumes $v_{\gamma_s}(Z_s)$ depend continuously on $s\in S$, the volume $v_{\gamma_s}(Z_s)$ stays {\it uniformly bounded} when $Z_s$ varies across any irreducible component of any fixed fibre $X_t$ with $t\in B$. Moreover, for every $t\neq 0$, the fibre $X_t$ is a $\partial\bar\partial$-manifold, by hypothesis. Hence, thanks to the $E_r$-sG metrics $\gamma_s$ being necessarily Gauduchon and to the De Rham cohomology class of all the currents of integration $[Z_s]$ being constant, Lemma \ref{Lem:volumes-equal_ddbar-G} implies that the volume $v_{\gamma_s}(Z_s)$ stays even {\it constant} when $Z_s$ varies across any irreducible component of any fixed fibre $X_t$ with $t\in B^\star$.

Thus, to prove the uniform boundedness of the family $(v_{\gamma_s}(Z_s))_{s\in\mu_{n-1}^{-1}(K)\cap S}$ of volumes we may assume without loss of generality that for each $t\in B^\star$, the fibre $X_t$ contains the support of at most one of the relative divisors $(Z_s)_{s\in\mu_{n-1}^{-1}(K)\cap S}$. Moreover, if either the set $\mu_{n-1}(S)$ contains $0$ or there exists a neighbourhood of $0$ in $B$ that does not meet $\mu_{n-1}(S)$, the uniform boundedness of the volumes follows at once from the above arguments. 

It then suffices to show the uniform boundedness of the volumes $$(v_{\gamma_s}(Z_s))_{s\in\mu_{n-1}^{-1}(K)\cap S}$$ in the case where $0\notin\mu_{n-1}(S)$ but $0$ is a limit point of $\mu_{n-1}(S)$. In other words, after recalling that $\mu_{n-1}(S)$ is connected (since $\mu_{n-1}$ is continuous) and after replacing the open ball $B$ about the origin of $\C^N$ with a smaller ball about the origin, we may assume that the family $(Z_s)_{s\in\mu_{n-1}^{-1}(K)\cap S}$ is a continuous family $(Z_t)_{t\in B^{\star}}$ of effective analytic divisors such that $Z_t\subset X_t$ for every $t\in B^\star$.

It has already been argued that the continuous, integral-class-valued map $$B^\star\ni t\mapsto\{[Z_t]\}_{DR}\in H^2(X,\,\Z)$$ must be constant, equal to an integral De Rham $2$-class that we denote by $\{\alpha\}_{DR}$. Moreover, the current of integration $[Z_t]$ is of bidegree $(1,\,1)$ for the complex structure $J_t$  of $X_t$, so $\{\alpha\}_{DR}\in H^{1,\,1}_{DR}(X_t,\,\C)$ for every $t\in B^\star$. By Theorem \ref{The:limits_11-classes}, $\{\alpha\}_{DR}\in H^{1,\,1}_{DR}(X_0,\,\C)$. Thus, $\{\alpha\}_{DR}$ satisfies the hypotheses of Proposition \ref{Prop:smooth-rep_d-closed-parts_11-classes}.

Therefore, after possibly shrinking $B$ about $0$ and then replacing it with any real curve $\Gamma\subset D$ through $0\in\C^N$, where $D:=B\cap L\subset B$ and $L$ is  a complex vector line in $\C^N$, Proposition \ref{Prop:smooth-rep_d-closed-parts_11-classes} ensures the existence of a $C^\infty$ family $(\widetilde\omega_t)_{t\in B}$ of $d$-closed, smooth, real $2$-forms on $X$ lying in the De Rham class $\{\alpha\}$ such that, for every $t\in B$, the $J_t$-pure-type components of $\widetilde\omega_t$ are $d$-closed. (We may replace $B$ with $\Gamma$ since ours is a problem for families of manifolds over $1$-dimensional bases $B$.) In particular, for every $t\in B^\star$, the current $[Z_t]$ is De Rham-cohomologous to $\widetilde\omega_t$.

Thus, all the hypotheses of Theorem \ref{The:boundedness} are satisfied. From (ii) of that theorem we get that the $\gamma_t$-volumes $(v_{\gamma_t}(Z_t))_{t\in B^\star}$ of the divisors $Z_t$ are uniformly bounded. This implies, thanks to [Bis64] and [Lie78, Theorem 1.1], that a limiting effective divisor $Z_0\subset X_0$ for the family of relative effective divisors $(Z_t)_{t\in B^\star}$ exists and the family is compact.   \hfill $\Box$

%Since this family has been chosen arbitrarily, it follows that $X_0$ has at least as many divisors as the nearby fibres $X_t$ with $t\neq 0$ and $t$ close to $0$. Meanwhile, we know (see, e.g., [CP94, Remark 2.22]) that the algebraic dimension of any compact complex manifold $X$ is the maximal number of effective prime divisors meeting transversally at a generic point of $X$. It follows that the algebraic dimension of $X_0$ is $\geq$ the algebraic dimension of the generic fibre $X_t$ with $t\in B^\star$ close to $0$.   \hfill $\Box$

\vspace{3ex}

Recall that the {\it algebraic dimension} $a(X)$ of a compact complex $n$-dimensional manifold $X$ is the maximal number of algebraically independent meromorphic functions on $X$. Equivalently, $a(X)$ is the transcendence degree over $\C$ of the field of meromorphic functions on $X$. It is standard that $a(X)\leq n$ and that $a(X) = n$ if and only if $X$ is Moishezon ([Moi67]). Since every meromorphic function gives rise to its divisor of zeros and poles, Moishezon manifolds can be regarded as the compact complex manifolds that carry ``many'' divisors.

On the other hand, we recall the following standard facts. The properness proved in Theorem \ref{The:dd-bar_def} under the assumptions therein guarantees that the images of the irreducible components of ${\mathcal C}^{n-1}({\mathcal X}/B)$ under $\mu_{n-1}$ are analytic subsets of $B$ thanks to Remmert's Proper Mapping Theorem. Let $\Sigma_{\nu}\subsetneq B$, for $\nu\in\Z$, be those such images (at most countably many) that are {\it strictly} contained in $B$. Each $\Sigma_{\nu}$ is thus a proper analytic subset of $B$.

Let $S$ be an arbitrary irreducible component of ${\mathcal C}^{n-1}({\mathcal X}/B)$. As recalled above, it gives rise to an analytic family (in the sense of [Bar75, Th\'eor\`eme 1, p. 38]) of relative effective divisors $(Z_s)_{s\in S}$ such that $Z_s\subset X_{\mu_{n-1}(s)}$ for every $s\in S$. We can either have: \begin{equation}\label{eqn:fam-cycles1}\mu_{n-1}(S)=B \hspace{3ex} \mbox{or}\end{equation} \begin{equation}\label{eqn:fam-cycles2}\mu_{n-1}(S)=\Sigma_\nu\subsetneq B \hspace{2ex} \mbox{for some}\,\,\, \nu\in\Z.\end{equation}

 Let $\Sigma=\bigcup_{\nu}\Sigma_{\nu}\subsetneq B$. Thus, every divisor $Z_{s_0}$ contained in a fibre $X_{t_0}$ lying above some point $t_0=\mu_{n-1}(s_0)\in B\setminus\Sigma$ (call such a fibre {\it generic}) stands in an analytic family of divisors $(Z_s)_{s\in S}$ covering the whole base $B$ as in (\ref{eqn:fam-cycles1}) (call these divisors {\it generic}), while the {\it exceptional} fibres $X_t$ (i.e. those above points $t\in\Sigma$) may have extra divisors (those standing in {\it isolated} families satisfying (\ref{eqn:fam-cycles2})) besides the {\it generic} divisors that ``sweep'' $B$ in families with the property (\ref{eqn:fam-cycles1}). Thus, intuitively, the properness over the base of the irreducible components of ${\mathcal C}^{n-1}({\mathcal X}/B)$ (proved in our case in Theorem \ref{The:dd-bar_def}) ensures that every fibre (in particular $X_0$) has ``at least as many'' divisors (at least the {\it generic} ones) as the {\it generic} fibres of the family.

\vspace{2ex}

Finally, recall the by now classical results of Fujiki and Campana in [Fuj78] and [Cam81, Theorem 1 and Corollaries 1, 2, 3]. In particular, in [Cam81, corollaire 2] the ``weakly K\"ahler'' assumption on the morphism $\pi:X\to S$ (that became $\pi:{\mathcal X}\rightarrow B$ in our notation) is made only to ensure the properness over the base of the relative space of divisors that was proved for our situation in Theorem \ref{The:dd-bar_def}. Indeed, Fujiki had earlier proved in [Fuj78, Theorem 4.5 and Proposition 4.8] the properness over the base of the relative space of cycles of any dimension under the assumption that the morphism $\pi$ is ``weakly K\"ahler''. Let us add that it was Lieberman who initiated this series of works on Chow compactness and that in [Cam81] ``algebraic'' means ``Moishezon''.

\vspace{3ex}

Thus, using these standard facts and our Theorem \ref{The:dd-bar_def}, we can finish the proof of our main result. 

%\begin{Cor}\label{Cor:semicont_alg-dim} Let $\pi:{\mathcal X}\rightarrow B$ be a complex analytic family of compact complex manifolds over an open ball $B\subset\C^N$ about the origin such that the fibre $X_t:=\pi^{-1}(t)$ is a {\bf $\partial\bar\partial$-manifold} for every $t\in B\setminus\{0\}$. Then $a(X_0)\geq a(X_t)$ for all $t\in B \setminus\{0\}$ sufficiently close to $0$, where $a(X_t)$ is the algebraic dimension of $X_t$.

%\end{Cor}

\vspace{3ex}

\noindent {\it Proof of Theorem \ref{The:moidef}.} Let $n=\mbox{dim}_\C X_t$ for all $t\in B$. The Moishezon property is well known to imply the $\partial\bar\partial$-property, so the fibre $X_t$ is a $\partial\bar\partial$-manifold for every $t\in B\setminus\{0\}$. Therefore, Theorem \ref{The:dd-bar_def} applies and ensures the properness (in the precise sense spelt out therein) over the base of the relative Barlet space of divisors associated with the family $(X_t)_{t\in B}$. 

By [Cam81, corollaire 2, p. 160], this properness guarantees that, whenever the fibre $X_t$ is Moishezon for every $t$ in a subset $B'\subset B$ which is not ``analytically meager'' in $B$ in the sense of [Cam 81, d\'efinition 1, p.158] (so, we can choose $B'=B\setminus\{0\}$ in our case), the fibre $X_t$ is Moishezon for all $t\in B$. In particular, $X_0$ must be Moishezon.  \hfill $\Box$ 

%tells us that $a(X_0)\geq a(X_t)$ for all $t\in B \setminus\{0\}$. Meanwhile, $a(X_t)=n$ for every $t\in B\setminus\{0\}$ by the Moishezon assumption on every $X_t$ with $t\in B\setminus\{0\}$. Since $a(X_0)\leq\mbox{dim}_\C X_0=n$, we must have $a(X_0)=n$. Hence, $X_0$ must be Moishezon.  \hfill $\Box$

\vspace{3ex}

Note that Theorem \ref{The:moidef} is also, implicitly, an upper semicontinuity result for the algebraic dimensions of the fibres of a holomorphic family of compact complex manifolds whose generic fibre is assumed to be Moishezon. (Actually, our method shows, when combined with classical results on the relative Barlet space of cycles, that semicontinuity holds under the weaker $\partial\bar\partial$-assumption on each fibre $X_t$ with $t\neq 0$, as a consequence of our properness Theorem \ref{The:dd-bar_def}, but we do not pursue this discussion here.) Without the $\partial\bar\partial$-assumption on $X_t$ with $t\neq 0$, the statement is known to fail even when the fibres are complex surfaces. An example of a family of compact complex surfaces of class VII (hence non-K\"ahler and even non-$\partial\bar\partial$), whose algebraic dimension drops from $1$ on the generic fibre $X_t$ to $0$ on the limiting fibre $X_0$, was constructed by Fujiki and Pontecorvo in [FP10].

\vspace{3ex}

\noindent {\bf References.} \\

\noindent [Aso06]\, A. Asok --- {Equivariant Vector Bundles on Certain Affine G-varieties} --- Pure Appl. Math. Q. {\bf 2} (2006), 1085-1102. arXiv:math/0604344v1

\vspace{1ex}

\noindent [Bar75]\, D. Barlet --- {\it Espace analytique r\'eduit des cycles analytiques complexes compacts d'un espace analytique complexe de dimension finie} --- Fonctions de plusieurs variables complexes, II (S\'em. Fran\c{c}ois Norguet, 1974-1975), LNM, Vol. {\bf 482}, Springer, Berlin (1975) 1-158.

\vspace{1ex}

\noindent [BM14]\, D. Barlet, J. Magn\'usson --- {\it Cycles analytiques complexes I : th\'eor\`emes de pr\'eparation des cycles} --- Cours sp\'ecialis\'es, collection SMF, {\bf 22} (2014).

\vspace{1ex}

\noindent [Bis64]\, E. Bishop --- {\it Conditions for the Analyticity of Certain Sets} --- Mich. Math. J. {\bf 11} (1964) 289-304.

\vspace{1ex}

\noindent [Cam81]\, F. Campana --- {\it R\'eduction alg\'ebrique d'un morphisme faiblement k\"ahl\'erien propre et applications} --- Math. Ann. {\bf 256} (1981), 157-189.

\vspace{1ex}

\noindent [CP94]\, F. Campana, T. Peternell --- {\it Cycle spaces} --- Several Complex Variables, VII,  319-349, Encyclopaedia Math. Sci., {\bf 74}, Springer, Berlin (1994).

\vspace{1ex}

\noindent [CFGU97]\, L.A. Cordero, M. Fern\'andez, A.Gray, L. Ugarte --- {\it A General Description of the Terms in the Fr\"olicher Spectral Sequence} --- Differential Geom. Appl. {\bf 7} (1997), no. 1, 75--84.

\vspace{1ex}

\noindent [Ehr47]\, C. Ehresmann ---{\it Sur les espaces fibr\'es diff\'erentiables} -- C. R. Acad. Sci. Paris {\bf 224} (1947), 1611-1612.

\vspace{1ex}

\noindent [Fro55]\, A. Fr\"olicher --- {\it Relations between the Cohomology Groups of Dolbeault and Topological Invariants} --- Proc. Nat. Acad. Sci. U.S.A. 41 (1955), 641–644. 

\vspace{1ex}

\noindent [Fuj78]\, A. Fujiki --- {\it Closedness of the Douady Space of a Compact K\"ahler Space} --- Publ. Res. Inst. Math. Sci. {\bf 14} (1978), 1-52.

\vspace{1ex}

\noindent [FP10]\, A. Fujiki, M. Pontecorvo --- {\it Non-Upper-Semicontinuity of Algebraic Dimension for Families of Compact Complex Manifolds} --- Math. Ann. {\bf 348}, no. 3 (2010), 593–599.

\vspace{1ex}

\noindent [Gau77]\, P. Gauduchon --- {\it Le th\'eor\`eme de l'excentricit\'e nulle} --- C.R. Acad. Sc. Paris, S\'erie A, t. {\bf 285} (1977), 387-390.

\vspace{1ex}

\noindent [Ger66]\, M. Gerstenhaber --- {\it On the Deformation of Rings and Algebras: II} --- Ann. of Math., {\bf 84}, no. 1 (1966), 1–19.

\vspace{1ex}

\noindent [Gra58]\, H. Grauert --- {\it Analytische Faserungen \"uber holomorph-vollst\"andigen R\"aumen} --- Math. Ann. {\bf 135} (1958), 263–273.

\vspace{1ex}

\noindent [Hat17]\, A. Hatcher --- {\it Vector Bundles and K-Theory} --- \url{https://pi.math.cornell.edu/~hatcher/VBKT/VB.pdf}

\vspace{1ex}

\noindent [Kly89]\, A. A. Klyachko --- {\it Equivariant Bundles over Toric Varieties} --- Izv. Akad. Nauk SSSR Ser. Mat. {\bf 53} (5) (1989), 1001–1039, 1135 (Russian); English transl., Math. USSR-Izv. 35 (1990), no. 2, 337–375.

\vspace{1ex}

\noindent [Kob87]\, S. Kobayashi --- {\it Differential Geometry of Complex Vector Bundles} --- Princeton University Press, 1987.

\vspace{1ex}

\noindent [Kod86]\, K. Kodaira --- {\it Complex Manifolds and Deformations of Complex Structures} --- Grundlehren der Math. Wiss. {\bf 283}, Springer (1986).

\vspace{1ex}

\noindent [KS60]\, K. Kodaira, D.C. Spencer --- {\it On Deformations of Complex Analytic Structures, III. Stability Theorems for Complex Structures} --- Ann. of Math. {\bf 71}, no.1 (1960), 43-76.

\vspace{1ex}

\noindent [Lie78]\, D. Lieberman --- {\it Compactness of the Chow Scheme: Applications to Automorphisms and Deformations of K\"ahler Manifolds} --- Lect. Notes Math. {\bf 670} (1978), 140-186.

\vspace{1ex}

\noindent [Mas18]\, M. Maschio --- {\it On the Degeneration of the Fr\"olicher Spectral Sequence and Small Deformations} --- arXiv e-print math.DG 1811.12877v1.

%\vspace{1ex}

%\noindent [Mic83]\, M. L. Michelsohn -- {\it On the Existence of Special Metrics in Complex Geometry} -- Acta Math. {\bf 143} (1983) 261-295.

\vspace{1ex}

\noindent [Moi67]\, B.G. Moishezon --- {\it On $n$-dimensional Compact Varieties with $n$ Algebraically Independent Meromorphic Functions} --- Amer. Math. Soc. Translations {\bf 63} (1967) 51-177.

\vspace{1ex}

\noindent [MM90]\, R. R. Mazzeo, R. B. Melrose --- {\it The Adiabatic Limit, Hodge Cohomology and Leray's Spectral Sequence} --- J. Diff. Geom. {\bf 31} (1990) 185-213.

\vspace{1ex}

\noindent [Pop09]\, D. Popovici --- {\it Limits of Projective Manifolds Under Holomorphic Deformations} -- arXiv e-print AG 0910.2032v3.

\vspace{1ex}

\noindent [Pop10]\, D. Popovici --- {\it Limits of Moishezon Manifolds Under Holomorphic Deformations} -- arXiv e-print AG 1003.3605v1.

\vspace{1ex}

\noindent [Pop13]\, D. Popovici --- {\it Deformation Limits of Projective Manifolds: Hodge Numbers and Strongly Gauduchon Metrics} --- Invent. Math. {\bf 194} (2013), 515-534.

\vspace{1ex}

\noindent [Pop14]\, D. Popovici --- {\it Deformation Openness and Closedness of Various Classes of Compact Complex Manifolds; Examples} --- Ann. Sc. Norm. Super. Pisa Cl. Sci. (5), Vol. XIII (2014), 255-305.

 \vspace{1ex} 

\noindent [Pop16]\, D. Popovici --- {\it Degeneration at $E_2$ of Certain Spectral Sequences} ---  Internat. J. Math. {\bf 27}, no. 14 (2016), DOI: 10.1142/S0129167X16501111.

\vspace{1ex}

\noindent [Pop17]\, D. Popovici --- {\it Adiabatic Limit and the Fr\"olicher Spectral Sequence} --- Pacific Journal of Mathematics, Vol. 300, No. 1, 2019, dx.doi.org/10.2140/pjm.2019.300.121.

\vspace{1ex}

\noindent [PU14]\, D. Popovici, L. Ugarte --- {\it Compact Complex Manifolds with Small Gauduchon Cone} --- Proc. London Math. Soc. (3) {\bf 116} (2018), no. 5, 1161-1186.

\vspace{1ex}

\noindent [PSU20]\, D. Popovici, J. Stelzig, L. Ugarte --- {\it Higher-Page Bott-Chern and Aeppli Cohomologies and Applications} -- J. reine angew. Math. (Crelle), DOI 10.1515/crelle-2021-0014.

\vspace{1ex}

\noindent [Ste20]\, J. Stelzig --- {\it Toric Vector Bundles: GAGA and Hodge Theory} --- Math. Z. (2020).

\noindent https://doi.org/10.1007/s00209-020-02605-6

\vspace{1ex}

\noindent [Wit85]\, E. Witten --- {\it Global Gravitational Anomalies} --- Commun. Math. Phys, {\bf 100}, 197-229 (1985).

\vspace{6ex}

\noindent Universit\'e de Toulouse, Institut de Math\'ematiques de Toulouse

\noindent 118, route de Narbonne, 31062, Toulouse Cedex 9, France

\noindent Email: popovici@math.univ-toulouse.fr

\end{document}